\documentclass[11pt,a4paper]{amsart}
\usepackage[utf8]{inputenc}
\usepackage[T1]{fontenc}
\usepackage{amsmath,amssymb,amsthm,mathtools}
\usepackage[margin=2.7cm]{geometry}
\usepackage{booktabs,longtable,array}
\usepackage[hidelinks]{hyperref}

\newtheorem{theorem}{Theorem}[section]
\newtheorem{proposition}[theorem]{Proposition}
\newtheorem{lemma}[theorem]{Lemma}
\newtheorem{corollary}[theorem]{Corollary}
\newtheorem{conjecture}[theorem]{Conjecture}
\theoremstyle{definition}
\newtheorem{definition}[theorem]{Definition}
\theoremstyle{remark}
\newtheorem{remark}[theorem]{Remark}

\DeclareMathOperator{\res}{res}

\newcommand{\Id}{\mathrm{Id}}
\newcommand{\C}{\mathbb{C}}
\newcommand{\A}{\mathcal{A}}
\newcommand{\beq}{\begin{equation}}
\newcommand{\eeq}{\end{equation}}

\hypersetup{
  pdftitle={Superpotentials for F-manifolds and bi-flat F-structures on the orbit spaces of reflection groups},
  pdfauthor={Alessandro Arsie and Paolo Lorenzoni}
}
\title[Superpotentials and reflection groups]{%
  \texorpdfstring{Superpotentials for F-manifolds and bi-flat F-structures on the orbit spaces of reflection groups}%
  {Superpotentials for F-manifolds and bi-flat F-structures on the orbit spaces of reflection groups}}

\author{Alessandro Arsie}
\address{Department of Mathematics and Statistics, The University of Toledo, 2801 W. Bancroft St., Toledo, OH 43606, USA}
\email{alessandro.arsie@utoledo.edu}

\author{Paolo Lorenzoni}
\address{Dipartimento di Matematica e Applicazioni, Universit\`a degli Studi di Milano-Bicocca, Via Roberto Cozzi 53, I-20125 Milano, Italy, and INFN Sezione di Milano-Bicocca}
\email{paolo.lorenzoni@unimib.it}


\date{}

\begin{document}
\raggedbottom
\begin{abstract}
{
We extend Dubrovin's construction of a superpotential from Frobenius manifolds to bi-flat F-manifolds.  The key step in this extension is choosing a flat coordinate of one of the connections studied in \cite{AL} and
invert this coordinate with respect to the spectral parameter and the construction does not require
an invariant metric.
Suppose that the manifold has dimension $n$ and that the superpotential admits a holomorphic
continuation to $n$ simple critical points with distinct critical values. These values give canonical
coordinates, and the residue formulas recover both the product and the flat connection from
the superpotential alone. Conversely, under suitable nondegeneracy conditions, a solution of the
superpotential equations for a given product defines a bi-flat F-manifold. For $B_n$ and $G(m,1,n)$ we obtain
the whole family of bi-flat structures from a polynomial superpotential by varying a twisted primary
differential. For $G(m,m,n)$ we obtain the natural bi-flat structure from a Laurent superpotential. Keeping
the Laurent polynomial fixed and varying the primary differential gives a new family of flat F-manifolds depending on one parameter.

We use the superpotential and the poles of the dual product along the reflecting hyperplanes to
prove the conjecture we formulated in \cite{ALrefl} for irreducible Coxeter groups and in
\cite{ALsurvey} for well generated complex reflection groups. The conjecture concerns bi-flat F-structures whose dual connection has the form studied by Dunkl and Kohno. When the reflecting
hyperplanes form a single orbit, we identify the bi-flat structure associated with the flat structure
of \cite{KMS1,KMS2,KMS}. We obtain the remaining families from the polynomial construction for $G(m,1,n)$
and an explicit computation for $F_4$. In every case the first flat structure extends polynomially to
the whole quotient, and the bi-flat structure is defined on the complement of the reflection
discriminant.}
\end{abstract}
\maketitle

\tableofcontents

\section{Introduction}

{
Let $G$ be a finite irreducible Coxeter group and let $V$ be the complexification of its reflection
representation. In \cite{Dub,DCoxeter}, Dubrovin used Saito's flat pencil of metrics \cite{Sa,SYS}
to construct a Frobenius manifold structure on $V/G$. The flat coordinates of Saito's flat metric give
distinguished basic invariants in which the coefficients of the product are polynomial. By almost
duality \cite{Dalm}, one obtains a second structure on the regular part of the quotient. Indeed, the
Levi-Civita connection of the intersection form pulls back to the trivial connection on the
complement of the reflecting hyperplanes, while the dual product has simple poles along these
hyperplanes and comes from a $\vee$ system in the sense of Veselov \cite{Ve,Ve2}.

For a complex reflection group, the natural Saito structure need not admit a compatible Frobenius
metric. This leads us to consider flat F-manifolds. Such a manifold carries a commutative
associative product whose unit is parallel with respect to a flat torsion free connection. In flat
coordinates the coefficients of the product are the second derivatives of the components of a local
vector potential.

In the case of a \emph{bi-flat F-manifold} \cite{ALDE,ALsurvey}, there are two flat structures $(\nabla,\circ,e)$ and
$(\nabla^*,*,E)$, related by the Euler field $E$. This notion was introduced to describe the pair formed by a
Frobenius manifold and its almost dual without introducing a metric. In \cite{ALrefl} we used it to
study the orbit spaces of complex reflection groups. Kato, Mano, and Sekiguchi \cite{KMS1,KMS2}
constructed a flat F-manifold, or equivalently a Saito structure without metric, with a polynomial
potential on the orbit space of every nontrivial irreducible well generated complex reflection
group. By developing an almost duality for Saito structures, Konishi, Minabe, and Shiraishi
\cite{KMS} obtained this natural structure from the trivial connection on the reflection
representation.

In \cite{ALrefl} we asked whether one could express the dual connection directly in terms of the
reflecting hyperplanes. For groups of ranks $2$ and $3$ we found bi-flat structures on the regular
quotient whose dual connection has the form studied by Dunkl and Kohno,
\[
\nabla^*=\nabla^0-\sum_{H}\frac{d\alpha_H}{\alpha_H}\otimes\tau_H\,\pi_H ,
\]
where $\alpha_H$ is a defining form for the hyperplane $H$ and $\pi_H$ is the projection onto its normal line. The
weights $\tau_H$ are constant on each orbit of reflecting hyperplanes. Moreover, the condition $\nabla^*E=0$ gives
one linear relation among them, so there are $k-1$ parameters when the hyperplanes form $k$ orbits.

This observation led us to formulate a conjecture for Coxeter groups in \cite{ALrefl} and for well
generated complex reflection groups in \cite{ALsurvey}. Let $V$ be the reflection representation and
denote by $V^{\mathrm{reg}}$ the complement of the reflecting hyperplanes. For every nontrivial irreducible group
that is well generated, the conjecture asks for a family of bi-flat structures on $V^{\mathrm{reg}}/G$ depending on
$|\mathcal H/G|-1$ parameters, where $\mathcal H/G$ denotes the set of orbits of reflecting hyperplanes. This is the largest
possible number of parameters for such structures \cite{ALsurvey}. In this family the dual
connection has the form given above and the natural flat structure extends polynomially to $V/G$. We
give the precise statement in Conjecture~\ref{conj:AL}.

The conjecture was already known for many groups of ranks $2$ and $3$, as well as for $B_4$ and $F_4$. In
this paper we construct the families for $B_n$ and $G(m,1,n)$ in arbitrary rank and give explicit formulas for
$F_4$. Our argument for groups with a single orbit applies in particular to $G(m,m,n)$ for $m,n\ge3$ and to $G_{29},G_{32},G_{33},G_{34}$.
Together with the cases already known, these results prove the following theorem.

\begin{theorem}\label{thm:intro}
Conjecture~\ref{conj:AL} holds for every nontrivial well generated irreducible complex reflection group.
\end{theorem}

To see how the number of orbits enters the proof, recall that the reflecting hyperplanes of each
group under consideration form either one or two orbits. In Theorem~\ref{thm:naturalDunklKohno} we
express the dual of the natural structure of \cite{KMS1,KMS2,KMS} in terms of these hyperplanes for
every such group. Indeed, in linear coordinates on $V$ the coefficients of the dual product are
homogeneous rational functions of degree $-1$. The Jacobian of the quotient map gives a bound on the
order of their poles along each reflecting hyperplane. The reflections fixing the hyperplane
restrict the terms that can occur in the polar part. Compatibility with the trivial connection then
leaves only a simple pole with constant residue. Once we subtract the polar parts, the remaining
coefficients are holomorphic and homogeneous of degree $-1$, so they vanish.

If there is a single orbit, the normalization determines the weight in this formula. If there are
two orbits, one parameter remains and the natural structure is a member of the resulting family. For
$F_4$ we compute this family explicitly in Proposition~\ref{prop:F4}. For $B_n$ and $G(m,1,n)$ we obtain it from a
superpotential.

Recall that Dubrovin associates a superpotential $\lambda(p,u)$ to a semisimple Frobenius manifold. Its critical
values as a function of $p$ give the canonical coordinates, and residue formulas recover the metric
and the product. Dubrovin's construction uses the Frobenius metric. Extensions of the Hurwitz
construction to flat F-manifolds were studied in \cite{Alm,PS}. In Remark~3.6 of \cite{PS},
Proserpio and Strachan suggest using the connections of Gauss and Manin type studied in \cite{AL} to
go beyond the Hurwitz case. The question of defining an analog of Dubrovin's superpotential for F-manifolds is also a thorny open issue arising in the study of $F$ topological recursion in Section~7 of \cite{BGU}.

Let us describe the construction. Fix a weight $\nu_0$ in the family of connections $\Gamma(\lambda,\nu_0)$ of \cite{AL}. By
Theorem~\ref{thm:ext}, the flow $\phi$ of the unit $e$ carries a flat function $P_0$ of the deformed dual
connection $\nabla^*+\nu_0\,*$ to a function
$P(x,\lambda)=P_0(\phi_{-\lambda}(x))$
which is flat on each regular slice where the formula is defined. Wherever $\partial_\lambda P\ne0$, we keep $x$ fixed and
invert $P(x,\lambda)$ with respect to $\lambda$. We call the resulting local inverse $\Lambda(P,x)$ a superpotential
(Definition~\ref{def:sp}). Proposition~\ref{prop:ext} shows that the functions obtained by
translation are flat for a single connection on the open subset of $M\times\mathbb C$ where $E\circ-\lambda\Id$ is invertible.

In Theorem~\ref{thm:char} we characterize this inverse by the unit condition $\partial_e\Lambda=1$, a homogeneity
condition, and the identity
\[
\Lambda_{hk}=\partial_P\left[\frac{\Lambda_h\Lambda_k-c^l_{hk}\Lambda_l}{\Lambda_P}\right].
\]
This identity is precisely Alcolado's criterion for extending the structure by one coordinate
\cite{Alc}. It thus gives the relation with the open WDVV equations discussed in
Remark~\ref{rem:PS}. If the superpotential continues holomorphically to $n$ simple critical points
with distinct critical values, these values give canonical coordinates by Theorem~\ref{thm:main}.

We next ask how to recover the connection once the product and the superpotential are given
(Section~\ref{subsec:reconstruction}). Evaluating the compatibility equation at $n$ simple critical
points with distinct critical values determines the Christoffel symbols of $\nabla$ and gives a residue
formula for them (Proposition~\ref{prop:christoffel}). In canonical coordinates the compatibility
equation reduces to $\binom n2$ scalar equations (Proposition~\ref{prop:S3can}). Moreover, the bilinear form
in the residue formula is one of the metrics of the Riemannian F-manifold with Killing unit
associated with the flat structure (Proposition~\ref{prop:riem}). Its Lam\'e coefficients are
therefore homogeneous of the same degree $\nu_0-\frac12$. If the degrees differ, a single superpotential cannot
realize all $n$ canonical coordinates as simple critical values (Corollary~\ref{cor:lame}). If the
degrees agree, a suitable choice of period realizes all of them (Proposition~\ref{prop:lamecrit}).
The critical values also give information about the product. Indeed, their differentials are algebra
homomorphisms and determine the product whenever they span the cotangent space.

Conversely, suppose we are given a commutative associative product with unit and a vector field $E$.
We seek a connection for this product. Assume that the superpotential satisfies the unit and
homogeneity conditions and that the compatibility equation admits a solution for the Christoffel
symbols. If $n$ local inverse branches have independent differentials with respect to the base
variables, this solution defines a unique flat torsion free connection $\nabla$ (Theorem~\ref{thm:flat}).
If, in addition, $N=\nu_0\Id-\nabla E$ is invertible, then $\nabla$ is compatible with the product, the unit is parallel, and
$E$ is a linear Euler field. By Proposition~\ref{prop:dualflat}, the resulting flat F-structure gives
a bi-flat pair wherever $E\circ$ is invertible. In Corollary~\ref{cor:flat} we recover flat functions of
the deformed dual connection on the regular part of the zero slice by considering the simple zeros
of the continued superpotential.

Table~\ref{tab:reconstruction} at the end of Section~\ref{sec:def} records the hypotheses under
which each part of the structure can be recovered locally.

Wherever $E\circ$ is invertible, the natural structure together with $E$ also determines the dual product
and the dual connection (Proposition~\ref{prop:dualflat}). Starting from a given bi-flat structure,
one still has to choose a superpotential. Indeed, at a fixed weight the local periods, up to the
affine changes $P\mapsto aP+b$ with $a\ne0$, correspond to lines of flat covectors $\xi$ at a regular reference point
satisfying $\xi(e)\ne0$ (Remark~\ref{rem:norm}).

If we start with a Frobenius manifold of charge $d$ and choose the weight $\nu_0=(1-d)/2$, the local inversion
described above agrees with Dubrovin's inversion of periods (Proposition~\ref{prop:Dub}).

Let us return to the reflection groups. For $G(m,1,n)$, including the case of $B_n$, we take as superpotential
the monic polynomial
\[
\lambda(z)=z^n+a_1z^{n-1}+\dots+a_n,\qquad z=p^m,
\]
and as twisted primary form the differential $z^{\beta-1}dz$. We choose a primitive that depends holomorphically
on $\beta$ and thus includes the logarithmic case $\beta=0$. The quotient $\mathbb C[z]/(z\lambda')$ defines the product, while the
weights of the dual connection depend on $\beta$. Taking twisted residues of powers of $\lambda$ gives flat
coordinates for $\nabla$. In this way we recover the families found in low rank in \cite{ALrefl,Bn} and
obtain the construction in every rank (Theorem~\ref{thm:Bn} and Proposition~\ref{prop:Gm1n}).

For $G(m,m,n)$ with $m\ge2$ and $n\ge3$, we consider the Laurent polynomial
\[
\lambda(z)=z^{n-1}+a_1z^{n-2}+\dots+a_{n-1}+\frac{a_n}z,\qquad a_n=(-1)^nt^m,\qquad t=p_1\cdots p_n.
\]
In this case there are two poles at which to take the twisted residues. By varying the parameter in
the primary form we obtain a new family of flat F-manifolds (Corollary~\ref{thm:Laurent}). On the
regular part of the Laurent model the dual product comes from $G(m,m,n)$, while the dual connection has the
form studied by Dunkl and Kohno for the arrangement of $G(m,1,n)$ (Remark~\ref{rem:Gmmnbeta}). We compute
the vector potentials for $2\le n\le4$.

The natural flat structure extends across $t=0$ to the whole quotient of $G(m,m,n)$ precisely when $\beta=1/m$. For this
value of the parameter the dual structure also extends across the regular points on $t=0$. Once we fix
the unit and the Euler field as in Theorem~\ref{thm:Gmmn}, this is the unique pair for which both
the dual product and the dual connection come from the arrangement of $G(m,m,n)$. The Laurent polynomial
gives a distinguished superpotential for this pair, namely $\lambda(p^m)$ with primary form proportional to
$dp$. When $m=2$ we recover Dubrovin's superpotential of type $D_n$. The Laurent model thus gives a
second proof of the conjecture for these groups.

In Section~\ref{sec:prelim} we recall the flat structures on orbit spaces and state the conjecture.
Table~\ref{tab:coverage} lists the cases that enter its proof. In Section~\ref{sec:def} we construct
the superpotential, starting with the family of connections of Gauss and Manin in
Section~\ref{sec:GM} and the translation along the unit in Section~\ref{sec:GMflat}. We define the
superpotential in Section~\ref{sec:defsp}, give an equivalent description in terms of the extended
connection in Section~\ref{sec:ext}, and characterize it by differential equations in
Section~\ref{sec:char}. In Section~\ref{subsec:reconstruction} we prove the reconstruction results,
and in Section~\ref{sec:Dub} we compare the construction with Dubrovin's. Section~\ref{sec:two}
deals with the families for $B_n$ and $G(m,1,n)$ and gives the formulas for $F_4$. In
Section~\ref{sec:ALconjectures} we consider the case of a single orbit. We first identify the dual
of the natural structure and then construct the Laurent family for $G(m,m,n)$.

\noindent\textbf{Acknowledgments.} The authors thank GNFM (Gruppo Nazionale di Fisica Matematica) for supporting activities that contributed to the research reported in this paper. P. Lorenzoni is supported by funds of INFN (Istituto Nazionale di Fisica Nucleare) under IS-CSN4 Mathematical Methods of Nonlinear Physics.

\noindent\textbf{AI Statement.} Parts of this work were developed in interaction with the AI systems Claude (Anthropic) and ChatGPT (OpenAI). They were used to explore possible definitions, to draft some proofs and to check computations symbolically and with exact arithmetic. The authors have checked all statements and proofs and take full responsibility for the content.
}

\section{Preliminary material}\label{sec:prelim}

{
We first recall the two flat structures that occur in the conjecture.
In Sections~\ref{sec:flatbi-flat} and~\ref{sec:almostdual} we describe the conditions that define them and explain
how their connections are related. In Section~\ref{sec:orbitspaces} we turn to the orbit spaces
of reflection groups. As we shall see, the number of orbits of reflecting
hyperplanes determines the number of parameters in the conjecture of
Section~\ref{sec:conjecture}.
}

\subsection{\texorpdfstring{Flat and bi-flat structures}{Flat and bi-flat structures}}\label{sec:flatbi-flat}

\begin{definition}

Let $M$ be a manifold endowed with a bilinear product $\circ:TM\times_M TM\to TM$,
a connection $\nabla$, and a vector field $e$. We denote by $c^i_{jk}$ the structure
constants of $\circ$ and by $\Gamma^i_{jk}$ the Christoffel symbols of $\nabla$. We say that
$(M,\circ,\nabla,e)$ is a flat F-manifold if the following conditions hold.
\begin{enumerate}
\item For each $\lambda$, the connection $\nabla_{(\lambda)}$ with
Christoffel symbols
\[\Gamma^i_{jk}+\lambda c^i_{jk}\]
is flat and torsion free.
\item The vector field $e$ is the unit of $\circ$.
\item The unit is parallel with respect to $\nabla$, that is, $\nabla e=0$.
\end{enumerate}
\end{definition}

{
If the first two conditions hold, the same data define an
\emph{almost flat F-manifold}. Let us examine the conditions that the
first axiom imposes on $\nabla$ and $\circ$. The torsion of $\nabla_{(\lambda)}$ is given by
\[
T^{(\lambda)k}_{ij}=\Gamma^k_{ij}-\Gamma^k_{ji}
+\lambda(c^k_{ij}-c^k_{ji}).
\]
The coefficients of $\lambda^0$ and $\lambda^1$ vanish precisely when $\nabla$ is torsion free
and $\circ$ is commutative. In this case the curvature of the pencil takes the form
\[
R^{(\lambda)k}_{ijl}=R^k_{ijl}+\lambda(\nabla_i c^k_{jl}-\nabla_j c^k_{il})
+\lambda^2(c^k_{im}c^m_{jl}-c^k_{jm}c^m_{il}),
\]
where $R$ denotes the curvature tensor of $\nabla$. The constant term vanishes
precisely when $\nabla$ is flat. The vanishing of the linear term is equivalent
to the symmetry of $\nabla_l c^k_{ij}$ in its lower indices, and the vanishing of the
quadratic term is equivalent to the associativity of $\circ$. Together with the
conditions obtained from the torsion, these are precisely the requirements
of the first axiom.

In flat coordinates for $\nabla$, the symmetry in the lower indices implies
the local existence of a vector potential $A=(A^i)$ such that
\[c^i_{jk}=\partial_j\partial_k A^i.\]
In terms of $A$, the condition of associativity becomes the system of
oriented associativity equations, also known as the generalized WDVV equations,
\[
\partial_j\partial_l A^i\,\partial_k\partial_m A^l
=\partial_k\partial_l A^i\,\partial_j\partial_m A^l .
\]
In the definition of a Frobenius manifold given by Dubrovin, the flat F-manifold is also endowed with a nondegenerate symmetric metric $\eta$ and an
Euler field $E$. The metric is parallel and invariant under the product, that is,
\[
\nabla\eta=0,\qquad \eta_{il}c^l_{jk}=\eta_{jl}c^l_{ik}.
\]
It follows that we can choose the local vector potential so that
$A^i=\eta^{il}\partial_l F$ for a scalar potential $F$. With this choice, the oriented associativity
equations reduce to the WDVV equations
\[
(\partial_a\partial_j\partial_l F)\eta^{ls}(\partial_s\partial_k\partial_m F)
=(\partial_a\partial_k\partial_l F)\eta^{ls}(\partial_s\partial_j\partial_m F).
\]
The action of the Euler field $E$ on the product and the metric is given by
\[\mathrm{Lie}_E\circ=\circ,\qquad {\rm Lie}_E\eta=(2-d)\,\eta\]
for some constant $d$. Since $\eta$ is parallel, the second identity implies
that $E$ is linear, in the sense that $\nabla\nabla E=0$.
}

\subsection{Dual and almost dual structures}\label{sec:almostdual}

{
Let $(M,\eta,\circ,e,E)$ be a Frobenius manifold. Multiplication by $E$ defines the
contravariant metric $g=(E\circ)\,\eta^{-1}$, which is called the \emph{intersection form}.
On the open set where $E\circ$ is invertible, we denote its Levi-Civita connection
by $\tilde\nabla$ and define the dual product by
\[
X*Y=(E\circ)^{-1}\,X\circ Y .
\]
Together with the unit $E$, this connection and product define the almost
dual structure, namely an almost flat F-structure with invariant metric $g^{-1}$.

In general $\tilde\nabla E\ne0$. However, we can make the dual unit parallel by choosing
a constant $\bar\lambda$ and setting $\nabla^*=\tilde\nabla+\bar\lambda*$. In this way we obtain a flat structure
$(\nabla^*,*,E)$ on the open set where $E\circ$ is invertible. Its connection and the connection
of the original flat structure $(\nabla,\circ,e)$ satisfy the compatibility conditions
\[(d_{\nabla}-d_{\nabla^{*}})(X\,\circ)=0,\quad(d_{\nabla}-d_{\nabla^{*}})(X\,*)=0\quad\text{for any vector field}\,X.\]
Here $d_{\nabla}$ denotes the exterior covariant derivative. In local coordinates it is given by
\[(d_{\nabla} \omega)^l_{i_0\dots i_k}=\sum_{j=0}^k (-1)^j \nabla_{i_j}\omega^l_{i_0\dots\hat{i}_j\dots i_k}.\]
Only the products and the connections enter these compatibility conditions.
This observation leads to the following definition.
}

\begin{definition}

A \emph{bi-flat F-manifold} is a manifold endowed with two flat structures
$(\nabla,\circ,e)$ and $(\nabla^*,*,E)$ on a common domain, where $E$ is a linear Euler field and $*$ is the
dual product induced by $E$. Their connections must satisfy $(d_\nabla-d_{\nabla^*})(X\circ)=0$ for every
vector field $X$ or, equivalently, $(d_\nabla-d_{\nabla^*})(X*)=0$. An equivalent expression for this
condition is
\begin{equation}\label{eq:dual}
\nabla^*_XY=\nabla_XY-\nabla_{X*Y}E.
\end{equation}
\end{definition}

{
Recall that the linearity of $E$ means that $\nabla\nabla E=0$, while the Euler conditions read
$[e,E]=e$ and $\mathrm{Lie}_E\circ=\circ$. The dual product $X*Y=(E\circ)^{-1}X\circ Y$ is defined wherever $E\circ$ is invertible, and its
unit is parallel precisely when $\nabla^*E=0$. We shall use \eqref{eq:dual} both to obtain the dual
connection from the natural one and to recover the natural structure from
the dual one. Let us therefore explain how this identity follows from
compatibility.
}

\begin{lemma}\label{lem:AHE}

Let $\nabla$ and $\nabla^*$ be torsion free connections, and let $\circ$ be a commutative,
associative product with unit $e$. Assume that $E$ is a vector field such that
$E\circ$ is invertible and $\nabla^*E=0$. Define $X*Y=(E\circ)^{-1}X\circ Y$. Then the following conditions are equivalent.
\begin{enumerate}
\item[(a)] $(d_\nabla-d_{\nabla^*})(X\circ)=0$ for every vector field $X$.
\item[(b)] $(d_\nabla-d_{\nabla^*})(X*)=0$ for every vector field $X$.
\item[(c)] $\nabla^*_XY=\nabla_XY-\nabla_{X*Y}E$ for all vector fields $X,Y$.
\end{enumerate}
\end{lemma}

\begin{proof}

To see that (a) and (b) are equivalent, observe that $X*Y=(E^{-1}\circ X)\circ Y$. Since multiplication
by $E^{-1}$ is invertible, the endomorphisms $X*$ and $X\circ$ describe the same family as $X$ varies.

Set $D(Y,Z)=\nabla^*_YZ-\nabla_YZ$. Since both connections are torsion free, $D$ is symmetric. For any
endomorphism $A$ we have $(\nabla^*_YA)Z=(\nabla_YA)Z+D(Y,AZ)-A\,D(Y,Z)$. Alternating in $Y$ and $Z$, we see that the terms
containing $A\,D$ cancel. This gives
\[
(d_{\nabla^*}-d_\nabla)(A)(Y,Z)=D(Y,AZ)-D(Z,AY).
\]
For $A=X*$, condition (b) therefore implies that $D(Y,X*Z)=D(Z,X*Y)$ for all $X,Y,Z$. In particular,
choosing $Z=E$, we obtain $D(Y,X)=D(E,X*Y)$. Since both connections are torsion free and $\nabla^*E=0$, we have
\[
D(E,W)=\nabla^*_EW-\nabla_EW=\big(\nabla^*_WE+[E,W]\big)-\big(\nabla_WE+[E,W]\big)=-\nabla_WE .
\]
It follows that $D(Y,X)=-\nabla_{X*Y}E$, which proves (c).

Conversely, suppose that (c) holds. Applying the same identity with $A=X\circ$, we obtain
\[
D(Y,X\circ Z)-D(Z,X\circ Y)=-\nabla_{E^{-1}\circ Y\circ X\circ Z}E+\nabla_{E^{-1}\circ Z\circ X\circ Y}E=0
\]
by the commutativity and associativity of $\circ$. This proves (a).
\end{proof}

Starting from the natural structure, one can determine the dual structure as follows.
This result was proved in \cite{AL} in the semisimple case and in \cite{KMS} in the general case.

\begin{proposition}\label{prop:dualflat}
Let $(M,\nabla,\circ,e)$ be a flat F-manifold with a linear Euler field
$E$. On the open set where $E\circ$ is invertible, define the product
$*$ by $X*Y=(E\circ)^{-1}(X\circ Y)$ and the connection $\nabla^*$ by
\eqref{eq:dual}. Then $(\nabla^*,*,E)$ is a flat F-structure, and
$(\nabla,\circ,e)$ and $(\nabla^*,*,E)$ form a bi-flat pair.
\end{proposition}

Thus, to specify a bi-flat pair, it suffices to give the natural flat
structure and a linear Euler field for which $E\circ$ is invertible.

Conversely, suppose that we are given $(\nabla^*,*,E)$ and that the natural unit $e$ is fixed.
To recover the natural structure, observe first that multiplication by $e$ with
respect to $*$ is $(E\circ)^{-1}$. Moreover, it follows from \eqref{eq:dual} and $\nabla e=0$ that $\nabla^*_{X\circ Y}e=-\nabla_{X*Y}E$, since $(X\circ Y)*e=X*Y$.
We therefore recover the natural product and connection from the formulas
\[
X\circ Y=(e*)^{-1}(X*Y),\qquad
\nabla_XY=\nabla^*_XY-\nabla^*_{X\circ Y}e .
\]
The choice of $e$ matters in the statement of uniqueness. Indeed, for any
nonzero constant $a$, the product $a\circ$ with unit $a^{-1}e$ gives the same dual data $(\nabla^*,*,E)$.

For later use, choose flat coordinates $u$ for $\nabla$ and denote the structure
constants of $\circ$ by $c^k_{ij}$. If $C_h$ denotes multiplication by $\partial_h$, then $(C_h)^l_k=c^l_{hk}$.
In these coordinates we have
\[
Q:=\nabla E\ (\text{constant}),\qquad L:=E\circ,\qquad L_\lambda:=L-\lambda\,\Id,\qquad
c^{*j}_{hk}=(L^{-1})^j_m\,c^m_{hk}.
\]
In these coordinates the Christoffel symbols of $\nabla$ satisfy $a^j_{hk}=0$. It follows
from \eqref{eq:dual} that the Christoffel symbols of the dual connection are $b^j_{hk}=-Q^j_p\,c^{*p}_{hk}$.

\subsection{\texorpdfstring{Flat and bi-flat structures on orbit spaces of reflection groups}{Flat and bi-flat structures on orbit spaces of reflection groups}}\label{sec:orbitspaces}

\subsubsection{Coxeter groups}

{
For a finite irreducible Coxeter group, Dubrovin gave a description of
the flat coordinates and the Euler field in terms of the basic invariants \cite{DCoxeter}.
}

\begin{theorem}\label{Dcoxeterth}

On the quotient of the complexified reflection representation of a finite
irreducible Coxeter group there is a Frobenius manifold structure $(\eta,\circ,e,E)$ with
the following properties.
\begin{enumerate}
\item The flat coordinates of $\eta$ form a system of basic invariants $(u_1,\dots,u_n)$ of
the group. They are called \emph{Saito flat coordinates}.
\item In these coordinates the unit and the Euler field are given by
\[e=\frac{\partial}{\partial u_n},\qquad E=\sum_{i=1}^n\left(\frac{d_i}{d_n}\right)u_i\frac{\partial}{\partial u_i},\]
\end{enumerate}
where $d_i$ denotes the degree of the invariant polynomial $u_i$. For rank $n\ge2$,
we have $2=d_1<d_2\le \dots\le d_n$, and the largest degree $d_n$ occurs only once. In rank one we have
$d_1=2$. In either case $d_n$ is the Coxeter number.
\end{theorem}

{
Dubrovin starts with a flat pencil of contravariant metrics on the
quotient. The pencil is generated by the Saito metric and the form induced
by the Euclidean metric on the reflection representation \cite{Sa,SYS}. In the linear
coordinates of the representation, we can express the dual product in terms
of the reflecting hyperplanes and the projections onto their normal lines.
}

\begin{theorem}\label{almostdualCoxeter}

Let $G$ be a finite Coxeter group acting on the Euclidean space $\mathbb E^n$ with
coordinates $(p^1,\dots,p^n)$. On the complexification of this space, denote by $g$ the
complex bilinear extension of the Euclidean metric and by $\widetilde\nabla$ its Levi-Civita
connection. On the complement of the reflecting hyperplanes an almost flat
structure with invariant metric $g$ is given by
\[\left(\widetilde\nabla,\quad*=\sum_{H\in \mathcal{H}}\frac{d\alpha_H}{\alpha_H}\otimes\sigma_H\pi_H,\quad E=\sum_k p^k\frac{\partial}{\partial p^k}\right).\]
Here we use the following notation.
\begin{itemize}
\item $\mathcal{H}$ is the set of reflecting hyperplanes of $G$,
\item $\alpha_H$ is a linear form defining the hyperplane $H$,
\item $\pi_H$ is the complex linear extension of the orthogonal projection onto $H^\perp$,
\item the weights $\sigma_H$ are invariant under $G$ and satisfy
\[\sum_{H\in \mathcal{H}}\sigma_H\pi_H=\mathrm{Id},\]
\end{itemize}
The data descend to the regular quotient.
\end{theorem}

\begin{proof}

By the normalization, the common fixed subspace is zero. If the reflection
representation is reducible, we decompose it as an orthogonal sum of
irreducible factors. The data decompose in the same way. On each factor
$\widetilde\nabla+s*$ is a Dunkl connection with invariant weights $-s\sigma_H$, so it is flat by \cite{CHL}.
In dimension one flatness follows directly. Each summand in the product $*$
is symmetric in $X$ and $Y$. Thus all the connections in the pencil are torsion
free. Evaluating the product at $E$ and using the normalization, we find $E*X=X$.
To check invariance of the metric, fix a hyperplane $H$ and pair the corresponding
summand of $X*Y$ with $Z$ using $g$. The result is a scalar multiple of
$\alpha_H(X)\alpha_H(Y)\alpha_H(Z)$, which is symmetric in the three arguments. This proves invariance of $g$.
Finally, all the data are invariant under $G$ and therefore descend to the
regular quotient. For the relation with $\vee$ systems, see \cite{ALjmp,FV,Ve,Ve2}.
\end{proof}

{
In the irreducible Coxeter case, Dubrovin's almost dual structure corresponds
to the choice of equal weights $\sigma_H$ in Theorem~\ref{almostdualCoxeter}. More precisely, if $h=d_n$, the
product in that theorem is $1/h$ times Dubrovin's almost dual product. The unit
is therefore $hE$, where $E$ is the Euler field of Theorem~\ref{Dcoxeterth}. The unshifted
connection is the same in the two descriptions \cite{Dalm}.
}

\subsubsection{Complex reflection groups}

{
Let $V=\C^n$ and let $W\subset GL(V)$ be a finite group. By a \emph{(complex) reflection} we
mean an element $s\in GL(V)$ of finite order whose fixed space $\ker(s-1)$ is a hyperplane $H_s$.
We call $W$ a \emph{complex reflection group} if it is generated by reflections.
We assume that $W$ is nontrivial and \emph{irreducible}, that is, the only
complex subspaces of $V$ that are invariant under $W$ are $0$ and $V$.
The dimension $n=\dim V$ is called the rank of the representation.
Since several reflections may fix the same hyperplane, we distinguish the arrangement
\[
  \A=\A(W)=\{H_s : s\in W \text{ a reflection}\}
\]
of reflecting hyperplanes, of cardinality $N^*=|\A|$, from the set of reflections,
whose cardinality we denote by $N$. The group $W$ acts on $\A$ by $w\cdot H=w(H)$.
The orbits of this action index the weights in the dual connection.

The theorem of Chevalley, Shephard, and Todd \cite{Ch,ST} states that $\C[V]^W$ is a polynomial
algebra generated by homogeneous invariants of degrees $d_1\le\dots\le d_n$. We order the
corresponding \emph{codegrees} so that $d_1^*\ge\dots\ge d_n^*$ \cite{OS,Br}. In this notation we have
\[
  |W|=\prod_i d_i,\qquad N=\sum_i (d_i-1),\qquad N^*=\sum_i (d_i^*+1)
\]
as recalled in \cite{LT} and \cite{OT}. For the groups that occur in our main theorem
there is a further relation between the degrees and codegrees. To describe
it, we first recall a condition on the generating reflections.
}

\begin{definition}

A nontrivial irreducible complex reflection group $W$ of rank $n$ is
\emph{well generated} if it can be generated by $n$ reflections.
\end{definition}

{
For any set of reflections that generates the group, the intersection
of their fixed hyperplanes is $V^W=0$. Such a set must therefore contain at least
$n$ reflections. For each group that is not well generated, the classification
provides a generating set of $n+1$ reflections \cite{LT}.

Orlik and Solomon \cite{OS} and Bessis \cite{Be2} give the following characterization in
terms of the degrees and codegrees.
}

\begin{theorem}

Let $W$ be nontrivial and irreducible. Then $W$ is well generated if and only if
\[
  d_i+d_i^* = d_n \qquad \text{for all } i=1,\dots,n .
\]
In this case $h:=d_n$ is the \emph{Coxeter number}, and $h=(N+N^*)/n$.
\end{theorem}

{
Every finite real reflection group, or Coxeter group, is well generated \cite{Hu}.
The class of well generated groups provides a natural setting for Coxeter
elements, for the lattice of noncrossing partitions $NC(W)$ \cite{Be1,Be2}, for the covering
of Lyashko and Looijenga \cite{Be2}, and for the Catalan numbers of $W$, given by
$\prod_i (d_i+h)/d_i$ \cite{Be1,Be2,Ar}.

According to the classification of Shephard and Todd \cite{ST}, every nontrivial
irreducible complex reflection group is a nontrivial cyclic group in rank
one, an imprimitive group
$G(m,p,n)$ in rank at least two, the symmetric group $\mathfrak S_n$ in its reflection representation,
or one of the $34$ exceptional groups $G_4,\dots,G_{37}$.
For $p\mid m$, the group $G(m,p,n)$ consists of the $n\times n$ monomial matrices for which each
nonzero entry is an $m$th root of unity and the product of these entries is
an $(m/p)$th root of unity. The group has order $m^n n!/p$. When $m,n\ge2$, its representation is
irreducible except in the case of $G(2,2,2)$. When $m=p=1$ and $n\ge2$, the monomial
representation of $\mathfrak S_n$ restricts to the hyperplane $\sum_i x_i=0$, where it gives the
irreducible reflection representation. If $n=1$, then $G(m,p,1)$ is cyclic of order $m/p$.

Applying the characterization in terms of degrees and codegrees to this
classification gives the following list \cite{ST,Be1}. See also \cite{LT}.
}

\begin{theorem}

A nontrivial irreducible complex reflection group is well generated if and
only if it occurs in the following list.
\begin{enumerate}
  \item $G(m,1,n)$ with $m\ge2$ and $n\ge1$. For $n=1$ this is the cyclic group $\mathbb Z/m$.
  \item $G(m,m,n)$ with $m\ge2$ and $n\ge2$, excluding the reducible group $G(2,2,2)$,
  and the symmetric group $G(1,1,n)\cong\mathfrak S_n$ in its reflection representation of dimension $n-1$.
  \item In rank $2$, the exceptional groups
        $G_4,G_5,G_6,G_8,G_9,G_{10},G_{14},G_{16},G_{17},G_{18},G_{20},G_{21}$.
  \item In ranks at least $3$, the exceptional groups
        $G_{23},\dots,G_{30},G_{32},\dots,G_{37}$.
\end{enumerate}
The remaining groups, which are not well generated, are $G(m,p,n)$ with $1<p<m$ and $n\ge2$,
the exceptional groups of rank $2$ given by $G_7,G_{11},G_{12},G_{13},G_{15},G_{19},G_{22}$, and $G_{31}$.
\end{theorem}

{
This list includes the complexifications of several real reflection groups.
}

\begin{remark}

In the notation above, the real reflection groups are given by $A_{n-1}=G(1,1,n)$, $B_n=G(2,1,n)$, $D_n=G(2,2,n)$,
$I_2(m)=G(m,m,2)$, $H_3=G_{23}$, $F_4=G_{28}$, $H_4=G_{30}$, and $E_6,E_7,E_8=G_{35},G_{36},G_{37}$. For the real case we refer to \cite{Hu}. Notice, for instance,
that $G(2,2,3)\cong G(1,1,4)$ and that $G(4,4,2)$ is conjugate to $G(2,1,2)$.
\end{remark}

{
Within the imprimitive families $G(m,p,n)$, three orbits of reflecting hyperplanes
can occur only for $G(m,p,2)$ with $p$ even and $p<m$, as described in \cite{LT} and \cite{OT}. In this
case the lines $x_1=\zeta x_2$ form two orbits and the coordinate lines form the third.
For well generated groups the classification gives the following result.
}

\begin{theorem}

Let $W$ be nontrivial, irreducible, and well generated. Then the action
of $W$ on $\A$ has either one or two orbits.
\end{theorem}

{
When there are two orbits, the pointwise stabilizers of hyperplanes in
different orbits may have different orders, as in $G_6$, or the same order,
as in $B_n$, $F_4$, and $G_5$ (see the tables below). For the infinite families we can
describe the orbits directly. In coordinates $x_1,\dots,x_n$, the reflecting hyperplanes
of $G(m,p,n)$ have equations
\[
  H_{ij}^{\zeta}:\ x_i=\zeta x_j\ (i<j,\ \zeta^m=1), \qquad
  H_i:\ x_i=0\ \ (\text{only if } p<m).
\]

For the cyclic group $G(m,1,1)=\mathbb Z/m$ with $m\ge2$, the only reflecting hyperplane is $\{0\}$.
For $G(m,1,n)$ with $m,n\ge2$, the $m\binom n2$ hyperplanes $H_{ij}^{\zeta}$ form one orbit, and the $n$ coordinate
hyperplanes $H_i$ form the other. This gives the two orbits for $B_n$ when $m=2$.
The symmetric group $G(1,1,n)=\mathfrak S_n$ with $n\ge2$ acts transitively on the reflecting hyperplanes
of its reflection representation. For $G(m,m,n)$ with $m\ge2$ and $n\ge3$, all $m\binom n2$ hyperplanes $H_{ij}^{\zeta}$
belong to the same orbit. In particular, we obtain the case of $D_n$ by taking $m=2$.
In rank two, the $m$ lines of $G(m,m,2)=I_2(m)$ with $m\ge3$ form a single orbit when $m$ is odd and
two orbits of size $m/2$ when $m$ is even.

For the exceptional groups, we give the sizes of the orbits in the two
tables below, following \cite{LT} and \cite{OT}. For a hyperplane in each orbit we record the order of
a generator of its pointwise stabilizer. The orders in the last column
correspond, in sequence, to the orbit sizes.
}

\medskip
\noindent\textbf{One orbit.}
\begin{center}
\begin{tabular}{lcccll}
\toprule
$W$ & rank & $|W|$ & degrees & $N^*$ & refl.\ order\\
\midrule
$G_4$ & 2 & 24 & 4, 6 & 4 & 3\\
$G_8$ & 2 & 96 & 8, 12 & 6 & 4\\
$G_{16}$ & 2 & 600 & 20, 30 & 12 & 5\\
$G_{20}$ & 2 & 360 & 12, 30 & 20 & 3\\
$G_{23}=H_3$ & 3 & 120 & 2, 6, 10 & 15 & 2\\
$G_{24}$ & 3 & 336 & 4, 6, 14 & 21 & 2\\
$G_{25}$ & 3 & 648 & 6, 9, 12 & 12 & 3\\
$G_{27}$ & 3 & 2160 & 6, 12, 30 & 45 & 2\\
$G_{29}$ & 4 & 7680 & 4, 8, 12, 20 & 40 & 2\\
$G_{30}=H_4$ & 4 & 14400 & 2, 12, 20, 30 & 60 & 2\\
$G_{32}$ & 4 & 155520 & 12, 18, 24, 30 & 40 & 3\\
$G_{33}$ & 5 & 51840 & 4, 6, 10, 12, 18 & 45 & 2\\
$G_{34}$ & 6 & 39191040 & 6, 12, 18, 24, 30, 42 & 126 & 2\\
$G_{35}=E_6$ & 6 & 51840 & 2, 5, 6, 8, 9, 12 & 36 & 2\\
$G_{36}=E_7$ & 7 & 2903040 & 2, 6, 8, 10, 12, 14, 18 & 63 & 2\\
$G_{37}=E_8$ & 8 & 696729600 & 2, 8, 12, 14, 18, 20, 24, 30 & 120 & 2\\
\bottomrule
\end{tabular}
\end{center}

\medskip
\noindent\textbf{Two orbits.}
\begin{center}
\begin{tabular}{lcccll}
\toprule
$W$ & rank & $|W|$ & degrees & orbit sizes & refl.\ orders\\
\midrule
$G_5$ & 2 & 72 & 6, 12 & $4+4$ & $3,\,3$\\
$G_6$ & 2 & 48 & 4, 12 & $6+4$ & $2,\,3$\\
$G_9$ & 2 & 192 & 8, 24 & $12+6$ & $2,\,4$\\
$G_{10}$ & 2 & 288 & 12, 24 & $8+6$ & $3,\,4$\\
$G_{14}$ & 2 & 144 & 6, 24 & $12+8$ & $2,\,3$\\
$G_{17}$ & 2 & 1200 & 20, 60 & $30+12$ & $2,\,5$\\
$G_{18}$ & 2 & 1800 & 30, 60 & $20+12$ & $3,\,5$\\
$G_{21}$ & 2 & 720 & 12, 60 & $30+20$ & $2,\,3$\\
$G_{26}$ & 3 & 1296 & 6, 12, 18 & $9+12$ & $2,\,3$\\
$G_{28}=F_4$ & 4 & 1152 & 2, 6, 8, 12 & $12+12$ & $2,\,2$\\
\bottomrule
\end{tabular}
\end{center}

{
We can now group the families according to the number of orbits.
}

\begin{center}
\renewcommand{\arraystretch}{1.25}
\begin{tabular}{p{0.44\textwidth}p{0.48\textwidth}}
\toprule
\textbf{One orbit on $\A$} & \textbf{Two orbits on $\A$}\\
\midrule
$\mathbb Z/m=G(m,1,1)$;\newline
$\mathfrak S_n=G(1,1,n)$;\newline
$G(m,m,n)$, $n\ge3$ (incl.\ $D_n$);\newline
$I_2(m)=G(m,m,2)$, $m$ odd;\newline
$G_4,G_8,G_{16},G_{20}$;\newline
$G_{23},G_{24},G_{25},G_{27},G_{29},G_{30}$,\newline
$G_{32},G_{33},G_{34},G_{35},G_{36},G_{37}$
&
$G(m,1,n)$, $m\ge2$, $n\ge2$ (incl.\ $B_n$);\newline
$I_2(m)=G(m,m,2)$, $m$ even;\newline
$G_5,G_6,G_9,G_{10},G_{14},G_{17},G_{18},G_{21}$;\newline
$G_{26},\ G_{28}=F_4$\\
\bottomrule
\end{tabular}
\end{center}

{
Among irreducible real reflection groups, the ones with two orbits are $B_n$
with $n\ge2$, $F_4$, and $I_2(m)$ with $m\ge4$ even. Each of the other types has one orbit.

Given a reflection group $G$ acting on $V=\mathbb C^n$, set
\[
V^{\mathrm{reg}}=V\setminus\bigcup_{H\in\mathcal H}H,\qquad
M=V/G,\qquad M^{\mathrm{reg}}=V^{\mathrm{reg}}/G .
\]
The complement of $M^{\mathrm{reg}}$ is called the reflection discriminant. Over $M^{\mathrm{reg}}$ the
quotient map is locally biholomorphic, so the linear coordinates of $V$
give local coordinates on the regular quotient. The natural structure
extends the construction of Dubrovin and Saito and is polynomial on the
whole of $M$. On $M^{\mathrm{reg}}$ it forms a bi-flat pair with the dual structure associated
with the arrangement. The dual structure is defined in terms of a connection
of Dunkl and Kohno and generalizes Veselov's $\vee$ systems. In this way one
obtains a bi-flat pair without assuming the existence of an invariant metric.

The description of the flat coordinates for the natural structure is due to
Kato, Mano, and Sekiguchi \cite{KMS1,KMS2}. See also \cite{KMS}.
}

\begin{theorem}\label{ALcomplexthm1}

Suppose that $G$ is a nontrivial irreducible complex reflection group and
that it is well generated. Its orbit space carries a polynomial flat F-structure $(\nabla,\circ,e,E)$ with a linear Euler field and the following properties.
\begin{enumerate}
\item The flat coordinates of $\nabla$ form a system of basic invariants $(u^1,\dots,u^n)$ of
the group, called generalized Saito coordinates.
\item In these coordinates the unit and the Euler field are given by
\[e=\frac{\partial}{\partial u^n},\qquad E=\sum_{i=1}^n\left(\frac{d_i}{d_n}\right)u^i\frac{\partial}{\partial u^i}.\]
\item By Lemma~7.4 and Theorem~7.5 of \cite{KMS}, the determinant of $E\circ$ gives,
up to a nonzero constant, a reduced equation for the reflection discriminant.
Thus $E\circ$ is invertible precisely on $M^{\mathrm{reg}}$.
\end{enumerate}
\end{theorem}

{
We next describe the second flat structure on the complement of the
reflecting hyperplanes. As in the Coxeter case, we express the product in terms of the
projections and their weights. A second choice of weights specifies the connection.
}

\begin{theorem}\label{ALcomplexthm2}

Let $G$ be a nontrivial irreducible complex reflection group acting on $\mathbb{C}^n$.
Consider the data
\[\left(\nabla^*
=\nabla^0-\sum_{H\in \mathcal{H}}\frac{d\alpha_H}{\alpha_H}\otimes\tau_H\pi_H,\quad *=\sum_{H\in \mathcal{H}}\frac{d\alpha_H}{\alpha_H}\otimes\sigma_H\pi_H,\quad E=\sum_k p^k\frac{\partial}{\partial p^k}\right),\]
where we use the following notation.
\begin{itemize}
\item $\mathcal{H}$ is the set of reflecting hyperplanes of $G$,
\item $\alpha_H$ is a linear form defining the hyperplane $H$,
\item $\pi_H$ is the projection onto the orthogonal complement of $H$ with respect to a positive definite Hermitian form invariant under $G$,
\item the weights $\sigma_H$ and $\tau_H$ are invariant under $G$ and satisfy
\begin{equation}\label{sums}
\sum_{H\in \mathcal{H}}\sigma_H\pi_H=\sum_{H\in \mathcal{H}}\tau_H\pi_H=\mathrm{Id},
\end{equation}
\item $\nabla^0$ is the standard flat connection on $\mathbb{C}^n$.
\end{itemize}
These data define a flat F-structure on $V^{\mathrm{reg}}$ and descend to $M^{\mathrm{reg}}$.
\end{theorem}

\begin{proof}

For every $s$, the connection $\nabla^*+s*$ is a Dunkl connection with invariant weights
$\tau_H-s\sigma_H$. Thus it is flat by \cite{CHL}, and in rank one this follows directly. Each
projection contributes a term that is symmetric in the two arguments of
the product, so every connection in the pencil is torsion free. As in the
curvature calculation for a flat F-manifold, the coefficients of $s$ and $s^2$
give, respectively, compatibility and associativity. To see that $E$ is the
unit, evaluate the product at this field and use the first sum in \eqref{sums}. Applying $\nabla^*$
to $E$ and using the second sum, we see that the unit is parallel. Finally,
invariance of the tensors and the connection under $G$ allows them to descend
to $M^{\mathrm{reg}}$.
\end{proof}

{
For the natural flat structure, a shift by the dual product gives the
connection whose flat coordinates are the linear coordinates of $V$.
We recall the characterization that we shall use below.
}

\begin{remark}\label{rem:KMSKMS}

Kato, Mano, and Sekiguchi obtained the structure of
Theorem~\ref{ALcomplexthm1} from systems of Okubo type and their isomonodromic deformations
\cite{KMS1,KMS2}. In the equivalent formulation of Konishi, Minabe, and Shiraishi \cite{KMS},
a \emph{Saito structure without metric} is a flat F-manifold endowed with
a linear Euler field. The corresponding \emph{almost Saito structure}
consists of the dual product with unit $E$, the connection $\nabla^*+\nu\,*$ for a fixed
value of $\nu$, and the distinguished field $e$. These connections occur as the
members $\Gamma(0,\nu)$ of the family of connections of Gauss and Manin in Section~\ref{sec:GM}.

We refer to this natural structure as the \emph{KMS flat structure}.
For the normalized unit and Euler field, this structure is characterized in \cite{KMS}
by two properties. Its flat coordinates are polynomial, and the linear
functions $p^1,\dots,p^n$ on the covering space provide local flat coordinates for $\nabla^*+\frac1{d_n}\,*$ on
the regular quotient. In other words, the trivial connection on $V$ is the almost
dual connection at weight $1/d_n$, which is the degree of the linear functions
with respect to $E$. We shall use this characterization in Section~\ref{sec:ALconjectures}.
\end{remark}

\subsection{The conjecture}\label{sec:conjecture}

{
We formulated the conjecture for Coxeter groups in \cite{ALrefl} and for well generated
complex reflection groups in \cite{ALsurvey}. It requires the natural structure on $M^{\mathrm{reg}}$
to extend polynomially to $M$ and to have basic invariants as flat coordinates.
We fix the unit $e$ as in Theorem~\ref{ALcomplexthm1} and define $E$ to be the radial field
divided by the largest invariant degree $h$. Starting with the dual data of
Theorem~\ref{ALcomplexthm2}, we multiply the product by $h$ so that the unit of $*$ is $E$. This
rescaling preserves the connection. In the families considered below, we
keep the dual product, $e$, and $E$ fixed and vary the weights of the connection.
}

\begin{conjecture}\label{conj:AL}

Let $G$ be a nontrivial finite complex reflection group acting irreducibly on
its reflection representation $V$, and suppose that $G$ is well generated.
Then the regular quotient $M^{\mathrm{reg}}$ carries a family of bi-flat structures of the
form described above. Their natural flat structures extend polynomially to
$M=V/G$. The family depends on $\mu$ parameters, where $\mu=|\mathcal{H}/G|-1$ and $|\mathcal{H}/G|$ denotes the number of
orbits of $G$ on the set of reflecting hyperplanes $\mathcal{H}$.
\end{conjecture}

{
To see how the number of parameters arises, recall that an invariant
connection of Dunkl and Kohno has one weight $\tau_{\mathcal O}$ for each orbit $\mathcal O$.
By irreducibility and Schur's lemma, $\sum_{H\in\mathcal O}\pi_H=(|\mathcal O|/n)\Id$. Therefore the condition that
$E$ be parallel with respect to $\nabla^*$ gives a single relation among the weights,
\[
\sum_{\mathcal O}|\mathcal O|\tau_{\mathcal O}=n.
\]
The space of weights satisfying this relation has dimension $\mu$. Once we
fix the dual product and $e$, the reconstruction following Lemma~\ref{lem:AHE} shows
that a given connection admits at most one compatible natural structure.
Thus, for fixed dual product, $e$, and $E$, a family of bi-flat structures of
the prescribed form depends on at most $\mu=|\mathcal H/G|-1$ parameters. We obtained this
bound in \cite{ALsurvey}. Conjecture~\ref{conj:AL} asserts that the bound is attained. The following
cases were established in \cite{ALrefl} and \cite{Bn}.
}

\begin{theorem}\label{ALComplexPrincipal}

Conjecture~\ref{conj:AL} holds for irreducible Weyl groups of ranks $2$ and
$3$, for $B_4$, and for the dihedral groups $I_2(m)$ with $m\ge3$.
It also holds for the exceptional well generated complex reflection
groups of ranks $2$ and $3$ and for $G(3,1,2)$ and $G(3,1,3)$.
\end{theorem}

{
For a Coxeter group with a single orbit on $\mathcal H$, the Frobenius structure of
Theorem~\ref{Dcoxeterth} and its dual with $\nabla^*E=0$ give the required bi-flat pair. Among the groups
with two orbits, we construct the family for $B_n$ in every rank and obtain
explicit formulas for $F_4$. Our construction for $G(m,1,n)$ applies whenever $m,n\ge2$.
The result for a single orbit also applies to $G_{29},G_{32},G_{33},G_{34}$ and to $G(m,m,n)$ with $m,n\ge3$.

In Table~\ref{tab:coverage} we relate these results to the classification recalled in
Section~\ref{sec:orbitspaces}. The table concerns groups of rank at least two. In some cases
the same group appears under different names. For the cyclic groups of rank
one, the conjecture follows from Theorem~\ref{thm:naturalDunklKohno}.
}

\begin{table}[ht]
\centering
\renewcommand{\arraystretch}{1.25}
\begin{tabular}{p{0.09\textwidth}>{\raggedright\arraybackslash}p{0.40\textwidth}>{\raggedright\arraybackslash}p{0.41\textwidth}}
\toprule
orbits & groups & reference\\
\midrule
1 & Coxeter groups $A_n$, $D_n$, $E_6,E_7,E_8$, $H_3,H_4$, $I_2(m)$ with $m$ odd &
Dubrovin \cite{DCoxeter,Dalm} and Theorem~\ref{thm:naturalDunklKohno}\\
1 & $G_4,G_8,G_{16},G_{20}$ (rank 2), $G_{24},G_{25},G_{27}$ (rank 3) &
\cite{ALrefl} and Theorem~\ref{thm:naturalDunklKohno}\\
1 & $G_{29},G_{32},G_{33},G_{34}$ & Theorem~\ref{thm:naturalDunklKohno}\\
1 & $G(m,m,n)$, $m\ge3$, $n\ge3$ & Theorems~\ref{thm:naturalDunklKohno} and~\ref{thm:Gmmn}\\
\midrule
2 & $I_2(m)$ with $m$ even & \cite{Bn}\\
2 & $G_5,G_6,G_9,G_{10},G_{14},G_{17},G_{18},G_{21}$ (rank 2), $G_{26}$ (rank 3) & \cite{ALrefl}\\
2 & $G(m,1,n)$, $m\ge2$, $n\ge2$ (including $B_n$) & Theorem~\ref{thm:Bn} and Proposition~\ref{prop:Gm1n}.
Earlier examples in \cite{ALrefl,Bn}\\
2 & $F_4=G_{28}$ & Reported in \cite{Bn}. Explicit formulas in Proposition~\ref{prop:F4}\\
\bottomrule
\end{tabular}
\caption{\texorpdfstring{Status of Conjecture~\ref{conj:AL} for the well generated irreducible complex reflection groups.}{Status of the conjecture for the well generated irreducible complex reflection groups.}}
\label{tab:coverage}
\end{table}

{
In the constructions below we distinguish the two cases according to the
number of orbits. If $\mu=0$, Theorem~\ref{thm:naturalDunklKohno} identifies the dual of the natural
structure of Theorem~\ref{ALcomplexthm1} and thus proves the conjecture. If $\mu=1$, the family
required by the conjecture depends on one parameter.
}

\section{\texorpdfstring{$F$ superpotentials}{F superpotentials}}\label{sec:def}

{
Let $M$ be a holomorphic bi-flat F-manifold endowed with the structures $(\nabla,\circ,e)$ and
$(\nabla^*,*,E)$. We work on a domain where $E\circ$ is invertible and use the notation of
Section~\ref{sec:almostdual}. In particular, $\nabla^*E=0$. Given a connection $D$, we say that a
function $f$ is flat if $D\,df=0$. When $D$ is flat and torsion free, any such function whose
differential is nonzero can be completed to a local system of flat coordinates.

In Dubrovin's construction one obtains the superpotential $\lambda(p,u)$ by inverting a function
$p(u,\lambda)$ with respect to $\lambda$. For each fixed $\lambda$ this function is a flat coordinate
of the pencil $g-\lambda\eta^{-1}$. On a bi-flat F-manifold the connections of Gauss and Manin type
studied in \cite{AL} provide the corresponding family of flat connections. We use the flow of the unit
to relate flat functions on different slices. We translate a flat function along this flow and invert the resulting function locally
wherever its derivative with respect to $\lambda$ is nonzero.

We recall the family of connections $\Gamma(\lambda,\nu)$ in Section~\ref{sec:GM}. In
Section~\ref{sec:GMflat} we prove the translation formula, which leads to the definition of the
superpotential in Section~\ref{sec:defsp}. We give an equivalent definition in terms of the extended
connection in Section~\ref{sec:ext} and characterize the inverse by three equations in
Section~\ref{sec:char}. The local inversion makes sense on any holomorphic bi-flat F-manifold. Suppose,
in addition, that the inverse extends holomorphically to $n$ simple critical points with distinct
critical values. In this case the critical values form a system of canonical coordinates and the
product is semisimple (Theorem~\ref{thm:main}). Under the same hypotheses we obtain residue formulas
for the flat connection in terms of the superpotential in Section~\ref{subsec:reconstruction}.
Finally, in Section~\ref{sec:Dub} we show how Dubrovin's construction arises in the Frobenius case.
}

\subsection{The connection of Gauss and Manin}\label{sec:GM}

{
We recall from \cite{AL} a family of connections of Gauss and Manin type which depends on two
parameters and interpolates between $\nabla^*$ and $\nabla$. Let $b$ and $a$ denote the Christoffel
symbols of these two connections, respectively, in arbitrary coordinates. For $\nu\in\mathbb C$, set
$b_\nu=b+\nu\,c^*$ and define
\begin{equation}\label{eq:GM}
\Gamma^j_{hk}(\lambda,\nu):=b^j_{\nu,hk}+\lambda\,(L_\lambda^{-1})^s_h\,
\big(b^j_{\nu,sk}-a^j_{sk}\big),
\end{equation}
on the complement of $\Sigma=\{(x,\lambda):\det(E\circ-\lambda\Id)=0\}$. In a semisimple chart with
canonical coordinates $r^i$, the graphs $\lambda=r^i$ form $\Sigma$. We recover $\nabla^*$ at
$\lambda=\nu=0$. As $\lambda\to\infty$ we have $\lambda L_\lambda^{-1}\to-\Id$, so the family tends to
$\nabla$. For every pair $(\lambda,\nu)$, the connection $\Gamma(\lambda,\nu)$ is flat and torsion
free on its domain, as proved in \cite{AL}. We write $X*_\lambda Y=(E-\lambda e)^{-1}\circ X\circ Y$.
In particular, $*_0=*$.

Since $e$ is the unit of $\circ$, we have $L_\lambda Z=(E-\lambda e)\circ Z$. Outside $\Sigma$, the
identity $L_\lambda(L_\lambda^{-1}e)=e$ gives $L_\lambda^{-1}e=(E-\lambda e)^{-1}$ in the tangent
algebra. Indeed, for any $Z$, associativity gives
\[
L_\lambda\bigl((E-\lambda e)^{-1}\circ Z\bigr)
=\bigl((E-\lambda e)\circ(E-\lambda e)^{-1}\bigr)\circ Z=Z.
\]
It follows that $L_\lambda^{-1}Z=(E-\lambda e)^{-1}\circ Z$. If we apply this formula to $X$ and to
$X\circ Y$, associativity gives
\begin{equation}\label{eq:aux1}
X*_\lambda Y=(L_\lambda^{-1}X)\circ Y=L_\lambda^{-1}(X\circ Y).
\end{equation}

To exhibit the dependence on $\lambda$, let us choose flat coordinates for $\nabla$. In these
coordinates $a=0$, and from \eqref{eq:dual} we obtain $b_\nu=(\nu\Id-Q)c^*$, or, in operator notation,
$b_{\nu,\cdot k}=(\nu\Id-Q)L^{-1}C_k$. Since $\Id+\lambda L_\lambda^{-1}=LL_\lambda^{-1}$ and the
multiplication operators commute, we obtain
\begin{equation}\label{eq:op}
\Gamma^m_{hk}(\lambda,\nu)=(\nu-Q)^m_p\,(L_\lambda^{-1})^p_l\,c^l_{hk}.
\end{equation}

\begin{remark}\label{rem:KMS}
It follows from \eqref{eq:dual} that \eqref{eq:GM} has the covariant expression
\[
\Gamma(\lambda,\nu)_XY=\nabla_XY+\nu\,X*_\lambda Y-\nabla_{X*_\lambda Y}E .
\]
This is the expression given in Proposition~3.1 of \cite{KMS}.
\end{remark}

We now describe how the eigenvalues of $\nabla E$ give rise to flat functions for particular members
of the family.
\begin{remark}\label{rem:eigen}
Suppose that $\nabla E$ is diagonalizable. Choose flat coordinates $u^1,\dots,u^n$ for $\nabla$ such
that $Q=\nabla E=\mathrm{diag}(d_1,\dots,d_n)$. In these coordinates \eqref{eq:op} gives
\[
\Gamma^i_{hk}(\lambda,\nu)=(\nu-d_i)\,(L_\lambda^{-1}C)^i_{hk},
\]
so that $\Gamma^i_{hk}(\lambda,d_i)=0$. To see what this means for the coordinate functions, recall
that a function $f$ is flat for a torsion free connection $\Gamma$ if and only if
$\partial_h\partial_kf=\Gamma^m_{hk}\,\partial_mf$. For $f=u^i$ the second derivatives in the
coordinates $u$ vanish, while $\partial_mu^i=\delta^i_m$. The flatness equation therefore reduces to
$\Gamma^i_{hk}=0$. Thus $u^i$ is a flat coordinate for $\Gamma(\lambda,d_i)$ for every $\lambda$. At
$\lambda=0$ this says that $u^i$ is a flat coordinate for $\nabla^*+d_i\,c^*$. Moreover,
$Qe=\nabla_eE=\nabla_Ee+[e,E]=e$ by torsion freeness and $\nabla e=0$. The unit is therefore an
eigenvector with eigenvalue $1$, and we may choose the coordinates so that $e=\partial_{u^n}$ and
$d_n=1$.
\end{remark}
}

\subsection{\texorpdfstring{Flat functions on the slices}{Flat functions on the slices}}\label{sec:GMflat}

For $\lambda_0\in\mathbb C$, the \emph{slice at $\lambda_0$} is the submanifold $M\times\{\lambda_0\}$
of $M\times\mathbb C$. Its \emph{regular part} is
\[
\bigl(M\times\{\lambda_0\}\bigr)\setminus\Sigma
=\{(x,\lambda_0):x\in M,\ \det L_{\lambda_0}(x)\ne0\}.
\]
Under the projection $(x,\lambda_0)\mapsto x$ the regular part of the slice is identified with the
open subset of $M$ where $\Gamma(\lambda_0,\nu)$ is defined. A function on this regular part is flat
at weight $\nu$ if its differential is parallel for $\Gamma(\lambda_0,\nu)$.

{
Fix $\nu_0\in\mathbb C$ and define
 \begin{equation}\label{eq:aux2}N:=\nu_0\Id-Q,\end{equation}
where $Q=\nabla E$ as in Section~\ref{sec:almostdual}. For each fixed $\lambda$ we may choose flat
coordinates for $\Gamma(\lambda,\nu_0)$ up to an affine transformation. In order to invert one of
these functions with respect to $\lambda$, we must first relate our choices on different slices. For
this purpose we use the local flow $\phi_t$ of the unit $e$.

\begin{theorem}\label{thm:ext}
{
Let $(x_0,\lambda_0)\notin\Sigma$ and choose a neighborhood of this point on which the map
$(x,\lambda)\mapsto(\phi_{\lambda_0-\lambda}(x),\lambda)$ is defined.
\begin{enumerate}
\item[(i)] The connections on the slices of this neighborhood satisfy
\begin{equation}\label{eq:GMtransl}
\Gamma(\lambda,\nu_0)=\phi_{\lambda_0-\lambda}^{*}\Gamma(\lambda_0,\nu_0).
\end{equation}
\item[(ii)] If $P_{\lambda_0}$ is a flat function of $\Gamma(\lambda_0,\nu_0)$ near $x_0$, then $P(x,\lambda)=P_{\lambda_0}(\phi_{\lambda_0-\lambda}(x))$ is flat on each slice and satisfies
\begin{equation}\label{eq:transl}
\partial_\lambda P=-\partial_eP.
\end{equation}
Moreover, in flat coordinates of $\nabla$ its differential in the base variables, written as the row
vector $\xi=d_uP$, satisfies
\begin{equation}\label{eq:xisys}
\partial_h\xi=\xi NL_\lambda^{-1}C_h,\qquad
\partial_\lambda\xi=-\xi NL_\lambda^{-1}.
\end{equation}
\item[(iii)] There is a constant $\kappa$ for which
$(E-\lambda_0 e)P_{\lambda_0}=\nu_0P_{\lambda_0}+\kappa$ and
$(E+\lambda\partial_\lambda)P=\nu_0P+\kappa$. If $\nu_0\ne0$, we may normalize $\kappa$ to zero by
adding a constant to $P_{\lambda_0}$.
\end{enumerate}
}
\end{theorem}

\begin{proof}
(i) Choose flat coordinates $u$ of $\nabla$ near $x_0$. Since $\nabla e=0$, the components $e^m$ of
the unit are constant. Thus its flow $\phi_t$ is the translation $u\mapsto u+te$, where $te$ denotes
the tuple $(te^m)_m$. For $(u,\lambda)$ near the reference point, set $v:=u+(\lambda_0-\lambda)e$. We
first compare the structure constants of the product and the components of the Euler field at $u$ and
$v$.

Using the symmetry $\partial_m c^\ell_{hk}=\partial_h c^\ell_{mk}$ and the identity
$e^mc^\ell_{mk}=\delta^\ell_k$ for the unit, we obtain
\[
\partial_e c^\ell_{hk}
=e^m\partial_m c^\ell_{hk}
=\partial_h(e^mc^\ell_{mk})
=\partial_h\delta^\ell_k=0.
\]
It follows that $c^\ell_{hk}(v)=c^\ell_{hk}(u)$. Moreover, $Q^m_j=\partial_j E^m$ is constant because
$\nabla Q=0$. We may therefore write $E^m(u)=Q^m_j u^j+b_0^m$, from which we obtain
\[
E^m(v)=E^m(u)+(\lambda_0-\lambda)Q^m_j e^j.
\]
By torsion freeness and $\nabla e=0$ we have $Qe=\nabla_eE=\nabla_Ee+[e,E]=e$. It follows that
$E(v)-\lambda_0 e=E(u)-\lambda e$ in our flat coordinates. Since the structure constants at $u$ and
$v$ agree, the operators of multiplication by these two elements also agree. Thus
\[
L_{\lambda_0}(v)=L_\lambda(u).
\]
In particular, translation sends the regular part of the slice at $\lambda$ into the regular part of
the slice at $\lambda_0$ wherever the translation is defined.

Since $N=\nu_0\Id-Q$ is constant, we can substitute these identities in \eqref{eq:op} to obtain
\begin{align*}
\Gamma^m_{hk}(\lambda_0,\nu_0)(v)
&=N^m_p\bigl(L_{\lambda_0}(v)^{-1}\bigr)^p_\ell c^\ell_{hk}(v)\\
&=N^m_p\bigl(L_\lambda(u)^{-1}\bigr)^p_\ell c^\ell_{hk}(u)
=\Gamma^m_{hk}(\lambda,\nu_0)(u).
\end{align*}
For fixed $\lambda$ the map $u\mapsto v$ has the identity as its Jacobian and all its second
derivatives vanish. The expression on the left is therefore the Christoffel symbol of the connection
pulled back by $\phi_{\lambda_0-\lambda}$. This proves (i).

(ii) Set $F=P_{\lambda_0}$, so that $P(u,\lambda)=F(v)$. Differentiating the translation at fixed
$\lambda$ gives $\partial_hP(u,\lambda)=(\partial_hF)(v)$ and
$\partial_h\partial_kP(u,\lambda)=(\partial_h\partial_kF)(v)$. Since $F$ is flat on the reference
slice, part (i) gives
\begin{align*}
\partial_h\partial_kP(u,\lambda)
&=\Gamma^m_{hk}(\lambda_0,\nu_0)(v)\,(\partial_mF)(v)\\
&=\Gamma^m_{hk}(\lambda,\nu_0)(u)\,\partial_mP(u,\lambda).
\end{align*}
Thus $P(\cdot,\lambda)$ is flat. Now keep $u$ fixed and observe that $\partial_\lambda v^m=-e^m$. By
the chain rule,
\[
\partial_\lambda P(u,\lambda)
=-e^m(\partial_mF)(v),
\qquad
\partial_eP(u,\lambda)=e^m(\partial_mF)(v).
\]
Comparing the two derivatives gives $\partial_\lambda P=-\partial_eP$.

Set $\xi_k=\partial_kP$. Using \eqref{eq:op} in the flatness equation gives
$\partial_h\xi_k=\xi_m(NL_\lambda^{-1})^m_\ell c^\ell_{hk}$, which is the first equation of
\eqref{eq:xisys}. To obtain the second equation, differentiate \eqref{eq:transl} with respect to
$u^k$. The components $e^h$ are constant and the partial derivatives commute, so
$\partial_\lambda\xi_k=-e^h\partial_h\xi_k$. If we contract the first equation with $e^h$ and use
$e^h C_h=\Id$, we obtain $\partial_\lambda\xi=-\xi NL_\lambda^{-1}$.

(iii) On the reference slice set $E_0=E-\lambda_0 e$ and $\zeta_k=\partial_kF$. Since
$\partial_h E_0^k=Q^k_h$ and $c^\ell_{hk}E_0^k=(L_{\lambda_0})^\ell_h$, the flatness equation for $F$
gives the following derivative of $E_0^k\zeta_k$.
\begin{align*}
\partial_h\big(E_0(F)\big)
&=Q^k_h\zeta_k+E_0^k\partial_h\zeta_k\\
&=Q^k_h\zeta_k+
  \zeta_m\bigl(NL_{\lambda_0}^{-1}\bigr)^m_\ell c^\ell_{hk}E_0^k\\
&=\zeta_m(Q+N)^m_h=\nu_0\partial_hF.
\end{align*}
It follows that $E_0(F)-\nu_0F$ is constant on a connected neighborhood in the reference slice. We
denote its value by $\kappa$.

To obtain the corresponding identity for $P$, observe that $E(u)=E_0(v)+\lambda e$ by part (i).
Applying the chain rule with $\lambda$ fixed gives
\[
E(P)(u,\lambda)=\big(E_0(F)\big)(v)+\lambda\,(\partial_eF)(v).
\]
By part (ii), the last term cancels when we add $\lambda\partial_\lambda P(u,\lambda)$. Therefore
\[
(E+\lambda\partial_\lambda)P(u,\lambda)
=\big(E_0(F)\big)(v)=\nu_0F(v)+\kappa=\nu_0P(u,\lambda)+\kappa.
\]
If $\nu_0\ne0$, we may replace $F$ by $F+\kappa/\nu_0$. This preserves the differential and gives the
homogeneous identity $E_0(F+\kappa/\nu_0)=\nu_0(F+\kappa/\nu_0)$.
\end{proof}

Equation \eqref{eq:transl} describes how flat coordinates on different slices are related. Once the
flat functions on one slice have been chosen, the equation determines them along the flow of
$e+\partial_\lambda$. In the case of Frobenius manifolds the corresponding formula is given in
\cite{Dalm}, namely $p(t,\lambda)=\hat p(t^1-\lambda,t^2,\dots,t^n)$. When $\nu_0\ne0$ we can choose
the additive constant so that $\kappa=0$. The homogeneity equation then reads
$P=\nu_0^{-1}\langle d_xP,E-\lambda e\rangle$, as in \cite{Dalm}.

\begin{remark}\label{rem:pullback}
Suppose that $\phi_{-\lambda}$ is defined and takes values in the regular part of the zero slice.
Equation \eqref{eq:GMtransl} then reads $\Gamma(\lambda,\nu_0)=\phi_{-\lambda}^*\Gamma(0,\nu_0)$. On
these slices flatness thus follows directly from the flatness of the deformed dual connection.
\end{remark}
}

\subsection{Definition of the superpotential}\label{sec:defsp}

Fix $\nu_0\in\mathbb C$. The connection on the zero slice is $\Gamma(0,\nu_0)=\nabla^*+\nu_0c^*$.
Choose a local flat coordinate $P_0$ of this connection and denote by $\kappa$ the constant determined
by $E(P_0)=\nu_0P_0+\kappa$ as in Theorem~\ref{thm:ext}(iii). It follows from
Theorem~\ref{thm:ext}(ii) that
\[
P(x,\lambda):=P_0\bigl(\phi_{-\lambda}(x)\bigr)
\]
gives a flat coordinate on each slice for the corresponding connection of Gauss and Manin
$\Gamma(\lambda,\nu_0)$. In flat coordinates $u$ of $\nabla$ the flow of the unit is a translation, so
that $P(u,\lambda)=P_0(u-\lambda e)$, where $e$ denotes the constant tuple of components of the unit.
Choose $(x_0,\lambda_0)\notin\Sigma$ such that $y_0=\phi_{-\lambda_0}(x_0)$ belongs to the domain of
$P_0$ and $(\partial_eP_0)(y_0)\ne0$. Set $q_0=P(x_0,\lambda_0)$. More generally,
Theorem~\ref{thm:ext} allows us to start the same construction on any regular slice.

\begin{definition}\label{def:sp}
We define the \emph{local superpotential of weight $\nu_0$} associated with $P_0$ at $(x_0,\lambda_0)$
as the holomorphic inverse $\Lambda(q,x)$ of $P(x,\lambda)$ with respect to $\lambda$. It is defined
near $(q_0,x_0)$ by the identity
\[
\Lambda\bigl(P(x,\lambda),x\bigr)=\lambda
\qquad\text{for $(x,\lambda)$ near $(x_0,\lambda_0)$}.
\]
\end{definition}
We also call the analytic continuation of this function a superpotential. The continuation may be
multivalued. We refer to its \emph{sheets} and use the term \emph{branch} for a local inverse
$P(x,\lambda)$ of $\Lambda$. Different choices of the flat function $P_0$ may lead to the same
continued superpotential. Indeed, in Sections~\ref{sec:Bn} and~\ref{sec:Gmmn} monodromy exchanges the
branches $P_0=w_j$, all of which give $\lambda(z)$ (Remark~\ref{rem:norm}). We call $\kappa$ the
\emph{constant} of the superpotential. Its \emph{primary form} is the differential of $P(x,\cdot)$
taken with $x$ fixed.

\begin{remark}\label{rem:spgeometry}
Let $U$ be an open neighborhood of $(x_0,\lambda_0)$ in $(M\times\mathbb C_\lambda)\setminus\Sigma$ on
which $P$ is defined. Consider the projection $\pi\colon U\to M$ given by $(x,\lambda)\mapsto x$. Its
fiber $U_x=U\cap(\{x\}\times\mathbb C_\lambda)$ is the local curve over $x$. Along this curve $x$ is
fixed and $\lambda$ varies. On the slices introduced in Section~\ref{sec:GMflat}, we fix $\lambda$ and
allow $x$ to vary.

Set $y=\phi_{-\lambda}(x)$. If we keep $x$ fixed, the flow parameter is $t=-\lambda$, so the curve
$\lambda\mapsto\phi_{-\lambda}(x)$ has velocity $-e_y$ at $y$. By the chain rule we obtain
\[
\partial_\lambda P(x,\lambda)
=\frac{d}{d\lambda}P_0\bigl(\phi_{-\lambda}(x)\bigr)
=dP_0\big|_y(-e_y)=-(\partial_eP_0)(y).
\]
To compute $\partial_eP$ we keep $\lambda$ fixed and let $x$ vary along $s\mapsto\phi_s(x)$. Since
$\phi_{-\lambda}(\phi_s(x))=\phi_s(y)$, we obtain
\[
\partial_eP(x,\lambda)
=\left.\frac{d}{ds}P(\phi_s(x),\lambda)\right|_{s=0}
=\left.\frac{d}{ds}P_0(\phi_s(y))\right|_{s=0}
=(\partial_eP_0)(y).
\]
Comparing these expressions gives $\partial_\lambda P=-\partial_eP$. Our assumption at $y_0$ therefore
implies that $\partial_\lambda P(x_0,\lambda_0)\ne0$. In the particular case $\lambda_0=0$ the
assumption reads $(\partial_eP_0)(x_0)\ne0$.

The map
\[
\Psi\colon U\longrightarrow\widetilde U,\qquad
\Psi(x,\lambda)=(x,q)=(x,P(x,\lambda))
\]
has invertible differential at $(x_0,\lambda_0)$ since it fixes $x$ and
$\partial_\lambda P(x_0,\lambda_0)\ne0$. After shrinking $U$, the holomorphic inverse function theorem
gives a biholomorphism $\Psi$ onto an open subset $\widetilde U\subset M\times\mathbb C_q$ with
inverse $\Psi^{-1}(x,q)=(x,\Lambda(q,x))$. Both maps commute with the projection to $M$. Thus
$\lambda$ and $q$ describe the same local curve over $x$. To interpret the inverse, fix $q$. The point
$\phi_{-\Lambda(q,x)}(x)$ belongs to the level hypersurface $P_0=q$. Since $\partial_eP_0$ is nonzero
there, the trajectory of the unit meets the hypersurface transversely. In particular, $\Lambda(q,x)$
is the unique complex time near $\lambda_0$ at which the flow of $-e$ starting at $x$ reaches this
level locally. When we write $\Lambda(P,x)$, its first argument denotes the coordinate $q$.

Fix $x$ and set $\Lambda_x(q)=\Lambda(q,x)$. Choose a disk $D\subset\mathbb C_\lambda$ small enough
that $\{x\}\times D\subset U$. The curve over this disk is parametrized either by $\lambda\in D$ or by
$q\in Q_x=P(x,D)$. In the product of the two coordinate domains it is the graph
\[
\bigl\{(q,\lambda)\in Q_x\times D:\lambda=\Lambda_x(q)\bigr\}.
\]
The function $\Lambda_x\colon Q_x\to D$ gives the projection of the graph to the $\lambda$ plane when
we use $q$ as a coordinate. It is a biholomorphism onto $D$, so this portion of the curve is a
covering with one sheet.

Let $j_x\colon U_x\hookrightarrow U$ denote the inclusion. The primary form is the relative one form
$\omega$ over $M$ whose restriction to the curve over $x$ is given by
\[
\omega_x=j_x^*(dP)=d(P\circ j_x)
=\partial_\lambda P(x,\lambda)\,d\lambda.
\]
Indeed, $d(P\circ j_x)$ is the differential taken with $x$ fixed. To form the pullback we evaluate
$dP$ at points of $U_x$ on vectors tangent to $U_x$. If we use $q$ as the coordinate, this gives
$\omega_x=dq$.

Differentiating the identity $\Lambda_x(P(x,\lambda))=\lambda$ at fixed $x$ on the inverse chart gives
\[
\Lambda_x'\bigl(P(x,\lambda)\bigr)\,\partial_\lambda P(x,\lambda)=1.
\]
Thus $\Lambda_x'$ is nonzero throughout the inverse chart. In Theorem~\ref{thm:main} we assume that,
in the $q$ coordinate, the superpotential extends holomorphically to a neighborhood of a simple
critical point $q_*$. At this point $\Lambda_x'(q_*)=0$ and $\Lambda_x''(q_*)\ne0$. Setting
$\lambda_*=\Lambda_x(q_*)$, we can write
\[
\Lambda_x(q)-\lambda_*=(q-q_*)^2h(q),\qquad
h(q_*)=\tfrac12\Lambda_x''(q_*)\ne0.
\]
Choose a holomorphic square root of $h$ near $q_*$. In the resulting coordinate $w=(q-q_*)\sqrt{h(q)}$
the map to the $\lambda$ plane takes the form $\lambda=\lambda_*+w^2$. Over a sufficiently small
punctured disk about $\lambda_*$, this local description gives a covering of degree two. Over the full
disk we obtain a branched covering with ramification index two at $q_*$. The same disk may have
further preimages elsewhere on the continued curve.

At $q_*$ the differential $d\lambda=\Lambda_x'(q)\,dq$ vanishes and $\lambda$ is no longer a
coordinate. The primary form $\omega_x=dq$ is holomorphic and nonzero at this point. Away from the
critical value we can write it as $\partial_\lambda P\,d\lambda$ on each inverse branch.
\end{remark}

\begin{remark}\label{rem:weight}
We can interpret the weight $\nu_0$ as the degree of the primary form. Indeed, set
$W=E+\lambda\partial_\lambda$. By Theorem~\ref{thm:ext}(iii) we have $\mathcal L_W(dP)=d(WP)=\nu_0dP$.
Since $W$ projects to $E$ on $M$, its flow takes fibers to fibers. Restricting the identity to the
fibers gives the corresponding homogeneity of the primary form.
\end{remark}

\begin{remark}\label{rem:norm}
Let us examine the dependence of the superpotential on the weight $\nu_0$ and the choice of flat
function $P_0$.

For any weight $\nu_0$ the connection $D=\nabla^*+\nu_0c^*$ is flat. At each point of its domain we
may choose a flat differential such that $dP_0(e)\ne0$. A local superpotential therefore exists for
every weight. To describe the freedom in this choice, take flat coordinates $y^1,\dots,y^n$ for $D$.
In these coordinates $D\,df=0$ means that all the second partial derivatives of $f$ vanish. Thus on a
connected flat chart the flat functions are precisely the functions $\sum_j c_j y^j+c_0$. They form a
vector space of dimension $n+1$, which includes the constants. If we replace $P_0$ by $aP_0+b$ with
$a\ne0$, the new inverse is $\Lambda((P-b)/a,x)$ and the new constant is $a\kappa-\nu_0 b$. When
$\nu_0\ne0$ we obtain the homogeneous normalization by taking $b=\kappa/\nu_0$ and $a=1$. At a fixed
base point the equivalence classes modulo affine changes of $P$ therefore correspond to the directions
$[c_1:\dots:c_n]$ of flat differentials for which $dP_0(e)\ne0$ at that point.

Suppose, for instance, that $\nabla E$ is diagonalizable and choose flat coordinates adapted to the
unit. By Remark~\ref{rem:eigen}, when $\nu_0=d_i$ the coordinate $u^i$ is flat for $D$. It gives a
superpotential provided $\partial_eu^i\neq0$. With $e=\partial_{u^n}$ this condition holds only for
$i=n$, in which case $\nu_0=d_n=1$. We then have $P=u^n-\lambda$ and obtain $\Lambda=u^n-P$.

At weight zero the constant $\kappa=dP_0(E)$ is invariant under the addition of a constant to $P_0$.
Write $E=\sum_j E^j\partial_{y^j}$ in flat coordinates $y^j$ of $\nabla^*$. The components $E^j$ are
constant, and for $P_0=\sum_j c_j y^j+c_0$ we have $\kappa=\sum_jc_jE^j$. Thus $\kappa=0$ defines a
hyperplane in the space of flat differentials. Wherever $e$ and $E$ are linearly independent, this
hyperplane contains a differential for which $dP_0(e)\ne0$. If $\kappa\ne0$, we may write
$P_0=\kappa\tau+f$, where $E(\tau)=1$ and $E(f)=0$. The logarithmic flat function $P_0=\log z$ for
$\beta=0$ in Section~\ref{sec:Bn} illustrates this case.

By Corollary~\ref{cor:flat} the inverse branch at a simple zero of a continued superpotential is flat
for $D$ and satisfies the affine homogeneity equation. In Sections~\ref{sec:Bn} and~\ref{sec:Gmmn} the
choice $P_0=w_j$ is the choice of a root. Monodromy exchanges these branches. They all give the same
rational function $\lambda(z)$ in the root coordinate $z$ and thus the same continued superpotential.
A general linear combination of the $w_j$ does not select a single root.
\end{remark}

\subsection{\texorpdfstring{The extended connection of Gauss and Manin}{The extended connection of Gauss and Manin}}\label{sec:ext}

{
We now show that the translated functions are flat for a single connection on the complement of
$\Sigma$ in $M\times\mathbb C$. To describe this connection, we use its restrictions to the slices and
the vector field $V=e+\partial_{\lambda}$ along which the translated functions are constant.

\begin{proposition}\label{prop:ext}
{
There is a unique torsion free connection $\hat\nabla$ on $(M\times\mathbb C)\setminus\Sigma$ which
preserves the tangent bundle of each slice, restricts to $\Gamma(\lambda,\nu_0)$ on that slice, and
makes $V=e+\partial_\lambda$ parallel. This connection is flat. In a neighborhood of any regular point
$(x_0,\lambda_0)$ it is the pullback of the product connection
$\Gamma(\lambda_0,\nu_0)\oplus d_\lambda$ by the map
\[
\Phi_{\lambda_0}(x,\lambda)=(\phi_{\lambda_0-\lambda}(x),\lambda),
\]
where $d_\lambda$ denotes the trivial connection in the $\lambda$ direction.
}
\end{proposition}

\begin{proof}
{
Regard $e$ as a vector field on the product tangent to the slices of constant $\lambda$. For small $s$
the local flow of $V=e+\partial_\lambda$ is given by
\[
\Theta_s(x,\lambda)=(\phi_s(x),\lambda+s).
\]
Choose a local flat function $P_{\lambda_0}$ for $\Gamma(\lambda_0,\nu_0)$ on the regular part of the
reference slice $\lambda=\lambda_0$. We identify this regular part with an open subset of $M$ by the
map $(y,\lambda_0)\mapsto y$ and write the function as $P_{\lambda_0}(y)$, with the subscript
$\lambda_0$ indicating the reference slice. For $(x,\lambda)$ near $(x_0,\lambda_0)$, set
$y=\phi_{\lambda_0-\lambda}(x)$. Then $(y,\lambda_0)$ is the point where the flow of $V$ through
$(x,\lambda)$ meets the reference slice. By definition the translated function is
$P(x,\lambda)=P_{\lambda_0}(y)$. Using the composition law for the local flow, we obtain
\[
P(\Theta_s(x,\lambda))=P(\phi_s(x), \lambda+s)
=P_{\lambda_0}\bigl(\phi_{\lambda_0-\lambda-s}(\phi_s(x))\bigr)
=P_{\lambda_0}\bigl(\phi_{\lambda_0-\lambda}(x)\bigr)=P(x,\lambda).\]
Thus $P$ is constant along the flow of $V$. In particular, $V(P)=0$, which is the relation
$\partial_eP+\partial_\lambda P=0$ obtained in \eqref{eq:transl}. For the coordinate $\lambda$ we have
$V(\lambda)=1$.}

After shrinking the neighborhood of $(x_0,\lambda_0)$, we may assume that $\Phi_{\lambda_0}$ maps it
biholomorphically onto a product neighborhood. {The inverse sends $(y,\lambda)$ to
$(\phi_{\lambda-\lambda_0}(y),\lambda)$. With $y$ as above, the composition law for the flow gives
$\Phi_{\lambda_0}(\Theta_s(x,\lambda))=(y,\lambda+s)$. Differentiating with respect to $s$ yields
$(\Phi_{\lambda_0})_*V=\partial_\lambda$.} Define $\hat\nabla$ to be the pullback of
$\Gamma(\lambda_0,\nu_0)\oplus d_\lambda$ by $\Phi_{\lambda_0}$. By part (i) of Theorem~\ref{thm:ext},
its restriction to each slice is $\Gamma(\lambda,\nu_0)$. Moreover, $\partial_\lambda$ is parallel for
the product connection, so $V$ is parallel for its pullback.

{
Set $D_0=\Gamma(\lambda_0,\nu_0)$ and $D=D_0\oplus d_\lambda$. On the first factor of the product
neighborhood choose flat coordinates $v^1,\dots,v^n$ for $D_0$. Their coordinate fields are
independent of $\lambda$. In the coordinates $(v^1,\dots,v^n,\lambda)$ we have
\[
D_{\partial_{v^i}}\partial_{v^j}=0,\qquad
D_{\partial_{v^i}}\partial_\lambda
=D_{\partial_\lambda}\partial_{v^i}
=D_{\partial_\lambda}\partial_\lambda=0.
\]
The first equality follows from the choice of flat coordinates for $D_0$. For the remaining equalities
we use the definition of the product connection and the fact that $d_\lambda$ is the trivial
connection on the second factor.

The functions $\widetilde v^i(x,\lambda)
=v^i(\phi_{\lambda_0-\lambda}(x))$ together with $\lambda$
form a coordinate system on the source, since $\Phi_{\lambda_0}$ is a biholomorphism. The differential
of $\Phi_{\lambda_0}$ sends their coordinate fields to those of $(v^1,\dots,v^n,\lambda)$. Thus the
Christoffel symbols of the pullback connection $\hat\nabla$ vanish in
$(\widetilde v^1,\dots,\widetilde v^n,\lambda)$. Since the coordinate fields commute, both the torsion
and the curvature vanish in this frame.

To extend $P_{\lambda_0}$ from the first factor to the product, we compose it with the projection
$\operatorname{pr}_1(y,\lambda)=y$ and obtain
\[
F=P_{\lambda_0}\circ\operatorname{pr}_1,\qquad
F(y,\lambda)=P_{\lambda_0}(y).
\]
Thus $F$ is independent of $\lambda$ for fixed $y$. To verify $D\,dF=0$, use the coordinates
$(v^1,\dots,v^n,\lambda)$. The components of $D\,dF$ are then the second derivatives of $F$. Those
with two $v$ indices vanish by $D_0\,dP_{\lambda_0}=0$. Those involving $\lambda$ vanish because $F$
is independent of $\lambda$. In the same coordinates $D\,d\lambda=0$, since $\lambda$ is one of the
flat coordinates.

For any function $G$ on the product neighborhood we have
\[
\hat\nabla\,d(G\circ\Phi_{\lambda_0})
=\Phi_{\lambda_0}^{*}(D\,dG).
\]
Indeed, the components on both sides are the same second derivatives in the two corresponding systems
of coordinates. For $G=F$ we obtain
\[
(F\circ\Phi_{\lambda_0})(x,\lambda)
=F(\phi_{\lambda_0-\lambda}(x),\lambda)
=P_{\lambda_0}(\phi_{\lambda_0-\lambda}(x))=P(x,\lambda).
\]
It follows that $\hat\nabla\,dP=0$. For $G=\lambda$ we obtain $\hat\nabla\,d\lambda=0$, since
$\lambda\circ\Phi_{\lambda_0}=\lambda$.}

To prove uniqueness, regard vector fields $X,Y$ on $M$ as vector fields tangent to the slices and
independent of $\lambda$. Their covariant derivative is then prescribed by the condition
$\hat\nabla_XY=\Gamma(\lambda,\nu_0)_XY$. {By Remark~\ref{rem:KMS} we have
$\hat \nabla_XY=\nabla_XY+N(X*_\lambda Y)$, where $N$ is defined by \eqref{eq:aux2}. Taking the unit
$e$ of $\circ$ in the second argument of \eqref{eq:aux1} gives $X*_\lambda e=L_\lambda^{-1}X$. Since
$\nabla_Xe=0$, it follows that
\[
\hat\nabla_X e=NL_\lambda^{-1}X.
\]} The condition $\hat\nabla_XV=0$ now gives $\hat\nabla_X\partial_\lambda=-NL_\lambda^{-1}X$.
Moreover, $[\partial_\lambda,X]=0$, so torsion freeness implies that
$\hat\nabla_{\partial_\lambda}X=\hat\nabla_X\partial_\lambda$. {In particular, for $X=e$ we obtain
$\hat\nabla_{\partial_\lambda}e=-NL_\lambda^{-1}e$. The condition $\hat\nabla_{\partial_\lambda}V=0$
then gives $\hat\nabla_{\partial_\lambda}\partial_\lambda=NL_\lambda^{-1}e$.}

{
In product coordinates $(x^1,\dots,x^n,\lambda)$, denote by $a^m_{hk}$ the Christoffel symbols of
$\nabla$ with respect to $x$. The preceding identities give
\begin{align*}
\hat\Gamma^m_{hk}
&=\Gamma^m_{hk}(\lambda,\nu_0)
=a^m_{hk}+(NL_\lambda^{-1})^m_\ell c^\ell_{hk},\\
\hat\Gamma^m_{h\lambda}
&=\hat\Gamma^m_{\lambda h}=-(NL_\lambda^{-1})^m_h,\qquad
\hat\Gamma^m_{\lambda\lambda}=(NL_\lambda^{-1}e)^m,\\
\hat\Gamma^\lambda_{AB}&=0
\end{align*}
for all lower indices $A,B$, including $\lambda$. The last equality follows since all the covariant
derivatives of coordinate vector fields computed above are tangent to the slices. These coefficients
determine the connection uniquely.} Thus the local pullbacks agree on their common domains and define
the extended connection of Gauss and Manin on the full complement of $\Sigma$.
\end{proof}

\begin{corollary}\label{cor:spext}
Let $U$ be a neighborhood of $(x_0,\lambda_0)$ which the map $\Phi_{\lambda_0}$ of
Proposition~\ref{prop:ext} takes biholomorphically onto $U_0\times D$. Here $U_0$ lies in the regular
part of the reference slice and $D$ is a disk about $\lambda_0$. On $U$, the translated functions of
Theorem~\ref{thm:ext} are precisely the flat functions $F$ of $\hat\nabla$ satisfying $V(F)=0$. Such a
function admits a local inverse with respect to $\lambda$ wherever $\partial_\lambda F\ne0$. On the
domain of Definition~\ref{def:sp} the inverse of such a normalized flat function is the superpotential
of weight $\nu_0$.
\end{corollary}

\begin{proof}
{
{
Set $G=F\circ\Phi_{\lambda_0}^{-1}$. Since $(\Phi_{\lambda_0})_*V=\partial_\lambda$, the condition
$V(F)=0$ implies $\partial_\lambda G=0$. By connectedness of $D$, we have
$G(y,\lambda)=G(y,\lambda_0)$ for every $y\in U_0$ and $\lambda\in D$. The function
$P_{\lambda_0}(y)=G(y,\lambda_0)$ on the reference slice is flat for $\Gamma(\lambda_0,\nu_0)$ by
Proposition~\ref{prop:ext}, and we recover $F(x,\lambda)=P_{\lambda_0}(\phi_{\lambda_0-\lambda}(x))$.
Conversely, every such translate is flat for $\hat\nabla$ by the same proposition and satisfies
$V(F)=0$ by \eqref{eq:transl}.} Its derivative with respect to $\lambda$ is $-\partial_eP_{\lambda_0}$
evaluated at the translated point. Wherever the flow to the regular part of the zero slice is defined,
the function $P_0(y)=P_{\lambda_0}(\phi_{\lambda_0}(y))$ is flat for $\Gamma(0,\nu_0)$ by
Theorem~\ref{thm:ext}(i). Its translate is $F$. Inverting this translate gives the local
superpotential of Definition~\ref{def:sp}.
}
\end{proof}

Adding $a\lambda$ to $F$ changes $V(F)$ by $a$. Once we have imposed $V(F)=0$, the condition for local
invertibility is $\partial_\lambda F\ne0$.
}

\subsection{\texorpdfstring{Properties of the superpotential}{Properties of the superpotential}}\label{sec:char}

{
{
Consider a translated flat function $P$ as in Theorem~\ref{thm:ext}, and set $\xi=d_uP$. In flat coordinates
for $\nabla$, system~\eqref{eq:xisys} takes the form
\[
\partial_h\partial_kP=(\xi NL_\lambda^{-1})_l\,c^l_{hk},
\qquad
\partial_\lambda\partial_lP=-(\xi NL_\lambda^{-1})_l .
\]
The first equality is the condition of flatness on the slice at fixed $\lambda$. To obtain the second, we
differentiate the translation equation $\partial_\lambda P=-\partial_eP$ with respect to $u$ and evaluate the
first equality on the unit in one argument. If we eliminate $\xi NL_\lambda^{-1}$ between the two equalities,
we obtain \eqref{eq:aux5} for $P$. We now express this equation in terms of the local inverse $\Lambda$.}

\begin{lemma}\label{lem:inv}
{ Let $u$ be flat coordinates for $\nabla$, and suppose that $P(u,\lambda)$ is holomorphic near a point at
which $\partial_\lambda P\ne0$. Denote its local inverse with respect to $\lambda$ by $\Lambda(P,u)$, so that
$P$ is now the first independent variable of $\Lambda$. Then the equation
\begin{equation}\label{eq:aux5}
\partial_h\partial_kP=-c^l_{hk}\,\partial_\lambda\partial_lP,
\end{equation}
in which we differentiate $P$ with respect to $u$ at fixed $\lambda$, is equivalent to the following identity,
in which we differentiate $\Lambda$ with respect to $u$ at fixed $P$.}
\begin{equation}\label{eq:S3}
\partial_P\!\left[\frac{\partial_h\Lambda\,\partial_k\Lambda-c^l_{hk}\,\partial_l\Lambda}{\partial_P\Lambda}\right]
=\partial_h\partial_k\Lambda .
\end{equation}
\end{lemma}

\begin{proof}

We denote the derivatives of $\Lambda$ by subscripts and take derivatives with respect to $u$ at fixed $P$.
For a vector field $X$ on $M$, we use the notation $\Lambda_X=X\Lambda$ with the same convention. In composing
with the inverse we evaluate these derivatives at $P=P(u,\lambda)$. Differentiating
$\Lambda(P(u,\lambda),u)=\lambda$ with respect to $u^h$ and $\lambda$ gives
\[
\partial_hP=f_h(P(u,\lambda),u),\quad f_h:=-\frac{\Lambda_h}{\Lambda_P},\qquad \partial_\lambda P=\frac1{\Lambda_P}.
\]
Applying the chain rule once more, we obtain $\partial_k\partial_hP=\partial_kf_h+f_k\,\partial_Pf_h$ and
$\partial_\lambda\partial_lP=\partial_Pf_l\,\partial_\lambda P=\partial_Pf_l/\Lambda_P$. More explicitly, we
have
\[
\partial_kf_h=-\frac{\Lambda_{hk}}{\Lambda_P}+\frac{\Lambda_h\Lambda_{kP}}{\Lambda_P^2},\qquad
\partial_Pf_h=-\frac{\Lambda_{hP}}{\Lambda_P}+\frac{\Lambda_h\Lambda_{PP}}{\Lambda_P^2},\qquad
f_k\,\partial_Pf_h=\frac{\Lambda_k\Lambda_{hP}}{\Lambda_P^2}-\frac{\Lambda_h\Lambda_k\Lambda_{PP}}{\Lambda_P^3}.
\]
Thus
\[
\Lambda_P\big(\partial_kf_h+f_k\,\partial_Pf_h\big)
=-\Lambda_{hk}+\frac{\Lambda_h\Lambda_{kP}+\Lambda_k\Lambda_{hP}}{\Lambda_P}-\frac{\Lambda_h\Lambda_k\Lambda_{PP}}{\Lambda_P^2}
=-\Lambda_{hk}+\partial_P\Big(\frac{\Lambda_h\Lambda_k}{\Lambda_P}\Big),
\]
and the independence of $c^l_{hk}$ from $P$ gives
$\Lambda_P\,c^l_{hk}\,\partial_Pf_l/\Lambda_P=-\partial_P\big(c^l_{hk}\Lambda_l/\Lambda_P\big)$. Combining the
two expressions, we obtain
\[
\Lambda_P\big(\partial_h\partial_kP+c^l_{hk}\,\partial_\lambda\partial_lP\big)
=-\Lambda_{hk}+\partial_P\Big[\frac{\Lambda_h\Lambda_k-c^l_{hk}\Lambda_l}{\Lambda_P}\Big],
\]
where, since $\Lambda_P\neq0$, the left side vanishes precisely when \eqref{eq:S3} holds.
\end{proof}

{We can now characterize a superpotential by differential equations near any point $(u_0,\lambda_0)$ in the
complement of $\Sigma$.}
\begin{theorem}\label{thm:char}
Fix $\nu_0$ and $\kappa$. Let $\Lambda(P,u)$ be holomorphic near $(P_*,u_0)$, with
$\partial_P\Lambda(P_*,u_0)\ne0$, and denote its local inverse by $P(u,\lambda)$. Set
$\lambda_0=\Lambda(P_*,u_0)$ and assume that $(u_0,\lambda_0)\notin\Sigma$. { We say that $\Lambda$ is a local
superpotential of weight $\nu_0$ and constant $\kappa$ based at $\lambda_0$ if its inverse $P$ is obtained, by
the translation of Theorem~\ref{thm:ext}, from a flat function $P_{\lambda_0}$ for $\Gamma(\lambda_0,\nu_0)$
satisfying
\[
(E-\lambda_0e)P_{\lambda_0}=\nu_0P_{\lambda_0}+\kappa.
\]
Then $\Lambda$ is such a superpotential if and only if the following three conditions hold near $(P_*,u_0)$.}
\begin{itemize}
\item[(S1)] $\partial_e\Lambda=1$;
\item[(S2)] $E\Lambda+(\nu_0P+\kappa)\,\partial_P\Lambda=\Lambda$;
\item[(S3)] $\Lambda$ satisfies \eqref{eq:S3}.
\end{itemize}

Set $\xi=d_uP$ and let $N$ be as in \eqref{eq:aux2}. In terms of $P$ and its differential $\xi$, the three
conditions have the following consequences and equivalent formulations.
\begin{enumerate}
\item[(i)] (S1) $\iff$ $\partial_\lambda P=-\partial_eP$ $\iff$
$P(x,\lambda)=P_{\lambda_0}(\phi_{\lambda_0-\lambda}(x))$, where $P_{\lambda_0}=P(\cdot,\lambda_0)$;
\item[(ii)] (S2) $\iff$ $(E+\lambda\partial_\lambda)P=\nu_0P+\kappa$;
\item[(iii)] (S3) $\iff$ $\partial_h\partial_kP=-c^l_{hk}\,\partial_\lambda\partial_lP$;
\item[(iv)] if (S2) and (S3) hold, then $\partial_h\xi=\xi\,NL_\lambda^{-1}C_h$ and
$\partial_\lambda\xi=-\xi\,NL_\lambda^{-1}$. In particular, $P(\cdot,\lambda)$ is a flat function for
$\Gamma(\lambda,\nu_0)$ for every $\lambda$. Conversely, these two equations imply (S3) and imply that
$(E+\lambda\partial_\lambda)P-\nu_0P$ depends only on $\lambda$. This calculation is valid on any regular
inverse chart and uses only the product, the field $E$, and a flat torsion free reference connection.
\end{enumerate}
In particular, when $\lambda_0=0$, the flat function in Definition~\ref{def:sp} is $P_0=P(\cdot,0)$. More
generally, suppose that $\phi_{-\lambda_0}$ takes the reference slice to the regular part of the zero slice,
and define $P_0(y):=P_{\lambda_0}(\phi_{\lambda_0}(y))$ wherever the flow has this property. It follows from
part~(i) of Theorem~\ref{thm:ext} that this function is flat for $\Gamma(0,\nu_0)$. The local inverse obtained
by translating this function agrees with Definition~\ref{def:sp}.
\end{theorem}

\begin{proof}
{ Shrinking the neighborhoods if necessary, we may assume that $\Lambda_P$ is nonzero and that
$\Lambda(P(u,\lambda),u)=\lambda$. For any vector field $X$ on $M$, we differentiate $P$ at fixed $\lambda$
and $\Lambda$ at fixed first argument. On the inverse chart we evaluate the derivatives of $\Lambda$ at
$(P(u,\lambda),u)$. With these conventions the chain rule gives
\[
\Lambda_P\,\partial_\lambda P=1,\qquad
X\Lambda+\Lambda_P\,X(P)=0.
\]

To prove (i), we take $X=e$ in the preceding identities and obtain
\[
(\partial_\lambda+\partial_e)P=\frac{1-\partial_e\Lambda}{\Lambda_P}.
\]
Thus (S1) is equivalent to $V(P)=0$, where $V=e+\partial_\lambda$. The local flow of $V$ takes $(x,\lambda)$
to $(\phi_t(x),\lambda+t)$, and every function satisfying $V(P)=0$ is constant along its trajectories. In
particular, taking $t=\lambda_0-\lambda$ gives
\[
P(x,\lambda)=P\bigl(\phi_{\lambda_0-\lambda}(x),\lambda_0\bigr)
=P_{\lambda_0}\bigl(\phi_{\lambda_0-\lambda}(x)\bigr),
\qquad P_{\lambda_0}=P(\cdot,\lambda_0).
\]
Conversely, for a function $P$ of this form the flow law gives $P(\phi_t(x),\lambda+t)=P(x,\lambda)$.
Differentiating at $t=0$, we obtain $V(P)=0$, which implies (S1).

To prove (ii), define $H:=(E+\lambda\partial_\lambda)P-\nu_0P$. Taking $X=E$ in the chain rule gives
$E(P)=-E\Lambda/\Lambda_P$. Since $\lambda=\Lambda(P(u,\lambda),u)$, it follows that
\[
H-\kappa
=\frac{\Lambda-E\Lambda-(\nu_0P+\kappa)\Lambda_P}{\Lambda_P}.
\]
Since $(u,\lambda)\mapsto(P(u,\lambda),u)$ is a local biholomorphism, $H=\kappa$ holds in the $(u,\lambda)$
chart if and only if the numerator vanishes in the $(P,u)$ chart. This is precisely condition (S2). Part (iii)
follows from Lemma~\ref{lem:inv}. Here we take derivatives with respect to $u$ at fixed $\lambda$ for $P$ and
at fixed $P$ for $\Lambda$.

For (iv), choose flat coordinates $u$ and write $\xi_k=\partial_kP$ and $Q^m_k=\partial_k E^m$. In this
notation, (iii) reads
\[
\partial_h\xi=-\partial_\lambda\xi\,C_h .
\]
Differentiating $H=E^m\xi_m+\lambda\partial_\lambda P-\nu_0 P$ with
respect to $u^k$ at fixed $\lambda$, we obtain
\[
\partial_k H
=E^m\partial_k\xi_m+\xi_mQ^m_k
+\lambda\partial_\lambda\xi_k-\nu_0\xi_k .
\]
By (iii) and the identity $L^l_k=E^m c^l_{km}$, the first term on the right is
$-E^m c^l_{km}\partial_\lambda\xi_l=-(\partial_\lambda\xi\,L)_k$. Thus
\[
d_uH=-\partial_\lambda\xi\,L_\lambda-\xi N.
\]
Under (S2) we have $H=\kappa$ and therefore $d_uH=0$. Since $L_\lambda$ is invertible on the regular chart, we
obtain $\partial_\lambda\xi=-\xi NL_\lambda^{-1}$. Substituting this expression into the equation for
$\partial_h\xi$ gives $\partial_h\xi=\xi NL_\lambda^{-1}C_h$. For each fixed $\lambda$, this is the flatness
equation $\Gamma(\lambda,\nu_0)\,d_uP=0$ by \eqref{eq:op}. Conversely, the two equations for $\xi$ imply
$\partial_h\xi=-\partial_\lambda\xi\,C_h$, which is (iii). Using the preceding expression for $d_uH$, we
obtain $d_uH=0$. Thus $H$ is locally a function of $\lambda$ alone.

Suppose now that $\Lambda$ is the inverse of a translated flat function with the given weight and constant. By
Theorem~\ref{thm:ext} we have $V(P)=0$, $H=\kappa$, and the system \eqref{eq:xisys}. Eliminating
$\xi NL_\lambda^{-1}$ between the two equations of this system gives (iii). Therefore parts (i), (ii), and
(iii) give (S1), (S2), and (S3).

Conversely, assume (S1), (S2), and (S3). By part (iv), $P_{\lambda_0}=P(\cdot,\lambda_0)$ is flat for
$\Gamma(\lambda_0,\nu_0)$, and by part (i), $P$ is its translate. Moreover,
\[
\partial_eP_{\lambda_0}(u_0)
=-\partial_\lambda P(u_0,\lambda_0)
=-\frac{1}{\Lambda_P(P_*,u_0)}\ne0.
\]
On evaluating $H=\kappa$ at $\lambda=\lambda_0$ and using $\partial_\lambda P=-\partial_eP$, we find
\[
(E-\lambda_0e)P_{\lambda_0}=\nu_0P_{\lambda_0}+\kappa.
\]
Thus the reference flat function has the required constant, and $\Lambda$ is a local superpotential of weight
$\nu_0$ and constant $\kappa$ based at $\lambda_0$.

Suppose that the flow takes the reference slice to the regular part of the zero slice, and set
$P_0=P_{\lambda_0}\circ\phi_{\lambda_0}$ wherever this holds. By the identity for the connections in part~(i)
of Theorem~\ref{thm:ext}, $P_0$ is flat for $\Gamma(0,\nu_0)$. Moreover, $[e,E]=e$ and invariance of $e$ under
the flow give $\phi_{\lambda_0}^*(E-\lambda_0 e)=E$. It follows that
\[
E(P_0)
=\bigl((E-\lambda_0e)P_{\lambda_0}\bigr)\circ\phi_{\lambda_0}
=\nu_0P_0+\kappa.
\]
Invariance of $e$ also gives $\partial_eP_0=(\partial_eP_{\lambda_0})\circ\phi_{\lambda_0}$, which is nonzero
at $\phi_{-\lambda_0}(u_0)$. Finally, the flow law gives
$P_0(\phi_{-\lambda}(x))
=P_{\lambda_0}(\phi_{\lambda_0-\lambda}(x))=P(x,\lambda)$ wherever both sides are
defined. Thus the inverse of $P$ agrees with the local superpotential in Definition~\ref{def:sp}.}
\end{proof}

{ We have written condition \eqref{eq:S3} in flat coordinates for $\nabla$. To express it in arbitrary
coordinates, fix $P$ and consider the covariant Hessian of the function $\Lambda(P,\cdot)$ on $M$. By
definition,
\[
\nabla^2\Lambda(X,Y):=X(Y\Lambda)-(\nabla_XY)\Lambda,\qquad\text{i.e.}\qquad
\nabla^2\Lambda(\partial_h,\partial_k)=\Lambda_{hk}-\Gamma^m_{hk}\Lambda_m
\]
in arbitrary coordinates. Since $\nabla$ is torsion free, the Hessian is symmetric. In flat coordinates
$\Gamma^m_{hk}=0$, so the Hessian reduces to $\partial_h\partial_k\Lambda$, which appears on the right side of
\eqref{eq:S3}. To give a tensorial interpretation of the left side, define
\[
R_\Lambda(X,Y):=X\Lambda\cdot Y\Lambda-(X\circ Y)\Lambda.
\]
This expression is symmetric in $X,Y$ and bilinear over $C^\infty(M)$. Dividing by $\partial_P\Lambda$ and
differentiating with respect to $P$ both commute with multiplication by functions on $M$. Thus (S3) is an
identity of tensors. To write it in local coordinates $x^1,\dots,x^n$, denote by $\Gamma^m_{hk}$ the
Christoffel symbols of $\nabla$ and use subscripts for derivatives with respect to $x^h$ at fixed $P$. In
these coordinates (S3) reads
\begin{equation}\label{eq:S3cov}
\Lambda_{hk}-\Gamma^m_{hk}\,\Lambda_m=\partial_P\Big[\frac{\Lambda_h\Lambda_k-c^m_{hk}\Lambda_m}{\Lambda_P}\Big] ,
\end{equation}
or equivalently $\nabla^2\Lambda(X,Y)=\partial_P\big[R_\Lambda(X,Y)/\partial_P\Lambda\big]$ for all vector
fields $X,Y$. In particular, $R_\Lambda$ vanishes precisely where the map $X\mapsto X\Lambda$ preserves
multiplication. In the proof of Theorem~\ref{thm:main} we use its vanishing at a simple critical point to show
that the differential of the critical value preserves multiplication. We shall also use \eqref{eq:S3cov} for
an arbitrary torsion free connection in Proposition~\ref{prop:christoffel} and Theorem~\ref{thm:flat}.}

Let us examine to what extent the three conditions determine the inverse function $P$.
\begin{remark}\label{rem:indep}
Once the weight is fixed, (S2) and (S3) determine $d_uP$ from its value on a single slice. In particular, two
solutions whose differentials agree on one slice differ locally by a function $f(\lambda)$ on a product
neighborhood with connected slices. If the constants in (S2) differ by $\Delta\kappa$, this function satisfies
$\lambda f'-\nu_0f=\Delta\kappa$. On a branch away from $\lambda=0$ we obtain
$f=c\lambda^{\nu_0}-\Delta\kappa/\nu_0$ for $\nu_0\ne0$ and $f=c+\Delta\kappa\log\lambda$ for $\nu_0=0$. At
weight zero, holomorphy at $\lambda=0$ implies $\Delta\kappa=0$, so $f$ is constant. Condition (S1) also
requires $f$ to be constant. The calculation in part~(iv) of Theorem~\ref{thm:char} involves only the product,
the flat coordinates of the reference connection, and the Euler field. We may therefore test (S1), (S2), and
(S3) on an explicit function before checking compatibility with the reference connection.
\end{remark}

Recall that a \emph{character} of a unital complex algebra $A$ is a complex linear map
$\chi\colon A\to\mathbb C$ that preserves multiplication and maps the unit of $A$ to $1$. Here we consider
$A=T_uM$ with the product $\circ$ and the unit $e_u$.

We next show how (S1) and (S3) give a character of the tangent algebra at a simple critical point, namely the
differential of the corresponding critical value.
\begin{theorem}\label{thm:main}
Let $\Lambda$ be a superpotential which extends holomorphically to neighborhoods of $n$ simple critical points
$P_*^1(u),\dots,P_*^n(u)$. These points may lie on different sheets of the continuation. Suppose that (S1),
(S2), and (S3) hold on the regular parts of these neighborhoods and that the critical values
$r^i=\Lambda(P_*^i(u),u)$ have distinct differentials. Then the tangent algebra is semisimple on the base
neighborhood, and the functions $r^i$ are canonical coordinates satisfying $E(r^i)=r^i$. Moreover,
\[
\partial_{r^j}\Lambda\,\big(P_*^i\big)=\delta_{ij},
\]
that is, at $P_*^i$ the dependence of $\Lambda$ on the canonical coordinates is, to first order, only through
$r^i$.
\end{theorem}

\begin{proof}

Since $\Lambda_{PP}\ne0$ at each critical point, the holomorphic implicit function theorem gives the branches
$P_*^i(u)$. Fix one of them and set $P_*=P_*^i(u)$. In (S3) write
$R_{hk}=\Lambda_h\Lambda_k-c^l_{hk}\Lambda_l$ and $B_{hk}=R_{hk}/\Lambda_P$. The only possible principal part
of $B_{hk}$ is $b_{-1}/(P-P_*)$. Now $\partial_P B_{hk}=\Lambda_{hk}$ is holomorphic, so the term
$-b_{-1}/(P-P_*)^2$ in the derivative must vanish. Therefore $b_{-1}=0$. Thus $R_{hk}=B_{hk}\Lambda_P$
vanishes at $P_*$ and we obtain
\begin{equation}\label{eq:hom}
\Lambda_h\Lambda_k=c^l_{hk}\Lambda_l\qquad\text{at }P=P_*,\quad\text{i.e.}\quad (X\Lambda)(Y\Lambda)=\big((X\circ Y)\Lambda\big)\ \text{at }P_* .
\end{equation}

For a vector field $X=X^h\partial_h$ on $M$, denote by $X\Lambda$ the derivative at fixed $P$. The chain rule
gives
\[
dr^i(X)=X^h\,\partial_h\big(\Lambda(P_*^i(u),u)\big)=\Lambda_P\big(P_*^i\big)\,dP_*^i(X)+(X\Lambda)\big(P_*^i\big)=(X\Lambda)\big(P_*^i\big).
\]

It follows from \eqref{eq:hom} that $dr^i(X\circ Y)=dr^i(X)\,dr^i(Y)$, and (S1) gives $dr^i(e)=1$. Thus the
$dr^i$ are characters of the tangent algebra. By Dedekind's lemma, distinct characters of an algebra are
linearly independent. Therefore the $dr^i$ are independent. It follows that $(r^1,\dots,r^n)$ are local
coordinates and that the characters identify the tangent algebra at each point with $\mathbb C^n$. In these
coordinates we have $dr^j(\partial_{r^i}\circ\partial_{r^k})=\delta_{ji}\delta_{jk}$. It follows that
$\partial_{r^i}\circ\partial_{r^k}=\delta_{ik}\partial_{r^i}$, as required for canonical coordinates.
Evaluating (S2) at $P_*^i$ gives $E(r^i)=r^i$. Taking $X=\partial_{r^j}$ in the preceding formula for $dr^i$
then gives $(\partial_{r^j}\Lambda)(P_*^i)=dr^i(\partial_{r^j})=\delta_{ij}$.
\end{proof}

{
\begin{remark}\label{rem:simplecrit}
For any simple critical point $P_*$ of $\Lambda$, the first part of the proof shows that the differential of
the critical value $r=\Lambda(P_*)$ is a character of the tangent algebra. In the semisimple case the
characters are the differentials of the canonical coordinates, so $dr=du^k$ for some canonical coordinate
$u^k$. Moreover, $E(r)=r$ by (S2) and $E(u^k)=u^k$. Evaluating $dr=du^k$ on $E$ therefore gives $r=u^k$. Thus, in the semisimple case,
the value of $\Lambda$ at every simple critical point is a canonical coordinate, and
$(\partial_{u^j}\Lambda)(P_*)=\delta_{jk}$. Several critical points can give the same coordinate. This happens
at $p=\pm\sqrt{-a_1/2}$ for the polynomial $p^4+a_1p^2+a_2$ in Remark~\ref{rem:lame} and at the $m$ lifts of
each critical point in Section~\ref{sec:Gmmn}. The hypothesis of Theorem~\ref{thm:main} requires that all $n$
canonical coordinates occur as simple critical values. Suppose, in particular, that $n$ simple critical values
$r^1,\dots,r^n$ are distinct as functions on a connected neighborhood. Their differentials are then distinct
on a dense open subset. Indeed, if $dr^i=dr^j$ on an open subset, evaluation on $E$ gives $r^i=r^j$ there, and
holomorphy extends the equality to the whole neighborhood. On the dense open subset where the differentials
are distinct, their independence follows from Dedekind's lemma and Theorem~\ref{thm:main} applies.
\end{remark}
}

Under the hypotheses of Theorem~\ref{thm:main}, the critical values $r^1,\dots,r^n$ give canonical
coordinates. Thus $\partial_i\circ\partial_j=\delta_{ij}\partial_i$ and $e=\sum_i\partial_i$, where
$\partial_i=\partial/\partial r^i$. We shall use these identities to write the natural connection and
condition (S3) in the coordinates $r^i$.

For a function $f(P,r)$ we use the notation
\[
f\big|_i:=f\big(P_*^i(r),r\big)
\]
for the value of $f$ at $P_*^i$, which depends only on $r$. The subscripts $j,k$ denote differentiation with
respect to $r^j,r^k$ \emph{at fixed $P$}, before evaluation. For instance, $\Lambda_{jk}|_i$ stands for
$(\partial_j\partial_k\Lambda)(P_*^i,r)$. If we evaluate first, the chain rule gives
\begin{equation}\label{eq:evalder}
\partial_j\big(f|_i\big)=f_j\big|_i+f_P\big|_i\;\partial_jP_*^i .
\end{equation}
In particular, $\Lambda|_i=r^i$, $\Lambda_P|_i=0$, and $\Lambda_{PP}|_i\neq0$. Moreover,
Theorem~\ref{thm:main} gives $\Lambda_j|_i=\delta_{ij}$.

We recall the form of the connection in canonical coordinates.

\begin{lemma}\label{lem:canform}
Let $(M,\nabla,\circ,e)$ be a semisimple flat F-manifold with canonical coordinates $r^1,\dots,r^n$. The
Christoffel symbols of $\nabla$ in these coordinates satisfy
\[
\Gamma^i_{jk}=0\ \ (i,j,k\ \text{distinct}),\qquad \Gamma^i_{jj}=-\Gamma^i_{ij}\ \ (i\neq j),\qquad
\Gamma^i_{ii}=-\sum_{l\neq i}\Gamma^i_{il} .
\]
In particular, the $n(n-1)$ functions $\Gamma^i_{ij}$, with $i\neq j$, determine $\nabla$.
\end{lemma}
A proof can be found in \cite{ALDE}.

\begin{proposition}\label{prop:S3can}
Under the hypotheses of Theorem~\ref{thm:main}, use canonical coordinates and denote by subscripts the
derivatives with respect to $r^i$ at fixed $P$.
\begin{enumerate}
\item[(i)] In the presence of (S1), condition (S3) is equivalent to the $\binom n2$ equations
\begin{equation}\label{eq:S3can}
\Lambda_{ij}-\Gamma^i_{ij}\,\Lambda_i-\Gamma^j_{ij}\,\Lambda_j=\partial_P\Big(\frac{\Lambda_i\Lambda_j}{\Lambda_P}\Big),\qquad i\neq j .
\end{equation}
\item[(ii)] The translated flat functions $P(\cdot,\lambda)$ of Theorem~\ref{thm:ext} satisfy, for every
$\lambda$, the linear system
\begin{equation}\label{eq:Pcan}
\partial_i\partial_jP=\Gamma^i_{ij}\,\partial_iP+\Gamma^j_{ij}\,\partial_jP,\qquad i\neq j ,
\end{equation}
which is independent of $\lambda$. Only the diagonal equations
$\nabla^2P(\partial_i,\partial_i)=-\partial_i\partial_\lambda P$ involve the parameter $\lambda$.
\end{enumerate}
\end{proposition}

\begin{proof}
(i) In canonical coordinates, take $h=i$ and $k=j$ in \eqref{eq:S3cov}. When $i\neq j$, we have $c^m_{ij}=0$,
so the numerator reduces to $\Lambda_i\Lambda_j$. By Lemma~\ref{lem:canform}, the only Christoffel symbols
$\Gamma^m_{ij}$ that can be nonzero are $\Gamma^i_{ij}$ and $\Gamma^j_{ij}$. This gives \eqref{eq:S3can}. It
remains to see why the diagonal equations $(i,i)$ follow. Both sides of \eqref{eq:S3cov}, written in the form
$\nabla^2\Lambda(X,Y)=\partial_P\big[R_\Lambda(X,Y)/\partial_P\Lambda\big]$, vanish for $Y=e$. Indeed, (S1)
and $\nabla e=0$ give $\nabla^2\Lambda(X,e)=X(e\Lambda)-(\nabla_Xe)\Lambda=0$, and
$R_\Lambda(X,e)=X\Lambda-X\Lambda=0$. Since $e=\sum_j\partial_j$, the equation $(i,i)$ is the negative of the
sum of the equations $(i,j)$ over $j\neq i$.

(ii) In flat coordinates, equation~\eqref{eq:aux5} in part~(iii) of Theorem~\ref{thm:char} reads
\[
\partial_h\partial_kP=-c^l_{hk}\,\partial_l(\partial_\lambda P).
\]
The left side is the Hessian of $P(\cdot,\lambda)$, whereas the right side is
$-(\partial_h\circ\partial_k)(\partial_\lambda P)$. Both expressions are tensorial, so
$\nabla^2P(X,Y)=-(X\circ Y)(\partial_\lambda P)$ in any coordinates. Now use canonical coordinates and take
$X=\partial_i$ and $Y=\partial_j$. For $i\neq j$ the right side vanishes, and Lemma~\ref{lem:canform} gives
\eqref{eq:Pcan}. For $i=j$ the right side reduces to $-\partial_i\partial_\lambda P$.
\end{proof}

\subsection{Reconstruction from a superpotential}\label{subsec:reconstruction}
In Dubrovin's construction the residue formula determines the Frobenius metric and therefore the unique
torsion free connection compatible with it. Under the hypotheses of Theorem~\ref{thm:main}, we can instead
evaluate (S3) at the chosen critical points to express the Christoffel symbols of $\nabla$ in terms of
$\Lambda$ and the product alone. We shall also use this equation to construct a connection from a given
function $\Lambda$, as in Theorem~\ref{thm:flat}.

\begin{proposition}\label{prop:christoffel}
Let $(M,\nabla,\circ,e,E)$ be a bi-flat F-manifold. Let $\Lambda$ be a superpotential of weight $\nu_0$ in the
sense of Definition~\ref{def:sp}, or a continuation of such a superpotential, on a domain where (S1), (S2),
and (S3) hold as in Theorem~\ref{thm:char}. In arbitrary coordinates $x^1,\dots,x^n$ on $M$, denote by
$\Gamma^m_{hk}$ the Christoffel symbols of $\nabla$. Subscripts denote differentiation with respect to $x^h$
at fixed $P$.
\begin{enumerate}
\item[(i)] The following identity holds as an equality of functions of $P$.
\begin{equation}\label{eq:chr}
\Gamma^m_{hk}\,\Lambda_m=\Lambda_{hk}-\partial_PB_{hk},\qquad B_{hk}:=\frac{\Lambda_h\Lambda_k-c^m_{hk}\Lambda_m}{\partial_P\Lambda} .
\end{equation}
\item[(ii)] Suppose in addition that $\Lambda$ satisfies the hypotheses of Theorem~\ref{thm:main}, with simple
critical points $P_*^1,\dots,P_*^n$. In the canonical coordinates $r^i=\Lambda(P_*^i)$ and with the notation
$f|_i=f(P_*^i)$, we have, for $i\neq j$,
\begin{equation}\label{eq:Gammacan}
\Gamma^i_{ij}=\Big[\Lambda_{ij}-\partial_P\Big(\frac{\Lambda_i\Lambda_j}{\Lambda_P}\Big)\Big]\Big|_i
=\Big[\Lambda_{ij}-\frac{\Lambda_{iP}\Lambda_{jP}}{\Lambda_{PP}}-\frac{\Lambda_{jPP}}{2\Lambda_{PP}}+\frac{\Lambda_{jP}\Lambda_{PPP}}{2\Lambda_{PP}^2}\Big]\Big|_i ,
\end{equation}
and Lemma~\ref{lem:canform} gives all the Christoffel symbols of $\nabla$ in terms of these functions.
\item[(iii)] Under the hypotheses of (ii), define the nondegenerate symmetric bilinear form
\[
\eta^\Lambda:=\sum_i\frac{dr^i\otimes dr^i}{\Lambda_{PP}(P_*^i)}.
\]
Its components are $\eta^\Lambda_{lm}=\sum_i\res_{P=P_*^i}\frac{\Lambda_l\Lambda_m}{\Lambda_P}\,dP$, and they
satisfy
\begin{equation}\label{eq:chrres}
\eta^\Lambda_{lm}\,\Gamma^m_{hk}=\sum_{i=1}^n\res_{P=P_*^i}\frac{\big(\Lambda_{hk}-\partial_PB_{hk}\big)\,\Lambda_l}{\Lambda_P}\,dP .
\end{equation}
\end{enumerate}
\end{proposition}

\begin{proof}
(i) Solving \eqref{eq:S3cov} in the coordinates $x$ for the term containing the Christoffel symbols gives the
stated identity between functions of $(P,x)$. Fix $x$ and regard the right side as a function of $P$. The left
side expresses it as a linear combination of the $n$ functions $\Lambda_1,\dots,\Lambda_n$ of $P$, with
coefficients $\Gamma^m_{hk}(x)$ independent of $P$. If these $n$ functions are linearly independent, the
coefficients are unique. We recover them by applying $n$ linear functionals on functions of $P$ dual to
$\Lambda_1,\dots,\Lambda_n$. In parts (ii) and (iii) we evaluate at the critical points. In
Remark~\ref{rem:chrBn} we take residues at the poles.

(ii) With the notation $f|_i$ introduced before Lemma~\ref{lem:canform}, let $i\neq j$ and set
$B_{ij}:=\Lambda_i\Lambda_j/\Lambda_P$. As we saw in the proof of Theorem~\ref{thm:main}, $B_{ij}$ is regular
at $P_*^i$. We may therefore evaluate \eqref{eq:S3can} at this point. Since $\Lambda_i|_i=1$ and
$\Lambda_j|_i=0$, its left side reduces to $\Lambda_{ij}|_i-\Gamma^i_{ij}$, giving
\[
\Gamma^i_{ij}=\Lambda_{ij}\big|_i-\partial_PB_{ij}\big|_i .
\]
This proves the first equality. To compute $\partial_PB_{ij}|_i$, set $R:=\Lambda_i\Lambda_j=B_{ij}\Lambda_P$
and differentiate once and twice with respect to $P$. We obtain
\[
R_P=\partial_PB_{ij}\,\Lambda_P+B_{ij}\,\Lambda_{PP},\qquad
R_{PP}=\partial_P^2B_{ij}\,\Lambda_P+2\,\partial_PB_{ij}\,\Lambda_{PP}+B_{ij}\,\Lambda_{PPP} .
\]
Since $\Lambda_P|_i=0$ and $\Lambda_{PP}|_i\neq0$, evaluating both identities at $P_*^i$ gives
\[
R_P\big|_i=B_{ij}\big|_i\,\Lambda_{PP}\big|_i,\qquad
R_{PP}\big|_i=2\,\partial_PB_{ij}\big|_i\,\Lambda_{PP}\big|_i+B_{ij}\big|_i\,\Lambda_{PPP}\big|_i .
\]
We can also use $R=\Lambda_i\Lambda_j$ with $\Lambda_i|_i=1$ and $\Lambda_j|_i=0$ to write
\[
R_P\big|_i=\Lambda_{jP}\big|_i,\qquad
R_{PP}\big|_i=\big(2\Lambda_{iP}\Lambda_{jP}+\Lambda_{jPP}\big)\big|_i .
\]
Comparing these expressions gives $B_{ij}|_i=\Lambda_{jP}|_i/\Lambda_{PP}|_i$ and
\[
\partial_PB_{ij}\big|_i=\frac{R_{PP}|_i-B_{ij}|_i\,\Lambda_{PPP}|_i}{2\Lambda_{PP}|_i}
=\Big[\frac{\Lambda_{iP}\Lambda_{jP}}{\Lambda_{PP}}+\frac{\Lambda_{jPP}}{2\Lambda_{PP}}-\frac{\Lambda_{jP}\Lambda_{PPP}}{2\Lambda_{PP}^2}\Big]\Big|_i ,
\]
which proves the second equality.

(iii) Near $P_*^i$ we have $\Lambda_P=\Lambda_{PP}|_i\,(P-P_*^i)+O\big((P-P_*^i)^2\big)$ with
$\Lambda_{PP}|_i\neq0$. Therefore, for any function $f$ regular at $P_*^i$, we have
\[
\res_{P=P_*^i}\frac{f}{\Lambda_P}\,dP=\frac{f|_i}{\Lambda_{PP}|_i} .
\]
Taking $f=\Lambda_l\Lambda_m$ and using $\Lambda_l|_i=\partial_lr^i$ as in the proof of
Theorem~\ref{thm:main}, we obtain the formula for $\eta^\Lambda_{lm}$ in terms of residues. In matrix
notation, $\eta^\Lambda=J^{\mathsf T}DJ$. Both $J=(\partial_lr^i)$ and
$D=\mathrm{diag}\big(1/\Lambda_{PP}|_i\big)$ are invertible, and the latter matrix is diagonal. Thus
$\eta^\Lambda$ is symmetric and nondegenerate.

In arbitrary coordinates, set $G:=\Lambda_{hk}-\partial_PB_{hk}$. The proof of Theorem~\ref{thm:main} shows
that $B_{hk}$ is regular at $P_*^i$, so $G$ is regular there as well. Part (i) gives
$G|_i=\Gamma^m_{hk}\,\Lambda_m|_i=\Gamma^m_{hk}\,\partial_mr^i$. Applying the residue formula with
$f=G\,\Lambda_l$, we obtain
\[
\sum_i\res_{P=P_*^i}\frac{G\,\Lambda_l}{\Lambda_P}\,dP=\sum_i\frac{G|_i\,\partial_lr^i}{\Lambda_{PP}|_i}
=\sum_i\frac{\Gamma^m_{hk}\,\partial_mr^i\,\partial_lr^i}{\Lambda_{PP}|_i}=\eta^\Lambda_{lm}\,\Gamma^m_{hk}.
\]
This proves \eqref{eq:chrres}. Since $\eta^\Lambda$ is nondegenerate, this identity determines the Christoffel
symbols $\Gamma^m_{hk}$.
\end{proof}

Let us express the natural connection in terms of the residue metric and compute the Lie derivatives of this
metric along the unit and Euler fields.
\begin{proposition}\label{prop:riem}
Suppose that the hypotheses of Theorem~\ref{thm:main} hold, and in the canonical coordinates $r^1,\dots,r^n$
write
\[
\eta^\Lambda=\sum_iH_i^2\,(dr^i)^2,\qquad H_i^2=\frac1{\Lambda_{PP}|_i} .
\]
Then the following assertions hold.
\begin{enumerate}
\item[(i)] The Christoffel symbols of $\nabla$ satisfy $\Gamma^i_{ij}=\partial_j\log H_i$ for $i\neq j$. It
follows from Lemma~\ref{lem:canform} that $\Gamma^i_{jj}=-\partial_j\log H_i$ for $i\neq j$, that
$\Gamma^i_{ii}=-\sum_{l\neq i}\partial_l\log H_i$, and that all the remaining symbols vanish. In particular,
the Lam\'e coefficients $H_i$ determine the natural connection.
\item[(ii)] The unit vector field is Killing, that is, $\mathrm{Lie}_e\eta^\Lambda=0$.
\item[(iii)] $\mathrm{Lie}_E\eta^\Lambda=(2\nu_0+1)\,\eta^\Lambda$.
\end{enumerate}
Thus $\eta^\Lambda$ is one of the metrics which, together with the product, define the Riemannian F-manifold
with Killing unit vector field associated with $(M,\nabla,\circ,e)$ in \cite{ABLR}.
\end{proposition}

\begin{proof}
We use the notation $f|_i$ introduced before Lemma~\ref{lem:canform} and the differentiation rule
\eqref{eq:evalder}. Applying \eqref{eq:evalder} to the two identities $\Lambda_P|_i=0$ and
$\Lambda_j|_i=\delta_{ij}$ from Theorem~\ref{thm:main} and then to the function $\Lambda_{PP}|_i$, we find the
following expressions, since the derivatives of $0$ and $\delta_{ij}$ vanish.
\begin{equation}\label{eq:critder}
\begin{aligned}
\partial_jP_*^i&=-\frac{\Lambda_{jP}|_i}{\Lambda_{PP}|_i},\\
\Lambda_{jk}\big|_i&=-\Lambda_{jP}\big|_i\;\partial_kP_*^i=\frac{\Lambda_{jP}|_i\,\Lambda_{kP}|_i}{\Lambda_{PP}|_i},\\
\partial_j\big(\Lambda_{PP}|_i\big)&=\Lambda_{jPP}\big|_i-\frac{\Lambda_{PPP}|_i\,\Lambda_{jP}|_i}{\Lambda_{PP}|_i} .
\end{aligned}
\end{equation}
Indeed, the first equality follows from
$0=\partial_j(\Lambda_P|_i)=\Lambda_{Pj}|_i+\Lambda_{PP}|_i\,\partial_jP_*^i$. For the second, use
$0=\partial_k(\Lambda_j|_i)=\Lambda_{jk}|_i+\Lambda_{jP}|_i\,\partial_kP_*^i$ and substitute the first
equality.

(i) Suppose $i\neq j$. Taking $k=i$ in the second identity of \eqref{eq:critder}, we see that the first two
terms in the second expression in \eqref{eq:Gammacan} cancel. The third identity of \eqref{eq:critder} then
gives
\[
\Gamma^i_{ij}=-\frac1{2\Lambda_{PP}|_i}\Big(\Lambda_{jPP}\big|_i-\frac{\Lambda_{PPP}|_i\,\Lambda_{jP}|_i}{\Lambda_{PP}|_i}\Big)
=-\frac{\partial_j\big(\Lambda_{PP}|_i\big)}{2\Lambda_{PP}|_i}=-\frac12\,\partial_j\log\Lambda_{PP}\big|_i ,
\]
which equals $\partial_j\log H_i$. Lemma~\ref{lem:canform} gives the remaining Christoffel symbols.

(ii) Condition (S1) gives $\sum_j\Lambda_j=1$ for every $P$. Differentiating this identity, we obtain
$\sum_j\Lambda_{jP}=\sum_j\Lambda_{jPP}=0$, in particular at $P_*^i$. It follows from the third identity in
\eqref{eq:critder} that $e\big(\Lambda_{PP}|_i\big)=\sum_j\partial_j\big(\Lambda_{PP}|_i\big)=0$. Since the
components of $e=\sum_j\partial_{r^j}$ are constant, we obtain
$\mathrm{Lie}_e\eta^\Lambda=\sum_ie\big(H_i^2\big)(dr^i)^2=0$.

(iii) {Differentiating (S2) once and twice with respect to $P$, we find
\[
E\Lambda_P+(\nu_0P+\kappa)\Lambda_{PP}=(1-\nu_0)\Lambda_P,\qquad
E\Lambda_{PP}+(\nu_0P+\kappa)\Lambda_{PPP}=(1-2\nu_0)\Lambda_{PP}.
\]
By \eqref{eq:evalder} we have $0=E(\Lambda_P|_i)=E\Lambda_P|_i+\Lambda_{PP}|_i\,E(P_*^i)$. The first identity
therefore gives $E(P_*^i)=\nu_0P_*^i+\kappa$. Applying \eqref{eq:evalder} once more and using the second
identity, we obtain $E(\Lambda_{PP}|_i)=E\Lambda_{PP}|_i+\Lambda_{PPP}|_i\,E(P_*^i)=(1-2\nu_0)\Lambda_{PP}|_i$
and consequently $E(H_i^2)=(2\nu_0-1)H_i^2$.} Since $E=\sum_ir^i\partial_{r^i}$, we have
$\mathrm{Lie}_E(dr^i)^2=2(dr^i)^2$, and it follows that $\mathrm{Lie}_E\eta^\Lambda=(2\nu_0+1)\eta^\Lambda$.

The last assertion follows from the description of semisimple Riemannian F-manifolds in \cite{ABLR}. In
canonical coordinates the metric is diagonal, and its Lam\'e coefficients determine the natural connection by
the formulas in (i).
\end{proof}

{
\begin{corollary}\label{cor:lame}
Let $(M,\nabla,\circ,e)$ be a semisimple flat F-manifold with linear Euler vector field $E$ and canonical
coordinates $r^1,\dots,r^n$. Let $H_1,\dots,H_n$ be Lam\'e coefficients with Killing unit, namely functions
satisfying $\partial_j\log H_i=\Gamma^i_{ij}$ for $j\neq i$ and $e(H_i)=0$. These functions are unique up to
constant factors. Then each $H_i$ is homogeneous of constant degree $d_i:=E(\log H_i)$. If $r^i$ is the value
at a simple critical point of a superpotential of weight $\nu_0$, then $d_i=\nu_0-\frac12$. In particular, for
a superpotential of weight $\nu_0$ satisfying the hypotheses of Theorem~\ref{thm:main}, we have
\[
d_1=d_2=\dots=d_n=\nu_0-\tfrac12 .
\]
\end{corollary}

\begin{proof}
Flatness of $\nabla$ together with $\nabla\nabla E=0$ implies $\mathrm{Lie}_E\nabla=0$. In canonical
coordinates we have $E=\sum_kr^k\partial_k$, so this identity reads $E(\Gamma^i_{jk})=-\Gamma^i_{jk}$. For
$j\neq i$, the commutator relation $[\partial_j,E]=\partial_j$ gives
\[
\partial_j\big(E\log H_i\big)=E\big(\partial_j\log H_i\big)+\partial_j\log H_i=E(\Gamma^i_{ij})+\Gamma^i_{ij}=0 .
\]
Moreover, $e\big(E\log H_i\big)=E\big(e\log H_i\big)+[e,E]\log H_i=e(\log H_i)=0$. Since $e=\sum_j\partial_j$,
it follows that $\partial_i(E\log H_i)=0$ as well. Thus $d_i$ is constant. Now suppose that
$r^i=\Lambda(P_*^i)$, where $P_*^i$ is a simple critical point. By Remark~\ref{rem:simplecrit} we have
$\Lambda_P|_i=0$ and $\Lambda_j|_i=\delta_{ij}$. For the index $i$, these are the only facts concerning the
critical values used in part~(ii) of Proposition~\ref{prop:christoffel} and in the proof of
Proposition~\ref{prop:riem}. We therefore have $\Gamma^i_{ij}=-\frac12\partial_j\log\Lambda_{PP}|_i$ for
$j\neq i$, together with $e(\Lambda_{PP}|_i)=0$ and $E(\Lambda_{PP}|_i)=(1-2\nu_0)\Lambda_{PP}|_i$. The first
two identities imply that $(\Lambda_{PP}|_i)^{-1/2}$ agrees with $H_i$ up to a constant factor. The third then
gives $d_i=\nu_0-\frac12$.
\end{proof}

\begin{proposition}\label{prop:lamecrit}
Suppose that on a chart the canonical coordinates $r^1,\dots,r^n$ are pairwise distinct. Then there exists a
period of weight $\nu_0$ whose continued superpotential has all the $r^i$ as simple critical values if and
only if $d_1=\dots=d_n=\nu_0-\frac12$. More generally, if $\nu_0-d_i\notin\mathbb Z$ and
$\operatorname{Re}(\nu_0-d_i)<1$, then for a generic period we have
$\Lambda-r^i\sim c\,(P-P_i)^{1/(\nu_0-d_i)}$ near $r^i$.
\end{proposition}

\begin{proof}
In canonical coordinates, the equation in $\lambda$ of system~\eqref{eq:xisys} has a simple pole at $r^i$ with
residue $R_i=N\Pi_i$, where $\Pi_i$ is the projection onto the direction indexed by $i$. This residue has rank
one and eigenvalue $\nu_0-1-d_i$, since $Q^i_i=1+E(\log H_i)=1+d_i$. If the eigenvalue is not an integer, the
normal form at a Fuchsian singularity in \cite{IY} gives
$\xi(e)=h(\lambda)+(\lambda-r^i)^{\nu_0-1-d_i}\varphi(\lambda)$, where $h,\varphi$ are holomorphic and
$\varphi(r^i)\ne0$ whenever the initial differential lies outside a hyperplane. Integrating
$\partial_\lambda P=-\xi(e)$, we find $P-P_i=(\lambda-r^i)^{\nu_0-d_i}f(\lambda)+O(\lambda-r^i)$ with
$f(r^i)\ne0$. When $\nu_0-d_i=\frac12$, the function $t=(\lambda-r^i)^{1/2}$ is then a local coordinate as a
function of $P$, and $\Lambda=r^i+t^2$. The necessity of the condition on the degrees follows from
Corollary~\ref{cor:lame}.
\end{proof}

\begin{remark}\label{rem:lame}
For $B_2$, as in Section~\ref{sec:Bn}, we have $\nu_0-d_1=\frac12$ and $\nu_0-d_2=\beta$. Indeed,
$\Lambda-a_2\sim a_1(1+\beta P_\beta)^{1/\beta}$. Thus the canonical coordinate $a_2$ is a simple critical
value only for $\beta=\frac12$, in which case $\Lambda=p^4+a_1p^2+a_2$ is Dubrovin's superpotential. It is a
critical value of multiplicity $k-1$ for $\beta=\frac1k$. When $\beta=0$ the exponent $\nu_0-d_2=0$ is an
integer and the behavior is logarithmic. More precisely, $\Lambda-a_2\sim a_1e^{P_0}$, so $a_2$ is an
asymptotic value as $\operatorname{Re}P_0\to-\infty$. For the period $q=\log z_1-\log z_2$ we obtain
$\Lambda=a_2-a_1^2e^q/(1+e^q)^2$, which never takes the value $a_2$.
\end{remark}
}

}

{ We now start from a commutative associative product with unit, a vector field $E$, and a function $\Lambda$
satisfying (S2) with $\partial_P\Lambda\ne0$. Suppose that \eqref{eq:chr} admits coefficients $\Gamma^m_{hk}$
independent of $P$. We shall see that if $n$ local inverse branches of $\Lambda$ have independent
differentials in the base variables, then these coefficients determine a unique flat torsion free connection
$\nabla$. If (S1) also holds and $\nu_0\Id-\nabla E$ is invertible, the product and $\nabla$ define a flat F-manifold with linear Euler field $E$. Proposition~\ref{prop:dualflat} then gives the dual structure wherever
$E\circ$ is invertible.

\begin{theorem}\label{thm:flat}
{Let $M$ be a manifold of dimension $n$ endowed with a commutative associative product $\circ$, a unit $e$,
and a vector field $E$. Set $L=E\circ$. Let $\Lambda(P,x)$ satisfy (S2) with $\partial_P\Lambda\neq0$. Suppose
that for $\lambda$ in an open set the equation $\Lambda(P,x)=\lambda$ has $n$ solutions
$P_1(x,\lambda),\dots,P_n(x,\lambda)$ with linearly independent differentials $d_xP_1,\dots,d_xP_n$. Suppose
further that in local coordinates $x^1,\dots,x^n$ there exist functions $\Gamma^m_{hk}(x)$, independent of
$P$, such that
\begin{equation}\label{eq:defconn}
\Lambda_{hk}-\partial_P\Big[\frac{\Lambda_h\Lambda_k-c^m_{hk}\Lambda_m}{\Lambda_P}\Big]=\Gamma^m_{hk}\,\Lambda_m
\end{equation}
identically in $P$. In other words, the left side, considered as a function of $P$, belongs to the span of
$\Lambda_1,\dots,\Lambda_n$. Then the coefficients $\Gamma^m_{hk}$ are unique and define the Christoffel
symbols of a torsion free connection $\nabla$ on $M$. This connection is independent of the choice of
coordinates, and $\Lambda$ satisfies (S3) in the form \eqref{eq:S3cov} with respect to it. Set
$N:=\nu_0\Id-\nabla E$. Then the following assertions hold.}
\begin{enumerate}
\item[(i)] $\nabla$ is flat;
\item[(ii)] if $N$ is invertible, the tensor $\nabla_lc^m_{hk}$ is symmetric in its lower indices;
\item[(iii)] if $\Lambda$ also satisfies (S1), then $\nabla e=0$;
\item[(iv)] $\nabla\nabla E=0$, that is, $E$ is linear in the flat coordinates of $\nabla$;
\item[(v)] if $N$ is invertible, then $\mathrm{Lie}_E\circ=\circ$, and if (S1) also holds, then $[e,E]=e$.
\end{enumerate}
In particular, if $\Lambda$ also satisfies (S1) and $N$ is invertible, then $(M,\nabla,\circ,e)$ is a flat F-manifold with linear Euler vector field $E$.
\end{theorem}

\begin{proof}
{ We shall derive from the conditions on $\Lambda$ a linear system of first order in $(x,\lambda)$ for the
differentials of the $n$ branches $P_j$. Since these differentials are independent, the system is integrable.
Its compatibility conditions are rational in $\lambda$, and their expansion at $\lambda=\infty$ will give the
assertions of the theorem. We use row vectors for covectors and let endomorphisms act on the right, so that
$(\xi A)_k=\xi_mA^m_k$. We write $(C_h)^l_k=c^l_{hk}$, $(\nabla E)^m_k=\nabla_kE^m$, $L=E\circ$, and
$L_\lambda=L-\lambda$.

\emph{Step 0.} We first construct the connection. For each fixed $x$, the functions
$\Lambda_1,\dots,\Lambda_n$ of $P$ are linearly independent. To see this, differentiate
$\Lambda(P_j(x,\lambda),x)=\lambda$ to obtain $\Lambda_m=-\Lambda_P\,\partial_mP_j$ at $P=P_j$. If
$a^m\Lambda_m=0$ identically in $P$, then $a^m\partial_mP_j=0$ for every branch. Since the covectors $d_xP_j$
form a basis, it follows that $a=0$. Thus \eqref{eq:defconn} determines the coefficients $\Gamma^m_{hk}$
uniquely. The symmetry of its left side in $h,k$ implies the same symmetry for the coefficients.

We next determine how the coefficients transform. Let $\tilde x$ be another system of coordinates and set
$J^h_a=\partial x^h/\partial\tilde x^a$. Use tildes to denote derivatives with respect to $\tilde x$ at fixed
$P$. The chain rule gives $\tilde\Lambda_a=J^h_a\Lambda_h$ and
\[
\tilde\Lambda_{ab}=J^h_aJ^k_b\,\Lambda_{hk}+\frac{\partial^2x^m}{\partial\tilde x^a\partial\tilde x^b}\,\Lambda_m ,
\]
whereas the numerator in brackets in \eqref{eq:defconn} transforms as a tensor. Indeed,
$\tilde\Lambda_a\tilde\Lambda_b-\tilde c^m_{ab}\tilde\Lambda_m=J^h_aJ^k_b\,(\Lambda_h\Lambda_k-c^m_{hk}\Lambda_m)$,
and both $\partial_P$ and division by $\Lambda_P$ commute with multiplication by functions of $x$. Therefore,
if \eqref{eq:defconn} holds in the coordinates $x$, its left side in the coordinates $\tilde x$ is
\[
J^h_aJ^k_b\,\Gamma^m_{hk}\Lambda_m+\frac{\partial^2x^m}{\partial\tilde x^a\partial\tilde x^b}\Lambda_m
=\Big(J^h_aJ^k_b\,\Gamma^m_{hk}+\frac{\partial^2x^m}{\partial\tilde x^a\partial\tilde x^b}\Big)\big(J^{-1}\big)^c_m\,\tilde\Lambda_c ,
\]
so that \eqref{eq:defconn} holds also in the coordinates $\tilde x$, with
$\tilde\Gamma^c_{ab}=(J^{-1})^c_m\big(J^h_aJ^k_b\Gamma^m_{hk}+\partial^2x^m/\partial\tilde x^a\partial\tilde x^b\big)$.
This is precisely the transformation law for the Christoffel symbols of a connection. Thus the coefficients
$\Gamma^m_{hk}$ define a torsion free connection $\nabla$ on $M$, and \eqref{eq:defconn} is the covariant form
\eqref{eq:S3cov} of (S3) for $\nabla$.

\emph{Step 1.} We now derive the equations for the differentials of the branches. Let $P$ be one of the
branches $P_j$, and set $\xi=d_xP$ and $\xi'=\partial_\lambda\xi$.
\begin{itemize}
\item The identity in the proof of Lemma~\ref{lem:inv}, namely
$\Lambda_P\big(\partial_h\partial_kP+c^l_{hk}\,\partial_\lambda\partial_lP\big)=-\Lambda_{hk}+\partial_PB_{hk}$,
follows from the chain rule and thus holds in arbitrary coordinates. Adding
$-\Lambda_P\Gamma^m_{hk}\partial_mP=\Gamma^m_{hk}\Lambda_m$ to both sides gives
$\Lambda_P\big(\nabla_h\xi_k+c^l_{hk}\xi'_l\big)=-\big(\Lambda_{hk}-\Gamma^m_{hk}\Lambda_m\big)+\partial_PB_{hk}$.
Therefore \eqref{eq:S3cov} is equivalent to
\begin{equation}\label{eq:covxi}
\nabla_h\xi=-\xi'\,C_h .
\end{equation}
\item By the calculation in part~(ii) of Theorem~\ref{thm:char}, (S2) is equivalent to
$(E+\lambda\partial_\lambda)P=\nu_0P+\kappa$. Differentiating and using \eqref{eq:covxi} in the identity
$\nabla_k(E^m\xi_m)=(\xi\nabla E)_k+E^m\nabla_k\xi_m=(\xi\nabla E)_k-(\xi'L)_k$, we obtain
$-\xi'L_\lambda=\xi\,N$. Thus, wherever the inverse exists, we have
\begin{equation}\label{eq:covxi2}
\xi'=-\xi\,NL_\lambda^{-1},\qquad\text{and thus}\qquad \nabla_h\xi=\xi\,NL_\lambda^{-1}C_h .
\end{equation}
\item The calculation in part~(i) of Theorem~\ref{thm:char} gives the equivalence between (S1) and
$\partial_\lambda P=-e^k\xi_k$.
\end{itemize}

\emph{Step 2.} We eliminate the parameter by the following observation. Let $A$ and $B$ be tensors on $M$
independent of $\lambda$. If $\xi A=\xi\,NL_\lambda^{-1}B$ holds for the differentials $\xi$ of all $n$
branches, then $A=0$ and $NB=0$. Indeed, fix $x$. Since the $n$ covectors $\xi$ form a basis,
$A=NL_\lambda^{-1}B$ holds for $\lambda$ in an open set. The right side is rational in $\lambda$ and agrees
with a constant on an infinite set. Thus the identity holds for every $\lambda$ outside the spectrum of $L$.
In particular, it holds for $|\lambda|$ sufficiently large, where
\[
L_\lambda^{-1}=-\frac1\lambda\Big(\Id-\frac L\lambda\Big)^{-1}=-\sum_{k\ge0}\frac{L^k}{\lambda^{k+1}} .
\]
It follows that $A=-\lambda^{-1}NB-\lambda^{-2}NLB-\dots$. Since the left side is constant in $\lambda$, its
value at infinity gives $A=0$, and comparison of the coefficients of $\lambda^{-1}$ gives $NB=0$. This
argument is pointwise in $x$ and requires no analyticity in $x$.

\emph{Step 3.} To prove flatness and the symmetry of $\nabla c$ in parts (i) and (ii), apply $\nabla_j$ to
equation~\eqref{eq:covxi}, namely $\nabla_h\xi_k=-c^l_{hk}\xi'_l$. Since $\partial_\lambda$ commutes with
$\nabla_j$, equation~\eqref{eq:covxi} also gives $\nabla_j\xi'_l=-c^m_{jl}\xi''_m$. We obtain
\[
\nabla_j\nabla_h\xi_k=-(\nabla_jc^l_{hk})\,\xi'_l+c^l_{hk}c^m_{jl}\,\xi''_m .
\]
Associativity implies that the last term is symmetric in $j,h$. Antisymmetrizing in $j,h$ and applying the
Ricci identity $(\nabla_j\nabla_h-\nabla_h\nabla_j)\xi_k=-R^m_{kjh}\xi_m$ for the torsion free connection
$\nabla$, we obtain, for fixed $j,h$,
\[
\xi\,R_{jh}=\xi'\,T_{jh}=-\xi\,NL_\lambda^{-1}T_{jh},\qquad
(R_{jh})^m_k=R^m_{kjh},\quad (T_{jh})^l_k=\nabla_jc^l_{hk}-\nabla_hc^l_{jk},
\]
where the last equality follows from \eqref{eq:covxi2}. Applying Step~2 with $A=R_{jh}$ and $B=-T_{jh}$ gives
$R_{jh}=0$, which proves (i). It also gives $NT_{jh}=0$. Therefore $T_{jh}=0$ when $N$ is invertible, which
proves (ii).

\emph{Step 4.} For part (iii), differentiate the third identity in Step~1. Using \eqref{eq:covxi} and the unit
identity $c^l_{hk}e^k=\delta^l_h$ gives
\[
\xi'_h=-\nabla_h(e^k\xi_k)=-(\nabla_he^k)\,\xi_k-e^k\nabla_h\xi_k=-(\nabla_he^k)\,\xi_k+\xi'_h ,
\]
so $(\nabla_he^k)\xi_k=0$ for each covector $\xi$ in a basis. This proves $\nabla e=0$.

\emph{Step 5.} To prove parts (iv) and (v), we compute the two sides of
$\nabla_h\xi'=\partial_\lambda(\nabla_h\xi)$. From $\partial_\lambda L_\lambda^{-1}=L_\lambda^{-2}$ and
$\xi'=-\xi NL_\lambda^{-1}$ we obtain
\[
\partial_\lambda\big(\xi\,NL_\lambda^{-1}C_h\big)=-\xi\,NL_\lambda^{-1}NL_\lambda^{-1}C_h+\xi\,NL_\lambda^{-2}C_h .
\]
On the other hand, from $\nabla_h(L_\lambda^{-1})=-L_\lambda^{-1}(\nabla_hL)L_\lambda^{-1}$ and
$\nabla_h\xi=\xi NL_\lambda^{-1}C_h$ we obtain
\[
\nabla_h\big(-\xi\,NL_\lambda^{-1}\big)=-\xi\,NL_\lambda^{-1}C_hNL_\lambda^{-1}-\xi\,(\nabla_hN)L_\lambda^{-1}+\xi\,NL_\lambda^{-1}(\nabla_hL)L_\lambda^{-1} .
\]
Equating the two expressions and multiplying on the right by $L_\lambda$, we use $C_hL_\lambda=L_\lambda C_h$
to obtain $\xi\,\nabla_hN=\xi\,NL_\lambda^{-1}Y_h$, where
\[
Y_h:=\nabla_hL-C_h+[N,C_h]=\nabla_hL-C_h-[\nabla E,C_h] .
\]
By Step~2 we have $\nabla_hN=0$, or equivalently $\nabla\nabla E=0$, as claimed in (iv). We also have
$NY_h=0$, so $Y_h=0$ when $N$ is invertible. Writing $L^m_k=c^m_{kl}E^l$, we can express $Y_h=0$, namely
$\nabla_hL=C_h+[\nabla E,C_h]$, in components as
\[
\nabla_hL^m_k=c^m_{hk}+c^l_{hk}\,\nabla_lE^m-c^m_{hl}\,\nabla_kE^l .
\]
On the left side we have $\nabla_hL^m_k=E^l\,\nabla_hc^m_{kl}+c^m_{kl}\,\nabla_hE^l$, and part~(ii) gives
$\nabla_hc^m_{kl}=\nabla_lc^m_{hk}$. Substituting these expressions, we obtain
\[
E^l\nabla_lc^m_{hk}-c^l_{hk}\nabla_lE^m+c^m_{lk}\nabla_hE^l+c^m_{hl}\nabla_kE^l=c^m_{hk} .
\]
Since $\nabla$ is torsion free, the left side is the expression for $(\mathrm{Lie}_Ec)^m_{hk}$ in terms of
covariant derivatives. Thus $\mathrm{Lie}_E\circ=\circ$.

Finally, suppose that (S1) also holds. By part~(iii) we have $\nabla e=0$, and the unit property gives
$L^m_ke^k=c^m_{kl}e^kE^l=E^m$. Contract the component identity for the covariant derivative of multiplication
by the Euler field with $e^k$. On the left we obtain $(\nabla_hL^m_k)\,e^k=\nabla_h(L^m_ke^k)=\nabla_hE^m$. On
the right, $c^m_{hk}e^k=\delta^m_h$ and $c^l_{hk}e^k=\delta^l_h$ give
$\delta^m_h+\nabla_hE^m-c^m_{hl}\,e^k\nabla_kE^l$. It follows that
\[
c^m_{hl}\,(\nabla_eE)^l=\delta^m_h\qquad\text{for all }h,
\]
that is, $\partial_h\circ\nabla_eE=\partial_h$. Thus $\nabla_eE=e$ and $[e,E]=\nabla_eE-\nabla_Ee=e$. }
\end{proof}

\begin{remark}\label{rem:flat}
Suppose that $\Lambda$ satisfies (S1) and (S2) and that the left side of \eqref{eq:defconn} lies in the span
of the $\Lambda_m$. Assume the algebraic hypotheses on the product and suppose there are $n$ independent
branches. If $N$ is invertible, the function then determines the connection, and all the axioms of a flat F-manifold with linear Euler vector field follow. In particular, the connection is flat, $\nabla c$ is
symmetric, the unit is parallel, and $E$ is linear in flat coordinates and satisfies the homogeneity
identities. By Proposition~\ref{prop:dualflat} this structure determines a bi-flat F-manifold wherever $E\circ$
is invertible. We use this construction in Sections~\ref{sec:Bn} and~\ref{sec:Gmmn}. There the superpotential
gives the product, and \eqref{eq:chr} defines the connection, so (S3) holds by construction. When $N$ is
invertible, the remaining properties follow from the theorem. At the remaining weights they follow by
continuation in $\beta$. Indeed, the product is independent of $\beta$, so $\mathrm{Lie}_E\circ=\circ$
continues to hold. The Christoffel symbols are polynomial in $\beta$, so the symmetry of $\nabla c$ also
persists.
\end{remark}
}

{ The hypotheses of Theorem~\ref{thm:main} can fail even when the product is semisimple.
\begin{remark}\label{rem:ss}
Both the continuation and the differentials of the critical values depend on the particular superpotential.
For instance, $\Lambda=u^n-P$ has no critical points even when the product is semisimple.
\end{remark}

Let us relate the coordinate $P$ on the fiber to Dubrovin's primary differential.
\begin{remark}\label{rem:primary}
In a model of Landau and Ginzburg one starts with a curve endowed with a function $\lambda$ and a primary
differential $\omega$. Dubrovin's periods are integrals of $\omega$ and of $\lambda^s\omega$, as in
\cite{Dub}. Wherever $P=\int\omega$ is a coordinate, we have $\omega=dP$ along the curve. The values of $P$ at
the preimages of $\lambda$ then give local flat functions for the connection of Gauss and Manin. In
Definition~\ref{def:sp} we write the inverse in the fiber coordinate $q$. If its expression simplifies in
another coordinate $p$, the primary form in this coordinate is $\omega=(\partial q/\partial p)\,dp$ at fixed
base point and may be multivalued. For instance, in the $B_n$ family, $\Lambda$ is polynomial in $z=p^2$ and
we have $\omega=z^{\beta-1}dz\propto p^{2\beta-1}dp$. At $\beta=1/2$ this form is proportional to $dp$.
\end{remark}

Under the hypotheses of Theorem~\ref{thm:main} the critical values give canonical coordinates. The simple
zeros, in turn, give flat functions on the zero slice.
\begin{corollary}\label{cor:flat}
Let $\Lambda$ be a superpotential of weight $\nu_0$ and constant $\kappa$. On the domain where $E\circ$ is
invertible, suppose that $P^{(j)}(u)$ is a simple zero on a sheet of the continuation of $\Lambda$. Thus
$\Lambda(P^{(j)}(u),u)=0$ and $\partial_P\Lambda\ne0$ at the zero. Then $P^{(j)}$ is flat for
$\Gamma(0,\nu_0)=\nabla^*+\nu_0c^*$ and satisfies $E(P^{(j)})=\nu_0P^{(j)}+\kappa$. If there are $n$ such
functions with independent differentials, they form a system of flat coordinates.
\end{corollary}

\begin{proof}

At a simple zero the continued superpotential admits a holomorphic inverse near $\lambda=0$. Thus
$P^{(j)}(u)=P(u,0)$ for the corresponding branch $P$ of the inverse. Parts (ii) and (iv) of
Theorem~\ref{thm:char} give the homogeneity and flatness equations, respectively. When the $n$ differentials
are independent, these functions give local coordinates.
\end{proof}

In Sections~\ref{sec:Bn} and~\ref{sec:Gmmn} we obtain the flat coordinates of $\nabla$ as residues of powers
of $\Lambda$ at its poles. These residues arise from the expansion of $P$ as $\lambda\to\infty$, as in
\cite{Dub}. }

\subsection{Reduction to Dubrovin's construction}\label{sec:Dub}

{
Let $M$ be a semisimple Frobenius manifold with invariant metric $\eta$ and charge $d$. Its intersection
form is $g^{\alpha\beta}=E^\varepsilon c^{\alpha\beta}_\varepsilon$. Set $\mathcal V=(1-\frac d2)\Id-\nabla E$, so that $\mathcal V-\frac12=\frac{1-d}2-Q$. Denote by $\widetilde\nabla$ the Levi-Civita connection of the intersection form and
by $\nabla^*$ the normalized dual connection satisfying $\nabla^*E=0$. Throughout this subsection we work on the
regular locus of the corresponding pencils.

To recover Dubrovin's construction, we first determine the weight for which the connection on the
slice $\lambda=0$ coincides with the Levi-Civita connection of the intersection form. For this weight we can
identify the equations on the product with the equations for Dubrovin's periods.

\begin{proposition}\label{prop:Dub}
\begin{enumerate}
\item[(i)] In flat coordinates of $\eta$, the coefficients of $\widetilde\nabla$ are given by $\Gamma^\alpha_{\beta\gamma}=(\nu_0\Id-Q)^\alpha_\varepsilon c^{*\varepsilon}_{\beta\gamma}$, where $\nu_0=(1-d)/2$. In
particular, $\widetilde\nabla=\nabla^*+\nu_0*$ and $\widetilde\nabla E=\nu_0\Id$.
\item[(ii)] For this weight, the $u$ and mixed components of $\hat\nabla$ give the period equations of the flat
pencil $g-\lambda\eta^{-1}$ of \cite{Dalm}. The translated periods are those flat functions $P$ of the extended connection that satisfy $(\partial_\lambda+e)P=0$. If $d\ne1$,
one can impose the homogeneous normalization of \cite{Dalm}.
\item[(iii)] At a regular point where $\partial_\lambda P\ne0$, the inverse of a translated period is the local inverse
of Theorem~\ref{thm:char}, with primary form $dp$ along the fiber. Whenever the unit flow to the
regular part of the zero slice is defined, this inverse is the superpotential of
Definition~\ref{def:sp}. In flat coordinates with $e=\partial_{t^1}$, conditions (S1) and (S2) read $\partial_{t^1}\lambda=1$ and $E\lambda+(\nu_0p+\kappa)\lambda_p=\lambda$,
respectively, while (S3) is equivalent to $\partial_\alpha\partial_\beta p=-c^\gamma_{\alpha\beta}
\partial_\lambda\partial_\gamma p$. If $d\ne1$, one can shift $p$ by a constant so that $\kappa=0$.
If $d=1$, one obtains the period by integrating its parallel differential and must keep the constant
$\kappa$ in (S2).
\end{enumerate}
\end{proposition}

\begin{proof}

To prove (i), recall from \cite{Dalm} that $\Gamma^\alpha_{\beta\gamma}=-g_{\beta\rho}\Gamma^{\rho\alpha}_\gamma$, where $\Gamma^{\rho\alpha}_\gamma=c^{\rho\varepsilon}_\gamma(\frac12-\mathcal V)^\alpha_\varepsilon$. Since $g^{\rho\sigma}=\eta^{\rho\kappa}L^\sigma_\kappa$ and $L$ is self adjoint with
respect to $\eta$, the inverse metric satisfies $g_{\beta\rho}=(L^{-1})^\kappa_\beta\eta_{\kappa\rho}$. It follows that
\[
g_{\beta\rho}\,c^{\rho\varepsilon}_\gamma=(L^{-1})^\kappa_\beta\,\eta^{\varepsilon\nu}c_{\nu\gamma\kappa}
=c^\varepsilon_{\gamma\kappa}(L^{-1})^\kappa_\beta=(C_\gamma L^{-1})^\varepsilon_\beta=c^{*\varepsilon}_{\beta\gamma},
\]

where we have used the commutativity of the multiplication operators. Thus $\Gamma^\alpha_{\beta\gamma}=(\mathcal V-\frac12)^\alpha_\varepsilon
c^{*\varepsilon}_{\beta\gamma}$. Finally, from $\widetilde\nabla_\beta E^\alpha=Q^\alpha_\beta+\Gamma^\alpha_{\beta\gamma}E^\gamma$ and
the unit identity $c^{*\varepsilon}_{\beta\gamma}E^\gamma=\delta^\varepsilon_\beta$ we obtain $\widetilde\nabla E=\nu_0\Id$.

For (ii), identify Dubrovin's $\mathcal U$ with $L$. In the notation of \cite{Dalm}, the period equations read
$\partial_\alpha\xi\,(L-\lambda)=\xi(\mathcal V-\frac12)C_\alpha$ and $\partial_\lambda\xi\,(L-\lambda)=\xi(\frac12-\mathcal V)$. With $N=\frac{1-d}2-Q$, they take the form $\partial_\alpha\xi=\xi\,NL_\lambda^{-1}C_\alpha$ and $\partial_\lambda\xi=-\xi\,NL_\lambda^{-1}$, which is precisely \eqref{eq:xisys}.

Conversely, let $\xi$ be a row solution. To see that the one form $\xi_h\,du^h-\xi(e)\,d\lambda$ is closed, observe that the $u$
equations are symmetric in the lower indices and that $\partial_h(\xi(e))=-\partial_\lambda\xi_h$ since $C_he=\partial_h$. A local primitive $P$ therefore
satisfies $\partial_\lambda P=-\partial_eP$. The row equations give the components of its Hessian in the base directions as well
as the mixed components. They also give $\partial_\lambda^2 P=-\partial_\lambda(\xi(e))=\xi NL_\lambda^{-1}e$, which shows that $P$ is flat for $\hat\nabla$. When $\nu_0\ne0$ and $\kappa=0$,
the homogeneity equation gives $P=\nu_0^{-1}\langle d_uP,E-\lambda e\rangle$. For $\nu_0=0$ the contraction equals $\kappa$, and we obtain $P$ locally by
integrating its differential.

Applying Corollary~\ref{cor:spext} and Theorem~\ref{thm:char} to the translated period proves (iii).
\end{proof}

We now explain how the differential characterization of the superpotential relates to the open WDVV
equations.
\begin{remark}\label{rem:PS}
The open WDVV equations describe extensions of a Frobenius manifold by one coordinate. Locally such
an extension has a vector potential on $M\times\mathbb C$, with a flat coordinate $x$ along the fibers, and the
projection to $M$ preserves multiplication. Theorem~3 of \cite{Alc} gives a criterion for such an
extension which applies to any almost flat F-manifold. See also Theorem~2.3 of \cite{PS}. More
precisely, if a function $h$ satisfies $\partial_xh\ne0$ and
\[
\partial_\alpha\partial_\beta h
=\partial_x\left[\frac{\partial_\alpha h\,\partial_\beta h-
c^\mu_{\alpha\beta}\partial_\mu h}{\partial_xh}\right],
\]
then it determines the product of a local extension in these coordinates. Indeed, the last component
$\Omega$ of the vector potential satisfies $\partial_x\Omega=h$, and its Hessian in the base directions is given by the
quotient in brackets. These derivatives determine the product on the extension.

For a Hurwitz Frobenius manifold, Almeida \cite{Alm} proves that its potential satisfies this
criterion with $h=\lambda$ in a local fiber coordinate given by the primary differential. Proserpio and
Strachan \cite{PS} prove the corresponding result for their almost dual Hurwitz structures with $h=\log\lambda$.
Taking $x=P$ and $h=\Lambda$, we see that the criterion is precisely (S3). Moreover, equation (3.9) of \cite{PS}
for homogeneity agrees with (S2) when $\nu_0=(1-d)/2$ and $\kappa=d_0$. Thus any local superpotential of a bi-flat F-manifold defines such an extension. Its unit is the lift of $e$ with zero $P$ component. Indeed, (S1)
gives $\Omega_{eP}=1$, while the formula for the Hessian in the base directions gives $\Omega_{eh}=0$ for each flat coordinate
direction $\partial_h$ on the base.
\end{remark}
}

\subsection{Summary}\label{subsec:summary}
Table~\ref{tab:reconstruction} collects the local reconstruction results proved in this section. Let
$P_i$ be a simple critical point of the continued superpotential in the primary coordinate $P$, and set
$r^i=\Lambda(P_i)$ and $h_i=1/\Lambda_{PP}(P_i)$. The differentials $dr^i$ are characters of the tangent algebra at the base point. Distinct
characters are linearly independent, although several critical points can give the same character.
On the semisimple locus each critical value $r^i$ obtained in this way is a canonical coordinate.

\begin{table}[htbp]
\caption{What the superpotential determines locally}
\label{tab:reconstruction}
\small
\renewcommand{\arraystretch}{1.12}
\begin{tabular}{@{}>{\raggedright\arraybackslash}p{0.325\linewidth}@{\hspace{0.03\linewidth}}>{\raggedright\arraybackslash}p{0.645\linewidth}@{}}
\toprule
\textbf{Data or hypothesis} & \textbf{Conclusion and uniqueness}\\
\midrule
$k$ simple critical points with independent differentials $dr^i$, where $1\le k<n$
&
The map $X\mapsto(dr^i(X))_{i=1}^k$ identifies the quotient by
$\bigcap_i\ker dr^i$ with $\mathbb C^k$ equipped with componentwise multiplication.
The sum $\sum_{i=1}^k h_i(dr^i)^2$ has rank $k$. In general these data do not determine
the product on the full tangent space
(Remark~\ref{rem:simplecrit}).\\
\addlinespace
$n$ simple critical points with distinct $dr^i$
&
The critical values $r^i$ are canonical coordinates and uniquely determine the
product, the unit, and the Euler field. Writing $\partial_i=\partial/\partial r^i$,
we have $\partial_i\circ\partial_j=\delta_{ij}\partial_i$,
$e=\sum_i\partial_i$ and $E=\sum_i r^i\partial_i$
(Theorem~\ref{thm:main}).\\
\addlinespace
The same $n$ critical values and the functions $h_i$ at the chosen critical points
&
The residue metric is $\eta^\Lambda=\sum_i h_i(dr^i)^2$, with Lam\'e coefficients satisfying $H_i^2=h_i$.
For $j\ne i$ we have $\Gamma^i_{ij}=\tfrac12\partial_j\log h_i$.
These symbols determine $\nabla$ uniquely. Multiplying each $h_i$ by
a nonzero constant can change the metric while preserving the connection
(Lemma~\ref{lem:canform} and Proposition~\ref{prop:riem}).\\
\addlinespace
Lam\'e degrees $d_i=E(\log H_i)$ on a semisimple chart, with $e(H_i)=0$
&
A superpotential of weight $\nu_0$ can realize $r^i$ as a simple critical
value only if $d_i=\nu_0-\tfrac12$. If all $d_i=d_0$ and the Euler
eigenvalues are pairwise distinct, a local period of weight
$d_0+\tfrac12$ realizing every $r^i$ exists. A generic choice of period has this property
(Corollary~\ref{cor:lame}, Proposition~\ref{prop:lamecrit} and Remark~\ref{rem:lame}).\\
\addlinespace
A commutative associative product with unit, a vector field $E$, and $n$ regular inverse branches
over an open set of levels, with independent differentials $d_xP_j$ in the base variables
&
Assume (S2). If \eqref{eq:defconn} admits a solution independent of $P$, it defines a
unique flat torsion free connection $\nabla$ for the given product.
If (S1) also holds and $\nu_0\Id-\nabla E$ is invertible, the data define
a flat F-manifold with linear Euler field (Theorem~\ref{thm:flat}).\\
\bottomrule
\end{tabular}
\end{table}

\section{The case of two orbits}\label{sec:two}

{
We first consider $B_n=G(2,1,n)$. In this case a monic polynomial $\lambda(z)$ of degree $n$ and
a primary form $z^{\beta-1}dz$ determine the entire bi-flat structure. More precisely, the
polynomial determines both products, while the primary form determines the weights of the dual
connection of Dunkl and Kohno and the residues which give flat coordinates for the natural
connection (Theorem~\ref{thm:Bn}). The substitution $p=z^{1/m}$ then gives the family for
$G(m,1,n)$ described in Section~\ref{sec:Gm1n}. In the case of $F_4$, we give the coordinates,
the vector potential, and the connection directly in Section~\ref{sec:F4}. In both cases there
are two orbits, so the family in Conjecture~\ref{conj:AL} depends on one parameter.
}

\subsection{\texorpdfstring{The $B_n$ families}{The Bn families}}\label{sec:Bn}
\subsubsection{The data}\label{sec:construction}

{
Let $p=(p_1,\dots,p_n)$ be Euclidean coordinates on $\mathbb C^n$. The action of $B_n$ consists
of permutations and sign changes of the coordinates $p_j$. Thus the coefficients of the following
polynomial give coordinates on the quotient.
{
\begin{equation}\label{eq:Bn}
\lambda(z)=\prod_{j=1}^n(z-z_j)=z^n+a_1z^{n-1}+\dots+a_n ,\qquad z_j=p_j^2 ,
\end{equation}
}
We may therefore identify the quotient with the space $M$ of monic polynomials of degree $n$. The
polynomial $\lambda$ is homogeneous of degree $1$ if we assign the degrees $\deg z=\frac1n$ and
$\deg a_i=\frac in$. Recall that in Dubrovin's construction for the Frobenius manifold $B_n$ the
superpotential is $\lambda(p^2)$ and the primary form is given by $\omega=dp$. For the choice
$\beta=1/2$ below, the extension of Remark~\ref{rem:PS} agrees, up to an affine change of the
fiber coordinate, with the solution of the open WDVV equations of type $B_n$ in formula (5.8) of
\cite{PS}. This solution is obtained by restriction from $A_{2n-1}$.

We now allow the primary form to vary and set $\omega=dP_\beta=z^{\beta-1}dz$, where
$\beta\in\mathbb C$. After choosing a local branch of the logarithm, we take the following
primitive, which is also defined at $\beta=0$.
\[
P_\beta(z)=
\begin{cases}
(z^\beta-1)/\beta,&\beta\ne0,\\
\log z,&\beta=0.
\end{cases}
\]
The primitive depends holomorphically on $\beta$ and, for $\beta\ne0$, differs from
$z^\beta/\beta$ by a constant. The primary form has weight $\nu_0=\beta/n$ and, when $\beta=1/2$,
is proportional to $dp$. To express the inverse branches in the coordinate $P_\beta$, we apply
$P_\beta$ to each of the roots of the polynomial defining the level. This gives
\begin{equation}\label{eq:ansatz}
P_j(a,\lambda)=P_\beta(z_j(a,\lambda)),\qquad
\lambda(z_j(a,\lambda),a)=\lambda,\qquad j=1,\dots,n,
\end{equation}
where we choose a local ordering of the roots $z_j(a,\lambda)$ on the locus where they are
distinct and nonzero.

Fix $a$ and use $z$ as a coordinate on the curve in Remark~\ref{rem:spgeometry}. On the chosen
branch, the coordinate $q=P_\beta(z)$ gives $\omega_a=z^{\beta-1}dz=dq$. If we use $\lambda$ as a
coordinate on the branch $z_j(a,\lambda)$, the same form is given by
$\omega_a=\partial_\lambda P_j\,d\lambda=z_j^{\beta-1}\lambda'(z_j)^{-1}\,d\lambda$. The latter
expression is singular at the critical points of $\lambda$. At any such point the form $\omega_a$
is holomorphic and nonzero in the coordinate $z$, provided $z\neq0$.

To obtain flat coordinates for the natural connection, we expand $P_\beta$ at $\lambda=\infty$.
With our choice of normalization, the coefficients of the expansion are given by the twisted
residues
{
\begin{equation}\label{eq:u}
u_k=\frac{1}{s_k}\res_{z=\infty}\lambda^{s_k}\,z^{\beta-1}dz=a_k+\dots,\qquad
s_k=\frac{k-\beta}{n},\quad k=1,\dots,n ,
\end{equation}
}
where $\res_{z=\infty}f\,dz$ denotes the coefficient of $z^{-1}$ in the expansion of $f$ at
infinity. By our choice of $s_k$, only integral powers of $z$ occur after multiplication by
$z^{\beta-1}$, so that $u_k=\frac1{s_k}[z^{-k}](\lambda/z^n)^{s_k}$. The coefficient
$[z^{-k}](\lambda/z^n)^{s}=g_k(a,s)$ is a polynomial in $a$ and $s$ which vanishes at $s=0$, as
we shall see in the proof of Lemma~\ref{lem:phi}. It follows that $u_k$ is a polynomial in $a$
and $\beta$. When $s_k=0$, that is, when $\beta=k$, its value is the coefficient of $z^{-k}$ in
$\log(\lambda/z^n)$.

\begin{remark}\label{rem:normalised}
In the notation of Definition~\ref{def:sp} the variable is $P_\beta$. If $\beta\ne0$, the
superpotential in this variable is $\Lambda(P_\beta,a)=\lambda((1+\beta P_\beta)^{1/\beta},a)$,
where we choose a branch locally. For $\beta=0$ we obtain $\Lambda(P_0,a)=\lambda(\exp P_0,a)$.
Using $\partial_{P_\beta}=z^{1-\beta}\partial_z$, we find
\[
z\lambda'=(1+\beta P_\beta)\partial_{P_\beta}\Lambda .
\]
This explains the occurrence of the polynomial $z\lambda'$ in the algebra and of the marked point
$z=0$ among its roots. For $\beta\ne0$ the change of variable to $P_\beta+1/\beta$ brings the
factor on the right into homogeneous form. The value at the marked point is $a_n$.
\end{remark}
}

\subsubsection{The theorem}

{
To compare the products and the connections, we express them in Euclidean coordinates. Let
$z_1,\ldots,z_n$ denote the roots of $\lambda(z)$ and choose $p_j$ locally with $p_j^2=z_j$.
Since the basic invariants are $a_k=(-1)^ke_k(p_1^2,\dots,p_n^2)$, the quotient map $p\mapsto a$
is locally biholomorphic away from the reflecting hyperplanes. In particular, $M^{\mathrm{reg}}$
is the space of polynomials with distinct nonzero roots.

The hyperplanes $p_i\mp p_j=0$ correspond to the long roots and those given by $p_j=0$ to the
short roots. These are the two orbits. To write sums over both orbits, we set
$\alpha_H=p_i\mp p_j$, $v_H=\partial_i\mp\partial_j$, and $w_H=1$ on the first orbit, and
$\alpha_H=p_j$, $v_H=\partial_j$, and $w_H=2$ on the second. With this choice we have
$\sum_Hw_H\,\alpha_H\otimes v_H=2n\,\Id$. The following theorem describes the bi-flat structures
for $B_n$ in terms of the superpotential.
\begin{theorem}\label{thm:Bn}
Let $n\ge2$, $\beta\in\mathbb C$ and $\nu_0=\beta/n$.
The superpotential $\lambda(z)=z^n+a_1z^{n-1}+\dots+a_n$ with
$\omega=z^{\beta-1}dz$ determines a polynomial flat
structure $(\nabla,\circ,e,E)$ on $M$ and a bi-flat structure
$(\nabla,\circ,e,\nabla^*,*,E)$ on $M^{\mathrm{reg}}$.
\begin{enumerate}
\item[(i)] The product $\circ$ is induced by multiplication in $\mathbb C[z]/(z\lambda')$ under the
identification $\partial_{a_i}\mapsto z^{n-i}$. Its unit is $e=\partial_{a_n}$ and its Euler
field is $E=\sum_i\frac ina_i\partial_{a_i}=\frac1{2n}\sum_jp_j\partial_{p_j}$. Neither the
product nor these fields depend on $\beta$.
\item[(ii)] In the coordinates $p$ the dual product $X*Y:=E^{-1}\circ X\circ Y$ takes the form
{
\[
X*Y=\sum_Hw_H\,\frac{\alpha_H(X)\,\alpha_H(Y)}{\alpha_H(p)}\,v_H ,
\]
}
which is the dual product of the Frobenius manifold $B_n$ \cite{Dalm}.
\item[(iii)] Let $\tilde\nabla$ denote the connection for which $w_j=P_\beta(z_j)$ are local flat coordinates
on the regular locus. When $\beta=0$ these coordinates are $w_j=\log z_j$. The connection
$\nabla^*=\tilde\nabla-\nu_0\,*$ is invariant under $B_n$ and takes the form studied by Dunkl and
Kohno,
{
\[
\nabla^*_XY=\partial_XY-\sum_H\kappa_H\,w_H\,\frac{\alpha_H(X)\,\alpha_H(Y)}{\alpha_H(p)}\,v_H ,\qquad
\kappa_{\rm long}=\frac\beta n,\quad\kappa_{\rm short}=\frac12-\frac{(n-1)\beta}n ,
\]
}
where the two weights correspond to the two orbits of reflecting hyperplanes. Moreover,
$\nabla^*E=0$ and $(\nabla^*,*,E)$ is a flat F-manifold.
\item[(iv)] The connection $\nabla$ has flat coordinates $u_k$ given by \eqref{eq:u}. These are polynomial
coordinates on the entire quotient. On $M^{\mathrm{reg}}$ we have
$\nabla_XY=\nabla^*_XY-\nabla^*_{X\circ Y}e$. The unit and Euler field of the flat F-manifold
$(\nabla,\circ,e)$ are given by $e=\partial_{u_n}$ and $E=\sum_k\frac knu_k\partial_{u_k}$. In
the coordinates $a$ the Christoffel symbols of $\nabla$ satisfy
\begin{equation}\label{eq:chra}
\sum_{l=1}^n\Gamma^l_{ij}(a)\,z^{n-l}
=-\Big(z\frac{d}{dz}+\beta\Big)A_{ij}(z),
\end{equation}
where $A_{ij}$ denotes the quotient in the division of $z^{2n-i-j}$ by $z\lambda'$. In
particular, the Christoffel symbols are polynomial in $a$ and affine in $\beta$.
\item[(v)] The two flat structures form a bi-flat F-manifold $(\nabla,\circ,e,\nabla^*,*,E)$ and satisfy
\eqref{eq:dual}. Locally, $\lambda$ considered as a function of $P_\beta$ is a superpotential of
weight $\nu_0$ and constant $\kappa=1/n$ for each of the branches $w_j$. For $\beta\ne0$ the
change of variable to $P_\beta+1/\beta$ gives the homogeneous normalization $\kappa=0$. On the
locus where $z\lambda'(z)$ has $n$ distinct zeros, the values of $\lambda$ at its $n-1$ nonzero
critical points, together with the marked value $\lambda(0)=a_n$, form a system of canonical
coordinates.
\end{enumerate}
When $\beta=\frac12$ the two weights coincide with $\frac1{2n}$ and $\tilde\nabla$ is flat in the
coordinates $p$. In this way we recover the Frobenius manifold $B_n$ with Dubrovin's
superpotential $\lambda(p^2)$ and primary form $\omega\propto dp$.
\end{theorem}

\begin{proof}
We compute the two products in Steps~1 and~2. In Step~3 we define the connection by
\eqref{eq:chra}, which gives (S3) by construction. We then use Theorem~\ref{thm:flat} to prove
flatness and the other properties of a flat F-manifold with linear Euler field. An integration by
parts shows that the twisted residues in \eqref{eq:u} are flat coordinates. In Step~4 we apply
Proposition~\ref{prop:dualflat} to construct the dual structure. The zeros of the superpotential
then allow us to identify its dual connection with the connection of Dunkl and Kohno. This
construction also proves flatness. A proof using the reflection arrangement can be found in
\cite{CHL}. In Step~5 we identify the superpotential and the canonical coordinates.

\emph{Step 1. Product, unit, Euler field.} Since $z\lambda'$ has degree $n$, the polynomials
$\partial_{a_i}\lambda=z^{n-i}$, for $i=1,\dots,n$, represent a basis of
$\mathbb C[z]/(z\lambda')$. Using this basis to identify $T_aM$ with the quotient, we define the
product by $(X\circ Y)\lambda\equiv X\lambda\,Y\lambda
\pmod{z\lambda'}$. The product is commutative and associative. Its unit is $e=\partial_{a_n}$
because $\partial_{a_n}\lambda=1$. With the degrees $\deg z=\frac1n$ and $\deg a_i=\frac in$,
every monomial of $\lambda$ has degree $1$. Thus
\[
E\lambda+\tfrac1n\,z\lambda'=\lambda,\qquad\text{i.e.}\qquad
E\lambda=\lambda-\tfrac1n\,z\lambda'\equiv\lambda\pmod{z\lambda'},
\]
which shows that $E\circ$ acts on the quotient by multiplication by $\lambda$. This endomorphism
is invertible precisely when $\lambda$ and $z\lambda'$ have no common zero, that is, when
$\lambda$ has simple roots and $a_n\ne0$.

We also obtain $\mathrm{Lie}_E\circ=\circ$ from the grading. Indeed, in the basis
$\partial_{a_i}$ we have
$(\mathrm{Lie}_E\,c)^k_{ij}=E(c^k_{ij})-(\deg a_k-\deg a_i-\deg a_j)c^k_{ij}$. In the division
$\partial_{a_i}\lambda\,\partial_{a_j}\lambda
=c^k_{ij}\,\partial_{a_k}\lambda+q_{ij}\,z\lambda'$, the degrees of $\partial_{a_i}\lambda$ and
$z\lambda'$ are $1-\deg a_i$ and $1$, respectively. Comparing the degrees on both sides, we
obtain $1+\deg a_k-\deg a_i-\deg a_j$ as the degree of $c^k_{ij}$. Finally, from $z_j=p_j^2$ we
obtain $\deg p_j=\frac1{2n}$, so that $E=\frac1{2n}\sum_jp_j\partial_{p_j}$ in the coordinates
$p$.

\emph{Step 2. Dual product.} In the quotient algebra the relation $E\circ(X*Y)=X\circ Y$ takes
the form $\lambda\cdot(X*Y)\lambda\equiv X\lambda\,Y\lambda\pmod{z\lambda'}$. On
$M^{\mathrm{reg}}$ it characterizes $(X*Y)\lambda$ as the unique polynomial of degree $\le n-1$
representing $X\lambda\,Y\lambda\,\lambda^{-1}$. In particular, when $X\lambda\,Y\lambda$ is
divisible by $\lambda$ and the ordinary quotient has degree $\le n-1$, that quotient itself is
the required representative.

We first apply this observation to the product of two distinct coordinate fields associated to
the roots. Away from the discriminant the roots $z_1,\dots,z_n$ are local coordinates.
Differentiating $\lambda=\prod_k(z-z_k)$ gives
{
\[
\partial_{z_i}\lambda=-\prod_{k\neq i}(z-z_k)=-\frac{\lambda}{z-z_i},
\]
}
which is a polynomial of degree $n-1$. When $i\ne j$ the two derivatives omit different factors,
so that
{
\[
\partial_{z_i}\lambda\;\partial_{z_j}\lambda\;\lambda^{-1}=\frac{\lambda}{(z-z_i)(z-z_j)}=\prod_{k\neq i,j}(z-z_k)
\]
}
is a polynomial of degree $n-2$ and consequently represents
$(\partial_{z_i}*\partial_{z_j})\lambda$. The partial fraction identity
{
\[
\frac{1}{(z-z_i)(z-z_j)}=\frac1{z_i-z_j}\Big(\frac1{z-z_i}-\frac1{z-z_j}\Big),
\]
}
multiplied by $\lambda$, gives
{
\[
(\partial_{z_i}*\partial_{z_j})\lambda=\frac1{z_i-z_j}\Big(\frac{\lambda}{z-z_i}-\frac{\lambda}{z-z_j}\Big)
=-\frac{\partial_{z_i}\lambda-\partial_{z_j}\lambda}{z_i-z_j},
\qquad\text{i.e.}\qquad
\partial_{z_i}*\partial_{z_j}=-\frac{\partial_{z_i}-\partial_{z_j}}{z_i-z_j},
\]
}
since the action of a vector field on $\lambda$ determines it uniquely. We now express this
product in the coordinates in which the reflecting hyperplanes are linear. The relation
$z_i=p_i^2$ gives $\partial_{p_i}=2p_i\,\partial_{z_i}$, so that for $i\ne j$ we obtain
{
\[
\partial_{p_i}*\partial_{p_j}=4p_ip_j\,\partial_{z_i}*\partial_{z_j}
=-\frac{4p_ip_j}{p_i^2-p_j^2}\Big(\frac{\partial_{p_i}}{2p_i}-\frac{\partial_{p_j}}{2p_j}\Big)
=-\frac{2\,(p_j\,\partial_{p_i}-p_i\,\partial_{p_j})}{p_i^2-p_j^2}.
\]
}

Denote by $F(X,Y)$ the sum in (ii). If $i\ne j$, only the forms $p_i-p_j$ and $p_i+p_j$
contribute to $F(\partial_{p_i},\partial_{p_j})$, since a nonzero term requires
$\alpha_H(\partial_{p_i})\,\alpha_H(\partial_{p_j})\ne0$. Adding their contributions gives
\[
\begin{gathered}
\frac{(1)(-1)}{p_i-p_j}\,(\partial_{p_i}-\partial_{p_j})+\frac{(1)(1)}{p_i+p_j}\,(\partial_{p_i}+\partial_{p_j}) \\
=\frac{-(p_i+p_j)(\partial_{p_i}-\partial_{p_j})+(p_i-p_j)(\partial_{p_i}+\partial_{p_j})}{p_i^2-p_j^2} \\
=-\frac{2\,(p_j\,\partial_{p_i}-p_i\,\partial_{p_j})}{p_i^2-p_j^2},
\end{gathered}
\]
in agreement with the expression for the mixed product.

It remains to determine the diagonal products. For this we use the unit. Since
$E*Y=E^{-1}\circ E\circ Y=Y$ and $E=\frac1{2n}\sum_kp_k\partial_{p_k}$, we have
{
\[
\partial_{p_j}=E*\partial_{p_j}=\frac1{2n}\Big(p_j\,\partial_{p_j}*\partial_{p_j}+\sum_{k\neq j}p_k\,\partial_{p_k}*\partial_{p_j}\Big),
\]
}
which determines $\partial_{p_j}*\partial_{p_j}$ in terms of the mixed products whenever
$p_j\ne0$. We can apply the same argument to $F$. Indeed, from
$\alpha_H(E)=\frac1{2n}\alpha_H(p)$ and $\sum_Hw_H\,\alpha_H\otimes v_H=2n\,\Id$ we obtain
{
\[
F(E,Y)=\frac1{2n}\sum_Hw_H\,\alpha_H(Y)\,v_H=Y .
\]
}
Thus the two symmetric products agree on the diagonal as well as on the mixed terms. The
expression for $*$ is independent of $\beta$ and at $\beta=\frac12$ coincides with the dual
product of the Frobenius manifold $B_n$ \cite{Dalm}.

{
\emph{Step 3. Natural connection.} We define $\nabla'$ by prescribing its Christoffel symbols in
the coordinates $a_1,\dots,a_n$ according to \eqref{eq:chra}. To see that this defines a torsion
free connection, observe that $z^{2n-i-j}$ has degree $\le2n-2$. Thus its quotient $A_{ij}$ by
$z\lambda'$ has degree $\le n-2$. The same bound applies to $\big(z\frac d{dz}+\beta\big)A_{ij}$,
which we can therefore express as a linear combination of $z^{n-l}$ for $l=2,\dots,n$. The
symmetry follows from $A_{ij}=A_{ji}$. Moreover, since $z\lambda'$ has constant leading
coefficient $n$, the Christoffel symbols are polynomial in $a$ and affine in $\beta$. In Step~4
we identify $\nabla'$ with the connection $\nabla$ in (iv).

\emph{(a) Conditions (S1) to (S3).} We use $P_\beta$ as the variable. Differentiating with
$P_\beta$ fixed is therefore equivalent to differentiating with $z$ fixed. In particular,
$\partial_{P_\beta}=z^{1-\beta}\partial_z$ and $\partial_{P_\beta}\Lambda=z^{1-\beta}\lambda'$,
also when $\beta=0$. Condition (S1) reads $\partial_{a_n}\lambda=1$. To verify (S2) with
$\nu_0=\frac\beta n$ and $\kappa=\frac1n$, we use the identity
$E\lambda=\lambda-\frac1nz\lambda'$ from Step~1 together with
\[
(\nu_0P_\beta+1/n)\,\partial_{P_\beta}\lambda=\frac{z^\beta}{n}\,z^{1-\beta}\lambda'=\frac1n\,z\lambda' .
\]
Both conditions are independent of the connection. To verify (S3), we use \eqref{eq:defconn}.
Since the superpotential is linear in the coordinates $a$, we have
$\partial_{a_i}\partial_{a_j}\lambda=0$. Moreover,
\begin{equation}\label{eq:S3z}
R_{ij}=\partial_{a_i}\lambda\,\partial_{a_j}\lambda-c^l_{ij}\,\partial_{a_l}\lambda=z^{2n-i-j}-c^l_{ij}z^{n-l}=z\lambda'\,A_{ij},
\end{equation}
since, by the definition of the product, $c^l_{ij}z^{n-l}$ is the remainder in the division of
$z^{2n-i-j}$ by $z\lambda'$. We obtain
\[
B_{ij}=\frac{R_{ij}}{\partial_{P_\beta}\Lambda}=A_{ij}\,z^\beta,\qquad
\partial_{P_\beta}B_{ij}=z^{1-\beta}\partial_z\big(z^\beta A_{ij}\big)=\Big(z\frac d{dz}+\beta\Big)A_{ij},
\]
so that \eqref{eq:defconn}, which here takes the form
$\Gamma^l_{ij}\,z^{n-l}=-\partial_{P_\beta}B_{ij}$, is satisfied by the symbols in
\eqref{eq:chra}. It follows that $\nabla'$ is the connection defined by $\lambda$ in
Theorem~\ref{thm:flat} and that $\lambda$ satisfies (S3) for this connection.

\emph{(b) Flatness.} If $a$ lies outside the discriminant of $\lambda(z)-\lambda$, the equation
$\lambda(z)=\lambda$ has $n$ simple roots $z_j(a,\lambda)$ which give local coordinates. On the
locus where these roots are nonzero, the branches $P_j=P_\beta(z_j)$ in \eqref{eq:ansatz} have
independent differentials $dP_j=z_j^{\beta-1}dz_j$. Together with (a), this verifies the
hypotheses of Theorem~\ref{thm:flat}. Thus $\nabla'$ is flat, $\nabla'e=0$, and
$\nabla'\nabla'E=0$. All three conditions are polynomial identities in $a$, so they hold on the
whole of $M$ for every $\beta$.

\emph{(c) Flat coordinates.} Let us first compute the Jacobian matrix of the change of variables
$a\mapsto u$ defined by the twisted residues $u$ in \eqref{eq:u}.

\begin{lemma}\label{lem:phi}
For $m=1,\dots,n$, let $s_m=\frac{m-\beta}n$ be the exponents in \eqref{eq:u}. For every
polynomial $f$ of degree $\le n-1$, define $\varphi_m(f)=\res\,\lambda^{s_m-1}f\,z^{\beta-1}dz$.
Let $\Phi(a)$ denote the $n\times n$ matrix whose entry $(m,i)$ is
\[
\Phi_{mi}:=\varphi_m(z^{n-i})=\res\,\lambda^{s_m-1}z^{n-i}\,z^{\beta-1}dz .
\]
The following identities hold.
\begin{enumerate}
\item[(i)] The matrix $\Phi$ is the Jacobian matrix of $a\mapsto u$. More precisely,
$\partial u_m/\partial a_i=\Phi_{mi}$, including when $s_m=0$. It is lower unitriangular, with
$\Phi_{mi}=0$ for $i>m$ and $\Phi_{mm}=1$. In particular, $u_m=a_m+f_m(a_1,\dots,a_{m-1})$ with
$f_m$ polynomial. Thus $u_1,\dots,u_n$ are polynomial coordinates on $M$ and the inverse change
of coordinates is polynomial.
\item[(ii)] $\varphi_m(1)=\delta_{mn}$. In particular, $e=\partial_{a_n}=\partial_{u_n}$.
\end{enumerate}
\end{lemma}

\begin{proof}[Proof of Lemma~\ref{lem:phi}]
To compute the residues, we set $x=z^{-1}$ and expand at infinity. For $|z|$ sufficiently large
we have
{
\[
\lambda=z^n\big(1+a_1x+\dots+a_nx^n\big),\qquad
\lambda^{s}=z^{ns}\big(1+a_1x+\dots+a_nx^n\big)^{s},
\]
}
where we denote the coefficient of $x^r$ in the binomial series by $g_r(a,s)$. More explicitly,
\[
\big(1+a_1x+\dots+a_nx^n\big)^s=\sum_{k\ge0}\binom sk\big(a_1x+\dots+a_nx^n\big)^k=\sum_{r\ge0}g_r(a,s)\,x^r,
\]
with coefficients
\[
g_r(a,s)=\sum_{k=0}^r\binom sk\sum_{\substack{i_1+\dots+i_k=r\\1\le i_j\le n}}a_{i_1}\cdots a_{i_k},
\]
where $\binom sk=\frac{s(s-1)\cdots(s-k+1)}{k!}$. For instance, $g_0=1$, $g_1=s\,a_1$, and
$g_2=s\,a_2+\binom s2a_1^2$. The series converges for $|x|$ sufficiently small, and each
$g_r(a,s)$ is a polynomial in $a_1,\dots,a_n$ and $s$. The only factor in which a noninteger
power of $z$ can occur is $z^{ns}$.

We now compute $\Phi_{mi}$ by substituting $s=s_m-1$ in this expansion. Since $ns_m=m-\beta$, the
leading exponent in $\lambda^{s_m-1}z^{n-i}z^{\beta-1}$ is
{
\[
(m-\beta-n)+(n-i)+(\beta-1)=m-i-1 ,
\]
}
which gives the Laurent series
\[
\lambda^{s_m-1}z^{n-i}\,z^{\beta-1}=z^{m-i-1}\sum_{r\ge0}g_r(a,s_m-1)\,z^{-r}
\]
in integer powers. The residue is $g_{m-i}(a,s_m-1)$ for $m-i\ge0$ and vanishes for $i>m$. When
$i=m$, the coefficient is $g_0=1$. Thus $\Phi$ is triangular with diagonal entries $1$. In
particular, taking $i=n$ gives $\res\,\lambda^{s_m-1}z^{\beta-1}dz=\varphi_m(1)=\delta_{mn}$.

Since the coefficients are polynomial in $a$, we may differentiate the residue with respect to
$a_i$. From $\partial_{a_i}\big(\tfrac1{s_m}\lambda^{s_m}\big)
=\lambda^{s_m-1}\partial_{a_i}\lambda=\lambda^{s_m-1}z^{n-i}$ we obtain
{
\[
\frac{\partial u_m}{\partial a_i}=\res\,\lambda^{s_m-1}z^{n-i}z^{\beta-1}dz=\Phi_{mi},
\]
}
which also holds for $s_m=0$. Thus $\Phi$ is the Jacobian matrix of $a\mapsto u$. Its triangular
form gives $u_m=a_m+f_m(a_1,\dots,a_{m-1})$, where $f_m$ is a polynomial. Its constant term
vanishes because at $a=0$ the form is $z^{m-1}dz$ and its residue is zero for $m\ge1$. Starting
with $a_1=u_1$ and $a_2=u_2-f_2(u_1)$, we can solve recursively for each coefficient. This gives
a polynomial inverse.

Finally, $du_m(e)=\partial u_m/\partial a_n=\varphi_m(1)=\delta_{mn}$, so $e=\partial_{a_n}=\partial_{u_n}$.
\end{proof}

We can now prove that the twisted residues $u_m$ are flat coordinates of $\nabla'$. For this we
verify the equations $\partial_{a_i}\partial_{a_j}u_m=\Gamma^l_{ij}\,\partial_{a_l}u_m$. Set
$s=s_m$ and $\omega=z^{\beta-1}dz$. By Lemma~\ref{lem:phi} we have
$\partial_{a_l}u_m=\varphi_m(z^{n-l})$, including when $s_m=0$. From \eqref{eq:chra} and the
linearity of $\varphi_m$ we obtain
\[
\Gamma^l_{ij}\,\partial_{a_l}u_m=\varphi_m\big(\Gamma^l_{ij}z^{n-l}\big)=-\varphi_m\big(z^{1-\beta}\partial_z(z^\beta A_{ij})\big)
=-\res\,\lambda^{s-1}\,\partial_z\big(z^\beta A_{ij}\big)\,dz .
\]
All the forms in this calculation have integer exponents at $z=\infty$, since
$\lambda^{s-1}z^\beta=z^{m-n}(1+O(z^{-1}))$. Integrating by parts, we obtain
\[
-\res\,\lambda^{s-1}\,\partial_z\big(z^\beta A_{ij}\big)\,dz=(s-1)\,\res\,\lambda^{s-2}\lambda'\,z^\beta A_{ij}\,dz=(s-1)\,\res\,\lambda^{s-2}R_{ij}\,\omega .
\]
The terms of $R_{ij}$ which involve the structure constants give zero upon taking the residue.
Indeed,
\[
(s-1)\,c^l_{ij}\,\res\,\lambda^{s-2}z^{n-l}\,\omega=c^l_{ij}\,\partial_{a_l}\big(\res\,\lambda^{s-1}\omega\big)=c^l_{ij}\,\partial_{a_l}\delta_{mn}=0
\]
by Lemma~\ref{lem:phi}. Therefore
\[
\Gamma^l_{ij}\,\partial_{a_l}u_m=(s-1)\,\res\,\lambda^{s-2}z^{n-i}z^{n-j}\,\omega
=\partial_{a_i}\big(\res\,\lambda^{s-1}z^{n-j}\omega\big)=\partial_{a_i}\partial_{a_j}u_m .
\]
By Lemma~\ref{lem:phi}, the functions $u_m$ give polynomial coordinates on $M$ and
$e=\partial_{u_n}$. Moreover, the residue formula gives $u_k$ degree $\frac kn$. It follows that
$E=\sum_k\frac knu_k\partial_{u_k}$ and $[e,E]=e$.

\emph{(d) Flat F-manifold.} Part (c) gives, in the flat coordinates $u$, the expression
$N=\nu_0-\nabla'E=\mathrm{diag}\big(\frac{\beta-k}n\big)$. This endomorphism is invertible for
$\beta\notin\{1,\dots,n\}$, in particular for $\beta=0$. For these values part (ii) of
Theorem~\ref{thm:flat} gives the symmetry of $\nabla'_lc^m_{hk}$. Since the structure constants
are independent of $\beta$ and the Christoffel symbols are affine in $\beta$, the symmetry
extends to every $\beta$. Together with Step~1, this proves that $(\nabla',\circ,e)$ is a flat F-manifold with linear Euler field $E$ and that the structure is polynomial on all of $M$.

\emph{Step 4. The bi-flat structure and the dual connection.} We apply
Proposition~\ref{prop:dualflat} on $M^{\mathrm{reg}}$. It follows that the connection
\[
D_XY:=\nabla'_XY-\nabla'_{X*Y}E
\]
is flat and that $(D,*,E)$ is a flat F-manifold with $DE=0$. Moreover, $(\nabla',\circ,e,D,*,E)$
is a bi-flat F-manifold. We now show that $D$ coincides with the connection $\nabla^*$ in (iii).

By Step~3(a), $\lambda$ satisfies (S1) to (S3) for this bi-flat structure. Applying part (iv) of
Theorem~\ref{thm:char} at $\lambda=0$, we find that the functions
$w_j=P_j(\cdot,0)=P_\beta(z_j)$, where $z_j$ are the roots of $\lambda$, form a system of flat
coordinates for $\Gamma(0,\nu_0)=D+\nu_0\,*$. Away from $p_j=0$, once we choose local branches of
the logarithms, these are also flat coordinates of $\tilde\nabla$. Since the Christoffel symbols
of both connections vanish in the same coordinates, we have $D+\nu_0\,*=\tilde\nabla$, or
equivalently
\[
D=\tilde\nabla-\nu_0\,*=\nabla^* .
\]
In particular, $(\nabla^*,*,E)$ is a flat F-manifold with $\nabla^*E=0$ for every $\beta$,
including $\beta=0$. We next determine the expression for $\nabla^*$ in the coordinates $p$.}

\emph{The form of the dual connection.} To express $\tilde\nabla$ in Euclidean coordinates, we
use its flat coordinates $w_j$. Choosing local branches of the logarithms, we obtain the
functions $w_j=P_\beta(p_j^2)$ as coordinates away from $p_j=0$. Each derivative
$\partial_{p_l}w_l$ is nonzero and satisfies
\[
\partial^2_{p_l} w_l=\frac{2\beta-1}{p_l}\,\partial_{p_l}w_l .
\]
Each $w_l$ depends only on $p_l$. Thus the equation
$\partial_i\partial_jw_l=\tilde\Gamma^k_{ij}\partial_kw_l$ for a flat coordinate gives
$\tilde\Gamma^l_{ll}=(2\beta-1)/p_l$, and all the other Christoffel symbols vanish. We obtain
{
\[
\tilde\nabla_XY=\partial_XY+(2\beta-1)\sum_j\frac{dp_j(X)\,dp_j(Y)}{p_j}\,\partial_{p_j}
=\partial_XY+\frac{2\beta-1}2\,(X*Y)_{\rm short},
\]
}
where $(X*Y)_{\rm short}=2\sum_j\frac{dp_j(X)\,dp_j(Y)}{p_j}\,\partial_{p_j}$ is the contribution
in (ii) from the hyperplanes $p_j=0$, which have $w_H=2$. Subtracting $\nu_0\,X*Y$, we obtain the
connection of Dunkl and Kohno with weights $\kappa_{\rm long}=\nu_0=\frac\beta n$ and
$\kappa_{\rm short}=\nu_0-\frac{2\beta-1}2=\frac12-\frac{(n-1)\beta}n$. The weights are constant
on the respective orbits, as required by invariance under $B_n$.

To verify the normalization, recall that $\alpha_H(E)=\frac1{2n}\alpha_H(p)$. For the two orbits
we have
{
\[
\sum_{\rm long}w_H\,\alpha_H\otimes v_H=2(n-1)\,\Id,\qquad
\sum_{\rm short}w_H\,\alpha_H\otimes v_H=2\,\Id
\]
}
To see the first equality, fix $i<j$. The two hyperplanes $p_i-p_j=0$ and $p_i+p_j=0$ together
contribute $(X^i-X^j)(\partial_i-\partial_j)+(X^i+X^j)(\partial_i+\partial_j)
=2X^i\partial_i+2X^j\partial_j$, and each index belongs to $n-1$ pairs. We therefore obtain
{
\[
\nabla^*_XE=\frac1{2n}\,X-\sum_H\kappa_Hw_H\,\frac{\alpha_H(X)\,\alpha_H(E)}{\alpha_H(p)}\,v_H
=\frac1{2n}\big(1-2(n-1)\kappa_{\rm long}-2\kappa_{\rm short}\big)X=0 .
\]
}

In this construction the flatness of $\nabla^*$ already follows from the duality. For a Coxeter
group it also follows from the classical result on connections of Dunkl and Kohno with invariant
weights \cite{CHL}.

To finish Step~4, we identify $\nabla'$ with the connection $\nabla$ in (iv). Indeed, from
$\nabla^*_XY=\nabla'_XY-\nabla'_{X*Y}E$ and $\nabla'e=0$ we obtain
$\nabla^*_Ze=-\nabla'_{Z*e}E=-\nabla'_{E^{-1}\circ Z}E$. Substituting $Z=X\circ Y$ and using
$E^{-1}\circ X\circ Y=X*Y$, we find
\[
\nabla^*_{X\circ Y}e=-\nabla'_{X*Y}E=\nabla^*_XY-\nabla'_XY
\]
(compare Proposition~3.7 of \cite{KMS}). Thus
$\nabla'_XY=\nabla^*_XY-\nabla^*_{X\circ Y}e=\nabla_XY$, and Step~3 gives all the assertions in
(iv).

\emph{Step 5. Superpotential and canonical coordinates.} By Step~4,
$(\nabla,\circ,e,\nabla^*,*,E)$ is a bi-flat F-manifold on $M^{\mathrm{reg}}$. Its connections
satisfy \eqref{eq:dual}, and the compatibility of the two flat structures follows from
Lemma~\ref{lem:AHE}.

We apply the construction of Section~\ref{sec:defsp} to the functions $w_j$. They are flat
coordinates of $\Gamma(0,\nu_0)=\nabla^*+\nu_0c^*$ and satisfy $E(w_j)=\nu_0w_j+1/n$. Using
$\partial_{a_n}\lambda=1$ and differentiating the equation for a root, we obtain
$\partial_ew_j=-z_j^{\beta-1}/\lambda'(z_j)\ne0$ on the regular locus. The flow of $-e$ subtracts
$\lambda$ from the constant term $a_n$. Thus the equation for a zero becomes the equation for a
root at level $\lambda$, which gives $w_j(\phi_{-\lambda}(a))=P_\beta(z_j(a,\lambda))$. Taking
the local inverse, we recover the superpotential in the variable $P_\beta$ with constant
$\kappa=1/n$.

To identify the canonical coordinates, let us evaluate the quotient algebra at a zero $z_*$ of
$z\lambda'$. The resulting character $X\mapsto(X\lambda)(z_*)$ of $T_aM$ is given by
$d\big(\lambda(z_*)\big)(X)$. Indeed, at a nonzero root of $z\lambda'$ we have $\lambda'(z_*)=0$,
while at $z_*=0$ the claim follows from $\lambda(0)=a_n$ and $X\lambda(0)=X(a_n)$. Wherever
$z\lambda'$ has simple roots, the $n$ characters are independent. The corresponding values of the
polynomial therefore give canonical coordinates. At $z=0$ we obtain the marked value $a_n$.
Whether the marked point is critical with respect to $P_\beta$ depends on $\beta$. For example,
at $\beta=1$ we have $P_\beta=z-1$ and $\lambda'(0)=a_{n-1}$, so that $z=0$ is a regular point
whenever $a_{n-1}\ne0$.
{
At each of the $n-1$ nonzero roots we have $dP_\beta/dz=z^{\beta-1}\ne0$, so each gives a simple
critical point in the primary coordinate. Its character is described in
Remark~\ref{rem:simplecrit}, and by Corollary~\ref{cor:lame} the corresponding Lam\'e degree is
$\beta/n-\tfrac12$. The residue form obtained from these $n-1$ characters has rank $n-1$. For
$\beta\ne\tfrac12$, evaluation at the marked root gives the remaining character of the quotient
algebra, and the entire connection is given by Step~3. When $\beta=\tfrac12$, the marked root
also gives a simple critical point in $P_\beta$, since $z=(1+P_\beta/2)^2$ and $a_{n-1}\ne0$. In
this case all $n$ differentials occur in the residue formula in part (iii) of
Proposition~\ref{prop:christoffel}.} Finally, for $\beta=\frac12$ the corresponding branches give
$w_j=2p_j-2$ and both weights equal $\frac1{2n}$.
\end{proof}

\begin{remark}\label{rem:chrBn}
Formula \eqref{eq:chra} expresses the covariant form \eqref{eq:S3cov} of (S3) in the coordinates
$a$, as we saw in Step~3(a) of the proof. We used this formula in Step~3 to define the connection
and then obtained its flatness from Theorem~\ref{thm:flat}.

Alternatively, we may compute the Christoffel symbols from the flat coordinates themselves.
Differentiating \eqref{eq:u}, we obtain
\[
\Gamma^l_{ij}\res_{z=\infty}\lambda^{s_m-1}\partial_{a_l}\lambda\,\omega
=(s_m-1)\res_{z=\infty}\lambda^{s_m-2}
\partial_{a_i}\lambda\,\partial_{a_j}\lambda\,\omega ,
\qquad 1\le m\le n .
\]
The coefficient matrix on the left is the unitriangular matrix $\Phi$ in Lemma~\ref{lem:phi}. For
$n=2$ the only symbol which can be nonzero is $\Gamma^2_{11}=-\beta/2$, as one can also see from
$u_1=a_1$ and $u_2=a_2-\beta a_1^2/4$.
\end{remark}

Throughout this family the orbit space has the same product, unit, and Euler field as the
Frobenius manifold $B_n$. The parameter in the dual connection enters through the choice of the
member of the pencil for which the functions $w_j$ are flat. The Euclidean choice corresponds to
$\beta=1/2$, while in the whole family we have $\beta=n\,\kappa_{\rm long}$. At $\beta=0$ the
weight of the long orbit vanishes and, on consistent branches, the dual flat functions are
$\log z_j=2\log p_j$. The flat coordinates of the natural connection are polynomial for every
$\beta$. In particular, for $\beta=k$ we obtain $u_k$ as the coefficient of $z^{-k}$ in
$\log(\lambda/z^n)$.

It remains to determine whether these data give all the normalized connections invariant under
$B_n$.
\begin{theorem}\label{conj:Bn}
Fix $e$, $E$, and the dual product as in Theorem~\ref{thm:Bn}. For every $n\ge2$, each normalized
connection of Dunkl and Kohno invariant under $B_n$ on $M^{\mathrm{reg}}$ is the dual connection
of exactly one bi-flat structure with these data. In ranks $2$, $3$, and $4$ the resulting
families agree with those of \cite{ALrefl,Bn} up to a reparametrization.
\end{theorem}

\begin{proof}
Following \cite{ALrefl}, we consider the product of the Coxeter arrangement together with a
connection of Dunkl and Kohno on the quotient. The weights on the two orbits are
$\kappa_{\rm long}$ for $p_i\mp p_j=0$ and $\kappa_{\rm short}$ for $p_j=0$. With our
conventions, the condition $\nabla^*E=0$ reads $2(n-1)\kappa_{\rm long}+2\kappa_{\rm short}=1$
\cite{ALrefl}. Thus $\kappa_{\rm long}$ determines the connection. By parts (iii) and (v) of
Theorem~\ref{thm:Bn}, every value is realized by taking $\kappa_{\rm long}=\beta/n$. Once we fix
a connection, the unit, and the dual product, uniqueness follows from the reconstruction formula
after Lemma~\ref{lem:AHE}.
\end{proof}

Thus the same pair $\big(z^n+a_1z^{n-1}+\dots+a_n,\,z^{\beta-1}dz\big)$ describes the families
for $B_n$ in every rank.
}

\subsection{\texorpdfstring{The groups $G(m,1,n)$}{The groups G(m,1,n)}}\label{sec:Gm1n}

{
We can use the same space $M$ of coefficients for the quotient by $G(m,1,n)$. Indeed, for $m,n\ge2$ the
group acts on $p\in\mathbb C^n$ by permuting the coordinates and multiplying them by roots of
unity of order dividing $m$. Its basic invariants are therefore the elementary symmetric
functions of $z_j=p_j^m$. Thus the quotient again parametrizes the polynomials
$\lambda(z)=\prod_j(z-z_j)$ with $\deg z=\frac1n$. When $m=2$ we recover $G(2,1,n)=B_n$. On the
covering space, however, we must now use the reflecting hyperplanes of this group and the
corresponding derivatives of $z_j$ with respect to $p_j$.

There are two orbits of reflecting hyperplanes. The first consists of the hyperplanes
$p_i=\zeta p_j$ with $\zeta^m=1$, and the second consists of the hyperplanes $p_j=0$. For
$\alpha_H=p_i-\zeta p_j$ set $v_H=\partial_i-\zeta^{-1}\partial_j$ and $w_H=1$, and for
$\alpha_H=p_j$ set $v_H=\partial_j$ and $w_H=m$. With these conventions,
$\sum_Hw_H\,\alpha_H\otimes v_H=mn\,\Id$ and the construction on the quotient takes the following
form.
\begin{proposition}\label{prop:Gm1n}
Let $m,n\ge2$. In the coordinates given by the coefficients of the polynomial, the natural
structures in Theorem~\ref{thm:Bn} define polynomial structures on $\mathbb C^n/G(m,1,n)$. The
regular quotient consists of the polynomials with distinct nonzero roots. In the linear
coordinates $p_j$ on $V^{\mathrm{reg}}$ we have
{
\[
E=\frac1{mn}\sum_jp_j\partial_{p_j},\qquad
X*Y=\sum_Hw_H\,\frac{\alpha_H(X)\,\alpha_H(Y)}{\alpha_H(p)}\,v_H ,
\]
}
{
\[
\nabla^*_XY=\partial_XY-\sum_H\kappa_H\,w_H\,\frac{\alpha_H(X)\,\alpha_H(Y)}{\alpha_H(p)}\,v_H ,
\]
}
where the weights are $\kappa=\frac\beta n$ on the hyperplanes $p_i=\zeta p_j$ and
$\kappa=\frac1m-\frac{(n-1)\beta}n$ on the hyperplanes $p_j=0$. They satisfy
$m(n-1)\kappa_{p_i=\zeta p_j}+m\kappa_{p_j=0}=1$, which gives $\nabla^*E=0$. With this dual
product and with $e$ and $E$ fixed, the family of Theorem~\ref{thm:Bn} realizes every normalized
connection of Dunkl and Kohno invariant under $G(m,1,n)$. In the coordinate $p=z^{1/m}$ the
superpotential and the primary form are given by
{
\[
\Lambda=\lambda(p^m)=p^{mn}+a_1p^{m(n-1)}+\dots+a_n,\qquad \omega=z^{\beta-1}dz=m\,p^{m\beta-1}dp .
\]
}
At $\beta=\frac1m$ both weights $\kappa_H$ equal $\frac1{mn}$. In this case $\omega\propto dp$
and $\nabla^*+\nu_0c^*$ is flat in the linear coordinates $p$. We thus recover the natural Saito
structure on the orbit space of $G(m,1,n)$ \cite{KMS}.
\end{proposition}

\begin{proof}
We first compute the mixed products in the coordinates given by the roots. Since the formulas in
these coordinates depend only on the coefficients of the polynomial, Step~2 of the proof of
Theorem~\ref{thm:Bn} gives
$\partial_{z_i}*\partial_{z_j}=-(\partial_{z_i}-\partial_{z_j})/(z_i-z_j)$ for $i\ne j$. Using
$\partial_{p_i}=mp_i^{m-1}\partial_{z_i}$, we find
{
\[
\partial_{p_i}*\partial_{p_j}=-m\,\frac{p_j^{m-1}\partial_{p_i}-p_i^{m-1}\partial_{p_j}}{p_i^m-p_j^m}
=\sum_{\zeta^m=1}\frac{-\zeta\,(\partial_{p_i}-\zeta^{-1}\partial_{p_j})}{p_i-\zeta p_j},
\]
}
where the second equality follows from $\sum_\zeta(\xi-\zeta\eta)^{-1}=m\xi^{m-1}/(\xi^m-\eta^m)$
and $\sum_\zeta\zeta(\xi-\zeta\eta)^{-1}=m\eta^{m-1}/(\xi^m-\eta^m)$. The hyperplanes
$p_i=\zeta p_j$ are the only ones contributing to this mixed product, and for each of them
$\alpha_H(\partial_{p_i})\alpha_H(\partial_{p_j})=-\zeta$. As in the proof for the Coxeter group,
the unit determines the diagonal products. Indeed, the proposed formula gives
$E*Y=\frac1{mn}\sum_Hw_H\alpha_H(Y)v_H=Y$, so it agrees on the diagonal as well.

To compute the connection in the coordinates $p$, we use the flat coordinates $P_\beta(p_j^m)$ of
$\nabla^*+\nu_0c^*$ given in part (iii) of Theorem~\ref{thm:Bn}. At $\beta=0$ they have the form
$m\log p_j$ on consistent branches. For every $\beta$, dividing the second derivative with
respect to $p_j$ by the first gives $(m\beta-1)/p_j$. The connection is therefore given by
$\partial_XY+(m\beta-1)\sum_jp_j^{-1}dp_j(X)dp_j(Y)\partial_{p_j}$. Its difference from
$\partial$ is $\frac{m\beta-1}m$ times the contribution to $X*Y$ from the hyperplanes $p_j=0$.
Subtracting $\nu_0\,X*Y$ with $\nu_0=\frac\beta n$, we obtain the stated weights.

At $\beta=\frac1m$ the flat coordinates are affine functions of $p_j$. The degrees of $G(m,1,n)$
are $m,2m,\dots,nm$, so that $h=mn$ and $\nu_0=\frac1{mn}=\frac1h$. Thus the linear functions
$p_j$ are flat coordinates of $\nabla^*+\frac1h\,*$. By the characterization recalled in
Remark~\ref{rem:KMSKMS}, we obtain the KMS flat structure. As in Step~4 of the proof of
Theorem~\ref{thm:Bn}, flatness follows directly from the construction. For another proof in the
case of invariant weights, see \cite{CHL}.
\end{proof}

Thus $B_n$ and $G(m,1,n)$ give the same bi-flat F-manifolds when we use the coefficients of the
polynomial as coordinates. The primary form is proportional to $dp$ at the distinguished value
$\beta=1/m$. For $m=2$ we recover the Frobenius structure of $B_n$. For $m\ge3$ we obtain the
natural Saito structure of $G(m,1,n)$, which admits no compatible Frobenius metric \cite{KM}.
Among connections with one weight on each orbit, the condition $\nabla^*E=0$ again leaves one
parameter \cite{ALrefl}. The argument used in Theorem~\ref{conj:Bn} therefore gives the entire
family.
}

\subsection{\texorpdfstring{The case of $F_4$}{The case of F4}}\label{sec:F4}

{
The group $F_4=G_{28}$ has rank $4$, degrees $2,6,8,12$, and Coxeter number $h=12$. Each of its
two orbits of reflecting hyperplanes contains $12$ elements. Thus the family in
Conjecture~\ref{conj:AL} again depends on one parameter. We write the flat coordinates of
$\nabla$ as polynomials in $p_1,\dots,p_4$ depending on $c_1$. From these coordinates and the
corresponding vector potential, we identify the weights of the dual connection in
Proposition~\ref{prop:F4}.

All $24$ reflections have order two. To describe the reflecting hyperplanes, we write
$\alpha_s(p)=\sum_i(\alpha_s)_ip_i$ and use the following vectors of coefficients.
\begin{eqnarray*}
&&\alpha_1=\left[1,0,0,0\right],\, \alpha_2=\left[0,1,0,0\right],\, \alpha_3=\left[0,0,1,0\right],\,\alpha_4=\left[0,0,0,1\right],\\
&&\alpha_5=\left[1,1,0,0\right],\, \alpha_6=\left[1,0,1,0\right],\, \alpha_7=\left[1,0,0,1\right],\,\alpha_8=\left[0,1,1,0\right],\\
&&\alpha_9=\left[0,1,0,1\right],\, \alpha_{10}=\left[0,0,1,1\right],\, \alpha_{11}=\left[1,-1,0,0\right],\,\alpha_{12}=\left[1,0,-1,0\right],\\
&&\alpha_{13}=\left[1,0,0,-1\right],\, \alpha_{14}=\left[0,1,-1,0\right],\, \alpha_{15}=\left[0,1,0,-1\right],\,\alpha_{16}=\left[0,0,1,-1\right],\\
&&\alpha_{17}=\left[1,1,1,1\right],\, \alpha_{18}=\left[1,1,1,-1\right],\, \alpha_{19}=\left[1,-1,1,1\right],\,\alpha_{20}=\left[1,1,-1,1\right],\\
&&\alpha_{21}=\left[1,1,-1,-1\right],\, \alpha_{22}=\left[1,-1,1,-1\right],\, \alpha_{23}=\left[1,-1,-1,1\right],\,\alpha_{24}=\left[1,-1,-1,-1\right].\\
\end{eqnarray*}
The forms $\alpha_1,\dots,\alpha_4,\alpha_{17},\dots,\alpha_{24}$ define the orbit $\mathcal O_1$
corresponding to the short roots, while $\alpha_5,\dots,\alpha_{16}$ define the orbit
$\mathcal O_2$ corresponding to the long roots. Starting with basic invariants $U_1,\dots,U_4$,
we obtain the generalized Saito flat coordinates by the following triangular change of
coordinates.
{
\begin{eqnarray*}
u_1&=&U_1\\
u_2&=&U_2+c_1U_1^3\\
u_3&=&U_3+c_2U_1U_2+c_3U_1^4\\
u_4&=&U_4+c_4U_1^2U_3+c_5U_2^2+c_6U_2U_1^3+c_7U_1^6
\end{eqnarray*}
}
The coefficients in this change of coordinates are the following polynomials in $c_1$.
\begin{eqnarray*}
c_2&=& \frac{1}{4}+18c_1\\
c_3&=&9c_1^2-\frac{1}{576}\\
c_4&=&\frac{1}{8}+9c_1\\
c_5&=&27c_1+\frac{3}{8}\\
c_6&=&54c_1^2-\frac{1}{96}\\
c_7&=&-\frac{1}{480}c_1-\frac{9}{20}c_1^2+\frac{54}{5}c_1^3+\frac{1}{11520}
\end{eqnarray*}
We choose the basic invariants as follows.
\begin{eqnarray*}
U_1&=&p_1^2+p_2^2+p_3^2+p_4^2,\\
U_2&=& \frac{1}{6}p_1^4p_2^2+\frac{1}{6}p_1^4p_3^2+\frac{1}{6}p_1^4p_4^2+\frac{1}{6}p_2^4p_1^2-\frac{1}{2}p_1^2p_2^2p_3^2-\frac{1}{2}p_1^2p_2^2p_4^2+\frac{1}{6}p_3^4p_1^2-\frac{1}{2}p_1^2p_3^2p_4^2\\
&&+\frac{1}{6}p_4^4p_1^2+\frac{1}{6}p_2^4p_3^2+\frac{1}{6}p_2^4p_4^2+\frac{1}{6}p_3^4p_2^2-\frac{1}{2}p_2^2p_3^2p_4^2+\frac{1}{6}p_4^4p_2^2+\frac{1}{6}p_3^4p_4^2+\frac{1}{6}p_4^4p_3^2,\\
U_3&=& \frac{1}{12}p_1^4p_2^4-\frac{1}{12}p_1^4p_2^2p_3^2-\frac{1}{12}p_1^4p_2^2p_4^2+\frac{1}{12}p_1^4p_3^4-\frac{1}{12}p_1^4p_3^2p_4^2+\frac{1}{12}p_1^4p_4^4-\frac{1}{12}p_2^4p_1^2p_3^2\\
&&-\frac{1}{12}p_2^4p_1^2p_4^2-\frac{1}{12}p_3^4p_1^2p_2^2+\frac{1}{2}p_1^2p_2^2p_3^2p_4^2-\frac{1}{12}p_4^4p_1^2p_2^2-\frac{1}{12}p_3^4p_1^2p_4^2-\frac{1}{12}p_4^4p_1^2p_3^2\\
&&+\frac{1}{12}p_2^4p_3^4-\frac{1}{12}p_2^4p_3^2p_4^2+\frac{1}{12}p_2^4p_4^4-\frac{1}{12}p_3^4p_2^2p_4^2-\frac{1}{12}p_4^4p_2^2p_3^2+\frac{1}{12}p_3^4p_4^4,\\
U_4&=& (\frac{1}{72}p_1^2p_2^2+\frac{1}{72}p_1^2p_3^2-\frac{1}{36}p_1^2p_4^2-\frac{1}{36}p_2^2p_3^2+\frac{1}{72}p_2^2p_4^2+\frac{1}{72}p_3^2p_4^2)
\times\\
&&(p_1^2p_2^2-2p_1^2p_3^2+p_1^2p_4^2+p_2^2p_3^2-2p_2^2p_4^2+p_3^2p_4^2)\times\\
&&(2p_1^2p_2^2-p_1^2p_3^2-p_1^2p_4^2-p_2^2p_3^2-p_2^2p_4^2+2p_3^2p_4^2)
\end{eqnarray*}
The invariants $U_1,\dots,U_4$ have degrees $2,6,8,12$, respectively. In the coordinates $u$ the
Hessians of the components of the following vector potential give the structure constants of the
product.
\begin{eqnarray*}
F^1&=&-\frac{1188}{35}u_1^7c_1^3-\frac{297}{140}u_1^7c_1^2-\frac{37}{840}u_1^7c_1-\frac{7}{23040}u_1^7+\frac{3}{64}u_2u_1^4+\frac{9}{2}u_2u_1^4c_1\\&&+108u_2u_1^4c_1^2-\frac{1}{8}u_3u_1^3-6u_3u_1^3c_1-54u_2^2u_1c_1-\frac{9}{8}u_2^2u_1+u_1u_4+3u_2u_3,\\
F^2&=&\frac{1}{663552}u_1^9+\frac{1}{6912}u_1^9c_1+\frac{1}{288}u_1^9c_1^2-\frac{81}{20}u_2u_1^6c_1^2-\frac{11}{120}u_2u_1^6c_1-\frac{17}{23040}u_2u_1^6\\
&&-\frac{324}{5}u_2u_1^6c_1^3+\frac{18}{5}u_3u_1^5c_1^2+\frac{3}{20}u_3u_1^5c_1+\frac{1}{576}u_3u_1^5+162u_2^2u_1^3c_1^2+\frac{27}{4}u_2^2u_1^3c_1\\&&+\frac{7}{96}u_2^2u_1^3-\frac{3}{8}u_2u_3u_1^2-18u_2u_3u_1^2c_1+\frac{1}{2}u_3^2u_1-18u_2^3c_1-\frac{3}{8}u_2^3+u_4u_2,\\
F^3&=&\frac{1}{995328}u_1^{10}-\frac{1}{24}u_3u_2u_1^3-\frac{9}{2}u_3u_2u_1^3c_1-108u_3u_2u_1^3c_1^2+\frac{1}{2880}u_3u_1^6+\frac{13}{240}u_3u_1^6c_1\\
&&+\frac{27}{10}u_3u_1^6c_1^2+\frac{216}{5}u_3u_1^6c_1^3-\frac{1}{8}u_2^3u_1-\frac{27}{2}u_2^3u_1c_1-324u_2^3u_1c_1^2+\frac{1}{64}u_2^2u_1^4\\
&&+\frac{39}{16}u_2^2u_1^4c_1+\frac{243}{2}u_2^2u_1^4c_1^2+1944u_2^2u_1^4c_1^3-\frac{23}{96768}u_2u_1^7-\frac{43}{896}u_2u_1^7c_1-\frac{57}{16}u_2u_1^7c_1^2\\
&&-\frac{810}{7}u_2u_1^7c_1^3-\frac{9720}{7}u_2u_1^7c_1^4+\frac{101}{414720}u_1^{10}c_1+\frac{3}{128}u_1^{10}c_1^2+\frac{9}{8}u_1^{10}c_1^3+27u_1^{10}c_1^4\\
&&+\frac{1296}{5}u_1^{10}c_1^5+u_4u_3,\\
F^4&=&5832u_1^3u_2^3c_1^3+\frac{729}{2}u_1^3u_2^3c_1^2+\frac{117}{16}u_1^3u_2^3c_1-\frac{61236}{5}u_2^2u_1^6c_1^4-\frac{5103}{5}u_2^2u_1^6c_1^3\\
&&-\frac{5031}{160}u_2^2u_1^6c_1^2-\frac{543}{1280}u_2^2u_1^6c_1-27u_3^2u_1^4c_1^2-\frac{9}{8}u_3^2u_1^4c_1+\frac{1}{2}u_4^2-486u_3u_2^2u_1^2c_1^2\\
&&-\frac{81}{4}u_3u_2^2u_1^2c_1+\frac{5832}{5}u_3u_2u_1^5c_1^3+\frac{729}{10}u_3u_2u_1^5c_1^2+\frac{117}{80}u_3u_2u_1^5c_1-\frac{1}{15360}u_3u_1^8\\
&&-\frac{2569}{26276659200}u_1^{12}+\frac{43}{1658880}u_2u_1^9-\frac{7641}{35200}u_1^{12}c_1^3-\frac{34421}{10137600}u_1^{12}c_1^2\\
&&-\frac{6881}{243302400}u_1^{12}c_1-\frac{3}{16}u_3u_2^2u_1^2+\frac{3}{320}u_3u_2u_1^5-\frac{1944}{5}u_3u_1^8c_1^4-\frac{162}{5}u_3u_1^8c_1^3\\
&&-\frac{159}{160}u_3u_1^8c_1^2-\frac{17}{1280}u_3u_1^8c_1+\frac{34992}{5}u_2u_1^9c_1^5+729u_2u_1^9c_1^4+\frac{1209}{40}u_2u_1^9c_1^3\\
&&+\frac{399}{640}u_2u_1^9c_1^2+\frac{883}{138240}u_2u_1^9c_1-\frac{4293}{550}c_1^4u_1^{12}-\frac{82377}{550}c_1^5u_1^{12}-\frac{329508}{275}c_1^6u_1^{12}\\&&-\frac{1}{96}u_3^2u_1^4-243u_2^4c_1^2-\frac{81}{8}u_2^4c_1
+\frac{3}{64}u_1^3u_2^3-\frac{97}{46080}u_2^2u_1^6-\frac{3}{32}u_2^4+\frac{1}{2}u_3^3.
\end{eqnarray*}
We now compare the product and connection obtained in this way with the dual structure defined by
the reflecting hyperplanes.
\begin{proposition}\label{prop:F4}
For each $c_1\in\mathbb C$, let $\nabla$ be the connection with flat coordinates $u_1,\dots,u_4$
and let $\circ$ be the product with structure constants
$c^i_{jk}=\partial_{u_j}\partial_{u_k}F^i$.
\begin{enumerate}
\item[(i)] The triple $(\nabla,\circ,e)$ defines a flat F-manifold on the orbit space of $F_4$ with unit
$e=\partial_{u_4}$ and linear Euler vector field
\[
E=\frac1{12}\big(2u_1\partial_{u_1}+6u_2\partial_{u_2}+8u_3\partial_{u_3}+12u_4\partial_{u_4}\big)=\frac1{12}\sum_kp_k\partial_{p_k},
\]
and for every value of $c_1$ the product coincides with that of the Frobenius manifold $F_4$.
\item[(ii)] On the regular quotient the dual product $X*Y=(E\circ)^{-1}X\circ Y$ and the connection defined
by \eqref{eq:dual} are given by
{
\[
X*Y=12\sum_{H}\frac16\,\frac{\alpha_H(X)}{\alpha_H(p)}\,\pi_H(Y),\qquad
\nabla^*_XY=\partial_XY-\sum_{H}\tau_H\,\frac{\alpha_H(X)}{\alpha_H(p)}\,\pi_H(Y),
\]
}
with
{
\[
\tau_H=\frac16-s\ \ (H\in\mathcal O_1),\qquad \tau_H=\frac16+s\ \ (H\in\mathcal O_2),\qquad
s=72\,c_1+\frac32 .
\]
}
The identity $\sum_H\tau_H\pi_H=\Id$ implies $\nabla^*E=0$. We thus obtain the formulas of
Theorem~\ref{ALcomplexthm2} with $\sigma_H=\frac16$ after the normalization $E\to12E$,
$*\to\frac1{12}*$.
\item[(iii)] The pair $(\nabla,\circ,e,\nabla^*,*,E)$ is bi-flat on the regular quotient. If we fix the dual
product together with $e$ and $E$, varying $c_1$ gives all weights compatible with invariance
under $G$ and with the normalization. Thus we obtain the family of Conjecture~\ref{conj:AL}. When
$c_1=-\frac1{48}$ we have $s=0$, so the two weights coincide. At this value we recover the
Frobenius manifold $F_4$ together with its normalized dual structure.
\end{enumerate}
\end{proposition}

\begin{proof}
We prove that the expressions in (ii) define the dual of the product and the connection in (i).
First observe that $\sum_{H\in\mathcal O_i}\pi_H$ is invariant under $G$. By irreducibility this
sum is a scalar multiple of the identity. Its trace is $12$, so the sum equals $3\,\Id$ on either
orbit. Thus $\sum_H\frac16\pi_H=\Id$ and
$\sum_H\tau_H\pi_H=\left[3\big(\frac16-s\big)+3\big(\frac16+s\big)\right]\Id=\Id$. It follows
from Theorem~\ref{ALcomplexthm2} that $(\nabla^*,*,E)$ is a flat F-manifold with $\nabla^*E=0$.

To compare the two products, we evaluate the product defined by the reflecting hyperplanes on the
linear coordinate fields and express the result in the frame associated to $u$. Set
$J_i=(\partial_i u_\alpha)_\alpha$. Denote by $C$ the product in this frame defined by the
Hessians of $F^\alpha$ and by $L$ the operator of multiplication by $E$ with respect to $C$. For
each form $\alpha_H$, let $\alpha_H^\sharp$ be the normal vector with the same Euclidean
coefficients. Define
\[
S^\alpha_{H,ij}
=\frac{(\alpha_H)_i(\alpha_H)_j}{\|\alpha_H\|^2}
  \frac{\partial_{\alpha_H^\sharp}u_\alpha}{\alpha_H(p)},
\qquad T^\alpha_{ij}=2\sum_H S^\alpha_{H,ij}.
\]
Since each $u_\alpha$ is invariant under reflection in $H$, its normal derivative vanishes on
$H$. Thus the division in $S^\alpha_{H,ij}$ is exact. Substituting the polynomials above and
expanding, we obtain
\begin{align*}
\partial_i\partial_j u_\alpha
  &=\sum_H\left(\frac{d_\alpha}6-\tau_H\right) S^\alpha_{H,ij},\\
C(J_i,J_j)&= L T_{ij}.
\end{align*}
In the second identity we evaluate the Hessians at $u=u(p)$. Both identities are polynomial in
$p_1,\ldots,p_4,c_1$. In fact, $T_{ij}$ gives the components of $\partial_{p_i}*\partial_{p_j}$
in the frame associated to $u$, so the second identity yields $X\circ Y=L(X*Y)$ on the regular
quotient. The terms containing $u_4$ in the potentials give the unit relation $X\circ e=X$.
Substituting $Y=e$, we obtain $L(e*)=\Id$. It follows that $E\circ$ is invertible everywhere on
the regular quotient and that $\circ=(e*)^{-1}*$.

The first identity gives the Christoffel symbols of $\nabla^*$ in the coordinates $u$, namely
$\Gamma^{*\alpha}_{\beta\gamma}=-\frac{d_\alpha}{12}c^{*\alpha}_{\beta\gamma}$ (see
Remark~\ref{rem:eigen}). Since $\nabla E=\mathrm{diag}\big(\frac{d_\alpha}{12}\big)$, this is
precisely the relation \eqref{eq:dual}. On the regular quotient the product $\circ=(e*)^{-1}*$ is
associative, because $*$ is associative and $e$ is invertible with respect to $*$. The structure
constants are polynomial, so associativity extends to the whole quotient. The vector potential
gives compatibility with $\nabla$. From the grading we obtain $\mathrm{Lie}_E\circ=\circ$, with
$E$ linear in $u$. Therefore $(\nabla,\circ,e)$ is a flat F-manifold with linear Euler field and
forms a bi-flat pair with $(\nabla^*,*,E)$ by Lemma~\ref{lem:AHE}.

Finally, the triangular change of coordinates gives $e=\partial_{U_4}$ for every $c_1$. Since $*$
is also independent of the parameter, the reconstruction formula gives the same product $\circ$
throughout the family. At $c_1=-1/48$ the dual product and the normalized dual connection agree
with those of Dubrovin's Coxeter construction. If we normalize the Frobenius unit to be
$\partial_{U_4}$, reconstruction therefore identifies the bi-flat pair at this value with the
Frobenius pair. The weights of the two orbits satisfy
$3\tau_{\mathcal O_1}+3\tau_{\mathcal O_2}=1$, and varying $s$ gives every pair satisfying this
relation.
\end{proof}
}

\section{The case of a single orbit}\label{sec:ALconjectures}

{
In this section we consider groups with a single orbit of reflecting hyperplanes. In this case
$\mu=0$, and the normalization fixes the weight of the dual connection. We show that the dual of
the KMS flat structure of Theorem~\ref{ALcomplexthm1} has the form prescribed by
Conjecture~\ref{conj:AL}. In fact, the formula for the dual in
Theorem~\ref{thm:naturalDunklKohno} applies to every nontrivial irreducible well generated
group. For $G(m,m,n)$ we give a second construction in Section~\ref{sec:Gmmn}, starting from a
Laurent superpotential. In this way we identify the member of the Laurent family for which the
natural structure extends to the whole quotient and the dual structure extends to the entire
regular quotient.
}

\subsection{\texorpdfstring{The natural dual connection}{The natural dual connection}}\label{sec:KMSdual}

{
Our first aim is to describe the dual of the KMS flat structure of Remark~\ref{rem:KMSKMS} in
terms of the reflecting hyperplanes of any nontrivial irreducible well generated group. The
cases that remain to be proved are $G_{29},G_{32},G_{33},G_{34}$ and $G(m,m,n)$ with $m\ge3$ and
$n\ge3$. Consider the normalized dual product $B$ in the linear coordinates $p$ on $V$. Its
coefficients are rational functions, homogeneous of degree $-1$, whose poles lie on the
reflecting hyperplanes. At a generic point of such a hyperplane $H$, the inverse Jacobian gives
the bound \eqref{eq:polesInitialBounds} on the orders of the poles. We then use the invariance
of $B$ under the cyclic stabilizer of $H$ and compatibility with the trivial connection to show
that the polar part of $B$ is a constant multiple of the elementary tensor $A_H$. Once we
subtract these contributions from all the hyperplanes, the remaining tensor is holomorphic
outside a set of codimension at least two. Its coefficients therefore extend holomorphically to
$V$, where their homogeneity of degree $-1$ implies that they vanish. 

\begin{theorem}\label{thm:naturalDunklKohno}
Let $G$ be a nontrivial finite well generated irreducible complex reflection group acting on
$V=\mathbb C^n$ by its standard reflection representation, and let $h$ be its largest invariant
degree. The dual of the KMS flat structure described in Theorem~7.5 and Section~9.2 of
\cite{KMS} is defined on $M^{\mathrm{reg}}$. We pull the dual product back to
{
\[
V^{\mathrm{reg}}=V\setminus\bigcup_{H\in\mathcal H}H,
\]
}
and set
{
\[
\mathcal E=\sum_i p_i\partial_{p_i},\qquad
E=\frac1h\mathcal E,\qquad
c^*_{\mathrm{nat}}(X,Y)
=(E\circ_{\mathrm{nat}})^{-1}(X\circ_{\mathrm{nat}}Y),\qquad
B=\frac1h c^*_{\mathrm{nat}}.
\]
}
Fix a positive definite Hermitian form invariant under $G$. For each $H\in\mathcal H$, choose a
nonzero linear form $\alpha_H$ with kernel $H$, and let $\pi_H$ denote the orthogonal projection
onto $H^\perp$. Then there exist constants $\rho_H\in\mathbb C$ such that
{
\[
\rho_{gH}=\rho_H\quad(g\in G),\qquad
\sum_{H\in\mathcal H}\rho_H\pi_H=\Id,
\]
}
and such that the normalized dual product and the dual connection are given by
{
\begin{equation}\label{eq:naturalDunklKohno}
\begin{aligned}
B_p(X,Y)&=\sum_{H\in\mathcal H}\rho_H
\frac{\alpha_H(X)}{\alpha_H(p)}\pi_H(Y),\\
(\nabla^*_{\mathrm{nat}})_XY
&=\partial_XY-B_p(X,Y),\qquad \nabla^*_{\mathrm{nat}}E=0,
\end{aligned}
\end{equation}
}
where $\partial$ denotes the trivial connection on $V$. The expression for $B$ as a sum over the
reflecting hyperplanes determines the constants $\rho_H$ uniquely.

If the reflecting hyperplanes of $G$ form a single orbit, then $\rho_H=n/|\mathcal H|$ for every
$H$ and
{
\begin{equation}\label{eq:keyfor4}
B_p(X,Y)=\frac{n}{|\mathcal H|}
\sum_{H\in\mathcal H}\frac{\alpha_H(X)}{\alpha_H(p)}\pi_H(Y).
\end{equation}
}
In particular, the natural structures for $G_{29},G_{32},G_{33},G_{34}$ have the form prescribed
by Conjecture~\ref{conj:AL} with $\mu=0$.
\end{theorem}

\begin{proof}
Choose homogeneous basic invariants $u^1,\ldots,u^n$ that form a system of flat coordinates for
the KMS flat structure. Denote their degrees by $d_1,\ldots,d_n$ and set
{
\[
J^\alpha_i=\frac{\partial u^\alpha}{\partial p_i}.
\]
}
By Theorem~7.5(4) and Section~6.2 of \cite{KMS}, the connection of the natural almost Saito
structure coincides with $\partial$ in the linear coordinates $p$. Since the authors list the
degrees in decreasing order, their parameter $1/d_1$ corresponds to $1/h$ in our notation. Set
$Q=\nabla E$. The covariant formula of Remark~\ref{rem:KMS} gives
$\partial_XY=\nabla_XY+(h^{-1}\Id-Q)c^*_{\mathrm{nat}}(X,Y)$. Using the fact that $du^\alpha$ is
parallel for $\nabla$ and that $du^\alpha\circ Q=(d_\alpha/h)du^\alpha$, we obtain
\[
(\partial^2u^\alpha)(X,Y)=(d_\alpha-1)du^\alpha(B(X,Y)).
\]
Since the reflection representation is nontrivial and irreducible, we have $d_\alpha>1$ for
every index. We can therefore divide by $d_\alpha-1$ and invert $J$ on $V^{\mathrm{reg}}$. This
gives the following expression for the components of $B$, as in \cite{ALrefl} and Section~7.4 of
\cite{KMS}.
{
\begin{equation}\label{eq:ALdualProduct}
B^k_{ij}
=
\sum_{\alpha=1}^n
\frac{1}{d_\alpha-1}
(J^{-1})^k_\alpha
\frac{\partial^2u^\alpha}{\partial p_i\partial p_j}.
\end{equation}
}
It follows from \eqref{eq:dual} that the normalized dual connection is $\partial-B$, since
$B=h^{-1}c^*_{\mathrm{nat}}$, and that $\nabla^*_{\mathrm{nat}}E=0$. In what follows we use the
compatibility of the natural almost Saito structure with the trivial connection in the form
{
\begin{equation}\label{eq:KMS2}
\partial_\ell B^k_{ij}
=
\partial_i B^k_{\ell j}.
\end{equation}
}
Indeed, the connection $\partial+\nu\,c^*_{\mathrm{nat}}$ is flat for every $\nu$ by
Remark~\ref{rem:KMS} and Proposition~3.1 of \cite{KMS}. The coefficient of $\nu$ in its
curvature is the antisymmetric part of $\partial c^*_{\mathrm{nat}}$ and must therefore vanish.

Since $u(gp)=u(p)$ for $g\in G$, differentiation gives $J(gp)g=J(p)$ and
$(\partial^2u)_{gp}(gX,gY)=(\partial^2u)_p(X,Y)$. Substituting these identities into
\eqref{eq:ALdualProduct} gives
\begin{equation}\label{eq:Binv}
B_{gp}(gX,gY)=g\,B_p(X,Y)\qquad(g\in G).
\end{equation}
The entries in row $\alpha$ of $J$ have degree $d_\alpha-1$, whereas the second derivatives of
$u^\alpha$ have degree $d_\alpha-2$. Thus the coefficients of $B$ are rational functions
homogeneous of degree $-1$. To see that they are holomorphic on $V^{\mathrm{reg}}$, recall that
formulas (6.3) and (6.4) of \cite{KMS} give $\det J=c\prod_{H\in\mathcal H}\alpha_H^{\ell_H-1}$,
where $c\neq0$ and $\ell_H=|G_H|$.

\emph{Local form along a hyperplane.} Fix $H\in\mathcal H$ and set $\ell=\ell_H\ge2$. Choose
linear coordinates $v=(t,y^2,\dots,y^n)$, with $t$ as coordinate on $H^\perp$ and $y$ as
coordinates on $H$. Thus $H=\{t=0\}$. The pointwise stabilizer $G_H$ fixes $H$ and acts
faithfully on the normal line. It is therefore cyclic, with generator $g(t,y)=(\zeta t,y)$,
where $\zeta$ is a primitive $\ell$th root of unity. The change from $p$ to $v$ is linear, so
\eqref{eq:ALdualProduct}, \eqref{eq:KMS2} and \eqref{eq:Binv} hold in the coordinates $v$ with
$J=\partial u/\partial v$. We denote by $B^k_{ij}$ the components of $B$ in these coordinates.
The index $1$ corresponds to $t$ and the indices $a,b\ge2$ correspond to $y$.

By the invariance $u(\zeta t,y)=u(t,y)$, the only powers of the normal coordinate that occur are
$t^{k\ell}$. We can thus express $u$ in the form $u(t,y)=F(t^\ell,y)$ for a polynomial map $F$.
The chain rule gives
\[
J(t,y)=K(t^\ell,y)\,\operatorname{diag}\bigl(\ell t^{\ell-1},1,\dots,1\bigr),\qquad K(s,y):=\frac{\partial F}{\partial(s,y)} .
\]
Denote by $H^\circ$ the set of points of $H$ that lie on no other reflecting hyperplane, and fix
$q=(0,y_0)\in H^\circ$. Choose a neighborhood $U=\{|t|<\varepsilon,\ |y-y_0|<\delta\}$ of $q$
that is disjoint from those hyperplanes, so that $U\setminus H\subset V^{\mathrm{reg}}$. Near
$q$ $\det J$ vanishes on $H$ to order exactly $\ell-1$. It follows that $\det K(0,y_0)\neq0$ and
that $J^{-1}=\operatorname{diag}\bigl((\ell t^{\ell-1})^{-1},1,\dots,1\bigr)K(t^\ell,y)^{-1}$
near $q$ whenever $t\neq0$. In particular, the rows of $J^{-1}$ with index $a\ge2$ are
holomorphic near $q$, while the entries of the first row are $t^{1-\ell}$ times holomorphic
functions. Since the second derivatives of $u$ are polynomial, we obtain from
\eqref{eq:ALdualProduct} the bounds
\begin{equation}\label{eq:polesInitialBounds}
B^a_{ij}\ \text{holomorphic near }q\quad(a\ge2),\qquad t^{\ell-1}B^1_{ij}\ \text{holomorphic near }q .
\end{equation}

To refine these bounds, we use the action of the stabilizer. The matrix of $g$ is
$\operatorname{diag}(\zeta,1,\dots,1)$. In components, \eqref{eq:Binv} therefore reads
$B^k_{ij}(\zeta t,y)\,\zeta^{\varepsilon_i+\varepsilon_j}=\zeta^{\varepsilon_k}B^k_{ij}(t,y)$,
where $\varepsilon_1=1$ and $\varepsilon_a=0$ for $a\ge2$. In particular,
\[
B^1_{11}(\zeta t,y)=\zeta^{-1}B^1_{11}(t,y),\qquad B^1_{1a}(\zeta t,y)=B^1_{1a}(t,y),\qquad B^1_{ab}(\zeta t,y)=\zeta\,B^1_{ab}(t,y).
\]
It follows from \eqref{eq:polesInitialBounds} that each component $B^1_{ij}$ has a Laurent
expansion $\sum_{k\ge1-\ell}f_k(y)t^k$ whose coefficients $f_k$ are holomorphic near $y_0$ and
which converges on $U\setminus H$. If a function satisfies $f(\zeta t,y)=\zeta^mf(t,y)$, the
only exponents that occur are those for which $k\equiv m\pmod\ell$. Since $\ell\ge2$, the only
possible negative exponent is $-1$ for $B^1_{11}$ and $1-\ell$ for $B^1_{ab}$. No negative
exponent occurs for $B^1_{1a}$. We thus obtain
\[
B^1_{11}=\frac{r(y)}t+\text{hol.},\qquad B^1_{1a}=\text{hol.},\qquad B^1_{ab}=c_{ab}(y)\,t^{1-\ell}+\text{hol.}
\]
We now apply the compatibility condition \eqref{eq:KMS2}. In the identity
$\partial_tB^1_{ab}=\partial_{y^a}B^1_{1b}$, the right hand side is holomorphic. The coefficient
of $t^{-\ell}$ on the left is $(1-\ell)c_{ab}(y)$, so $c_{ab}=0$. Similarly, the right hand side
of $\partial_{y^a}B^1_{11}=\partial_tB^1_{a1}$ is holomorphic, and the term
$\partial_{y^a}r(y)\,t^{-1}$ on the left must vanish. It follows that $r$ is locally constant.
Moreover, $r(y)=\lim_{t\to0}tB^1_{11}(t,y)$ is defined on all of $H^\circ$. Since its complement
in $H$ is a finite union of proper complex subspaces, this domain is connected. Thus $r$ takes a
single value $r_H$ on $H^\circ$.

For $p\in V^{\mathrm{reg}}$, define
\[
A_{H,p}(X,Y):=\frac{\alpha_H(X)}{\alpha_H(p)}\pi_H(Y).
\]
The tensor is independent of the choice of scale for $\alpha_H$. Taking $\alpha_H=t$, we have
$\pi_H(Y)=dt(Y)\,\partial_t$, so the only nonzero component of $A_H$ in the coordinates $v$ is
$(A_H)^1_{11}=1/t$. We have thus proved that
\begin{equation}\label{eq:polesLocalTensor}
B-r_HA_H\ \text{extends holomorphically to a neighborhood of every point of }H^\circ .
\end{equation}

\emph{Global form.} We now subtract the polar parts at all the reflecting hyperplanes. Define
\begin{equation}\label{eq:defR}
R:=B-\sum_{H\in\mathcal H}r_HA_H .
\end{equation}
By construction, $R$ is holomorphic on $V^{\mathrm{reg}}$. Near any point of $H^\circ$, the
difference $B-r_HA_H$ is holomorphic by \eqref{eq:polesLocalTensor}, and each term $A_{H'}$ with
$H'\neq H$ is holomorphic as well. Thus $R$ is holomorphic on the complement in $V$ of the union
of the intersections $H\cap H'$ with $H\neq H'$. This union has codimension two, so Riemann's
extension theorem gives a holomorphic extension of the coefficients of $R$ to $V$. These
coefficients are homogeneous of degree $-1$, as are those of $B$ and $A_H$. In particular,
$R^k_{ij}(cp)=c^{-1}R^k_{ij}(p)$ for $c\in\mathbb C^*$. Any nonzero coefficient would therefore
be unbounded near $0$, which contradicts holomorphicity. Thus $R=0$ and
\begin{equation}\label{eq:polesGlobalProduct}
B_p(X,Y)=\sum_{H\in\mathcal H}r_H\frac{\alpha_H(X)}{\alpha_H(p)}\pi_H(Y)
\end{equation}
on the complement of the reflecting hyperplanes. We set $\rho_H=r_H$.

\emph{Uniqueness, invariance and normalization.} Suppose that two families $\rho_H$ and
$\rho'_H$ give the same sum. Then $(\rho_{H_0}-\rho'_{H_0})A_{H_0}$ is holomorphic near a point
of $H_0^\circ$, which implies $\rho_{H_0}=\rho'_{H_0}$. The invariance of the Hermitian form
gives $\pi_{gH}g=g\pi_H$. Moreover, $\alpha_{gH}\circ g$ is a multiple of $\alpha_H$, so
$A_{gH,gp}(gX,gY)=gA_{H,p}(X,Y)$. By \eqref{eq:Binv}, applying $g$ to
\eqref{eq:polesGlobalProduct} gives a second expression for $B$ in which the weights are
permuted. By uniqueness,
\begin{equation}\label{eq:polesOrbitWeights}
\rho_{gH}=\rho_H .
\end{equation}
The fact that $E$ is the unit of $c^*_{\mathrm{nat}}$ and the identity $\mathcal E=hE$ give
$B(\mathcal E,X)=X$. Evaluating \eqref{eq:polesGlobalProduct} on $(\mathcal E,X)$ and using
$\alpha_H(\mathcal E)=\alpha_H(p)$, we find
\begin{equation}\label{eq:polesNormalization}
\sum_{H\in\mathcal H}\rho_H\pi_H=\Id .
\end{equation}
Suppose now that $G$ has a single orbit of reflecting hyperplanes. By
\eqref{eq:polesOrbitWeights} we have $\rho_H=\rho$ for all $H$. Taking the trace in
\eqref{eq:polesNormalization} gives $\rho|\mathcal H|=n$, which proves \eqref{eq:keyfor4}.

For the groups $G_{29},G_{32},G_{33},G_{34}$ the reflecting hyperplanes form a single
orbit~\cite{AW}. This proves the last statement. In Corollary~\ref{cor:oneorbit} we describe the
consequences for Conjecture~\ref{conj:AL}.

\end{proof}

Combining this expression for the dual product with the flat pencils of
Theorem~\ref{ALcomplexthm2}, we obtain the following result.
\begin{corollary}\label{cor:oneorbit}
Let $G$ be a nontrivial well generated irreducible complex reflection group and let
$(\nabla,\circ,e,E)$ be the KMS flat structure on its orbit space given by
Theorem~\ref{ALcomplexthm1}, with $[e,E]=e$. Set $X*Y=(E\circ)^{-1}X\circ Y$ and define
$\nabla^*$ by \eqref{eq:dual}. Then $(\nabla,\circ,e,\nabla^*,*,E)$ is a bi-flat F-manifold on
$M^{\mathrm{reg}}$. In the linear coordinates on $V^{\mathrm{reg}}$ the dual product and
connection are given by
{
\[
X*Y=h\sum_{H\in\mathcal H}\rho_H\frac{\alpha_H(X)}{\alpha_H(p)}\pi_H(Y),\qquad
\nabla^*_XY=\partial_XY-\sum_{H\in\mathcal H}\rho_H\frac{\alpha_H(X)}{\alpha_H(p)}\pi_H(Y),
\]
}
where the weights are invariant under $G$ and satisfy $\sum_H\rho_H\pi_H=\Id$. In particular,
Conjecture~\ref{conj:AL} holds when $G$ has one orbit of reflecting hyperplanes, in which case
$\rho_H=n/|\mathcal H|$. Once $e$ and $E$ are fixed, the bi-flat pair in this case is unique
among those for which the dual product and connection have the form of
Theorem~\ref{ALcomplexthm2}, with the product rescaled as above.
\end{corollary}

\begin{proof}
The Euler vector field of the KMS flat structure is linear, and $E\circ$ is invertible on
$M^{\mathrm{reg}}$ by Theorem~\ref{ALcomplexthm1}. Proposition~\ref{prop:dualflat} then shows
that $(\nabla^*,*,E)$ is a flat F-manifold and that the two structures define a bi-flat F-manifold on $M^{\mathrm{reg}}$. The linear coordinates $p^i$ are flat for
$\Gamma(0,\frac1h)=\nabla^*+\frac1h\,*$ by Remark~\ref{rem:KMSKMS}. Thus
$\nabla^*=\partial-\frac1h\,*=\nabla^*_{\mathrm{nat}}$, and Theorem~\ref{thm:naturalDunklKohno}
gives the stated formulas with $*=c^*_{\mathrm{nat}}=hB$. This is the structure of
Theorem~\ref{ALcomplexthm2} after multiplying its product by $h$ and dividing its unit by $h$.

If the hyperplanes form one orbit, invariance and normalization determine the common weight
$\rho_H=n/|\mathcal H|$. The same normalization fixes the weight of the rescaled dual product.
Once the unit is fixed, the reconstruction following Lemma~\ref{lem:AHE} determines the natural
product and connection from the dual data.
\end{proof}

To compare this formula with the one obtained from a Laurent superpotential below, we express
the weights in terms of the corresponding normal vectors.
\begin{remark}\label{rem:consistency}
For $G(m,m,n)$ with $m\ge2$ and $n\ge3$, we have $|\mathcal H|=m\binom n2$ and therefore
$\rho_H=\frac{2}{m(n-1)}=\frac2h$. The vectors used in Theorem~\ref{thm:Gmmn} are
$v_H=\partial_i-\zeta^{-1}\partial_j$. Since $\alpha_H(v_H)=2$, the orthogonal projection is
$\pi_H=\frac12\alpha_H\otimes v_H$, which accounts for the weight $\frac1h$ in that theorem. For
an irreducible Coxeter group every reflection has order two. In this case
$|\mathcal H|=\frac{nh}2$ gives $\rho_H=\frac2h$, in agreement with Dubrovin's almost duality.
\end{remark}

We now give the construction for $G(m,m,n)$ starting from a Laurent superpotential.
}

\subsection{\texorpdfstring{The groups $G(m,m,n)$}{The groups G(m,m,n)}}\label{sec:Gmmn}

{
Let $m\ge2$, $n\ge3$ and $N=n-1$. The hyperplanes $p_i=\zeta p_j$ with $\zeta^m=1$ form a single
orbit under $G(m,m,n)$. For $m=2$ the group is the Coxeter group $D_n$. As there is only one
orbit, the normalization determines the weights in the expressions for the dual product and
connection. Once we fix the unit $e=\partial_{a_N}$ and the normalized Euler field, the dual
data determine the bi-flat pair on the regular quotient. We now construct a distinguished
superpotential for this pair.

Set $z_j=p_j^m$ and $t=p_1\cdots p_n$. As basic invariants we take the signed elementary
symmetric functions $a_1,\dots,a_{n-1}$ of the $z_j$, together with $t$. The respective degrees
are $m,2m,\dots,(n-1)m$ and $n$, and the Coxeter number is $h=m(n-1)$. In terms of these
invariants we describe the quotient by the family of Laurent polynomials
\begin{equation}\label{eq:Gmmn}
\lambda(z)=\frac{\prod_j(z-z_j)}{z}=z^{N}+a_1z^{N-1}+\dots+a_{N}+\frac{a_n}z,\qquad a_n=(-1)^nt^m .
\end{equation}
As for $B_n$, we reconstruct the bi-flat structure from the superpotential
(Theorem~\ref{thm:Gmmn}). Here we begin on the space $M'$ of coefficients for which $a_n\ne0$.
For every $\beta\ne0$ the construction gives a natural structure on $M'$, as we show in
Corollary~\ref{thm:Laurent}. Using $a_n=(-1)^nt^m$, we pull this structure back to the part of
the quotient of $G(m,m,n)$ where $t\ne0$. It extends to the whole quotient precisely when
$\beta=\frac1m$, as explained in Remark~\ref{rem:Gmmnbeta}. For this value of the parameter the
dual structure also extends to the entire regular quotient.
}

\subsubsection{The data}\label{sec:Laurent}

{
Choose local branches $p_j=z_j^{1/m}$ with product $t$, and fix $\epsilon$ satisfying
$\epsilon^m=(-1)^n$ so that $a_n^{1/m}=\epsilon t$. The space $M'$ of Laurent polynomials
\eqref{eq:Gmmn} satisfying $a_n\ne0$ is the quotient by $S_n$ of the tuples of nonzero roots
$z_j$ of $z\lambda(z)$. The relation $a_n=(-1)^nt^m$ gives a map from the part of the quotient
of $G(m,m,n)$ where $t\ne0$ to $M'$. We assign the degrees $\deg z=\frac1N$ and
$\deg a_k=\frac kN$ for $k=1,\dots,n$.

Fix $\beta\ne0$ and set $\omega=z^{\beta-1}dz=dP$, where $P=z^\beta/\beta$. For $\beta=\frac1m$,
the substitution $z=p^m$ gives $\omega=m\,dp$. We identify the tangent space to $M'$ with
$V=\mathrm{span}(z^{N-1},\dots,z,1,z^{-1})$ by differentiating $\lambda$ in the coordinate
directions. The resulting basis consists of $\partial_{a_i}\lambda=z^{N-i}$ for $i\le N$,
together with $\partial_{a_n}\lambda=z^{-1}$. To construct flat coordinates for the natural
connection, we take residues at the two poles of $\lambda$ and define
{
\begin{equation}\label{eq:TL}
u_k=\frac1{s_k}\res_{z=\infty}\lambda^{s_k}z^{\beta-1}dz\ \ \Big(s_k=\frac{k-\beta}N,\ k=1,\dots,N\Big),\qquad
u_n=\frac1\beta\res_{z=0}\lambda^{\beta}z^{\beta-1}dz=\frac{a_n^\beta}\beta ,
\end{equation}
}
where at infinity we follow the convention after \eqref{eq:u}, with $N$ in place of $n$. If
$s_k=0$, we define $u_k$ as the coefficient of $z^{-k}$ in $\log(\lambda/z^N)$. In particular,
$u_k-a_k$ is polynomial in $a_1,\dots,a_{k-1}$. The powers $z^\beta$ and $a_n^\beta$ are defined
using local branches of the logarithm on $z\ne0$ and $a_n\ne0$, respectively.

}

\subsubsection{The theorem}\label{sec:descent}

{
For each hyperplane $p_i=\zeta p_j$, where $i<j$ and $\zeta^m=1$, define
$\alpha_H=p_i-\zeta p_j$ and $v_H=\partial_i-\zeta^{-1}\partial_j$. These hyperplanes form the
single orbit, and the associated endomorphisms satisfy
$\sum_H\alpha_H\otimes v_H=m(n-1)\Id=h\,\Id$.

We now describe the bi-flat structure that corresponds to $\beta=1/m$.

\begin{theorem}\label{thm:Gmmn}
Let $n\ge3$, $m\ge2$, $\beta=\frac1m$ and $\nu_0=\frac\beta N=\frac1h$. In the coordinate $p$,
consider the Laurent superpotential $\Lambda$ and the primary form $\omega$ given by
\[
\Lambda=\lambda(p^m)=p^{m(n-1)}+a_1p^{m(n-2)}+\dots+a_{n-1}+(-1)^n\frac{t^m}{p^m},\qquad\omega=z^{\beta-1}dz=m\,dp .
\]
These data define a natural flat structure $(\nabla,\circ,e,E)$ which is polynomial on the whole
quotient $M=\mathbb C^n$. On $M^{\mathrm{reg}}$ this structure belongs to a bi-flat pair
$(\nabla,\circ,e,\nabla^*,*,E)$ with the following properties.
\begin{enumerate}
\item[(i)] On the locus $t\ne0$, the product $\circ$ is induced by multiplication in
$\mathbb C[z,z^{-1}]/(z^2\lambda')$ under the identification $\partial_{a_i}\mapsto z^{N-i}$ for
$i\le N$ and $\partial_{a_n}\mapsto z^{-1}$. In the coordinates $(a_1,\dots,a_N,t)$ its
structure constants are polynomial and define the product on all of $M$. The unit is
$e=\partial_{a_N}$, and the Euler field is
\[
E=\sum_{k=1}^N\frac{k}{N}a_k\partial_{a_k}+\frac{n}{h}t\partial_t
=\frac1h\sum_jp_j\partial_{p_j}.
\]
\item[(ii)] In the coordinates $p$, the dual product $X*Y:=E^{-1}\circ X\circ Y$ takes the form
{
\[
X*Y=\sum_{p_i=\zeta p_j}\frac{\alpha_H(X)\,\alpha_H(Y)}{\alpha_H(p)}\,v_H .
\]
}
\item[(iii)] Choose the branches so that on $M^{\mathrm{reg}}\cap\{t\ne0\}$ the twisted zeros of
the superpotential are $w_j=z_j^\beta/\beta=m\,p_j$. The connection $\tilde\nabla$ for which
these functions are flat coordinates coincides with $\partial$ in the linear coordinates $p$ and
extends to the entire regular locus. The formula
{
\[
\nabla^*:=\tilde\nabla-\nu_0\,*=\partial-\frac1h\,*
\]
}
gives the connection of Dunkl and Kohno invariant under $G(m,m,n)$, with weight $\frac2h$ when
the projections $\pi_H$ are used. The connection satisfies $\nabla^*E=0$, and $(\nabla^*,*,E)$
is a flat F-manifold.
\item[(iv)] The functions $u_1,\dots,u_N$ of \eqref{eq:TL}, together with $u_n=\epsilon\,m\,t$,
form a system of polynomial coordinates on $M$ which are flat for $\nabla$. On
$M^{\mathrm{reg}}$ we have $\nabla_XY=\nabla^*_XY-\nabla^*_{X\circ Y}e$. The unit is
$e=\partial_{u_N}$, and $(\nabla,\circ,e)$ is a flat F-manifold with polynomial vector
potential.
\item[(v)] The bi-flat F-manifold $(\nabla,\circ,e,\nabla^*,*,E)$ admits $\Lambda$ as a local
superpotential in the coordinate $P=z^\beta/\beta=m\,p$, with weight $1/h$ and constant
$\kappa=0$, on the inverse charts in $a_n\ne0$. At generic points of the base the $n$ critical
values on the $z$ curve form a system of canonical coordinates. Each nonzero simple critical
point lifts to $m$ points on the $p$ curve, all with the same critical value. For $e$ and $E$ as
in (i), the bi-flat pair on $M^{\mathrm{reg}}$ is unique among pairs whose dual product and
connection have the form of Theorem~\ref{ALcomplexthm2}, with the product rescaled by $h$.
\end{enumerate}
For $m=2$ we recover the Frobenius manifold $D_n$. In this case $\lambda(p^2)$ with
$\omega\propto dp$ is the superpotential of Dubrovin \cite{Dub}, and $\frac1h=\frac{1-d}2$.
\end{theorem}

\begin{proof}
The argument follows the proof of Theorem~\ref{thm:Bn}. We first construct the product and its
dual in Steps~1 and~2. In Step~3 we use \eqref{eq:chrL} to define the natural connection, verify
(S3) and apply Theorem~\ref{thm:flat}. In Step~4 we use Proposition~\ref{prop:dualflat} to
obtain the dual structure. In Step~5 we identify the superpotential and the canonical
coordinates and prove uniqueness. A feature of the Laurent case is that $\lambda$ has a second
pole at $z=0$. This pole contributes to the definition of the connection and supplies one
additional flat coordinate. By using the invariant $t$ we can prove polynomiality in coordinates
that are defined across $t=0$.

\emph{Step 1.} On $a_n\ne0$ the critical points of the superpotential on the $z$ curve are the zeros of
$z^2\lambda'(z)$. This polynomial has degree $n$ and constant term $-a_n$, so all its zeros are
nonzero. Thus $z$ is invertible modulo $z^2\lambda'$ and the map
$V\to\mathbb C[z,z^{-1}]/(z^2\lambda')$ is an isomorphism. The identity
$(X\circ Y)\lambda\equiv X\lambda\,Y\lambda$ modulo $z^2\lambda'$ defines a commutative
associative product with unit $\partial_{a_N}$, since $\partial_{a_N}\lambda=1$. With the
grading $\deg z=\frac1N$, $\deg a_k=\frac kN$, every monomial in \eqref{eq:Gmmn} has degree $1$.
As in Step~1 of the proof of Theorem~\ref{thm:Bn}, we obtain
$E\lambda=\lambda-\frac1Nz\lambda'$. Since $z\lambda'=z^{-1}\cdot z^2\lambda'\equiv0$, the
endomorphism $E\circ$ acts on the quotient algebra by multiplication by $\lambda$. The grading
also gives $\mathrm{Lie}_E\circ=\circ$. In the coordinates $p$ we have
$\deg p=\frac1{mN}=\frac1h$ and therefore $E=\frac1h\sum_jp_j\partial_{p_j}$.

To prove that the product extends polynomially across $t=0$, we express it in the basis
$\partial_{a_1},\dots,\partial_{a_N},\partial_t$. The functions corresponding to these vectors
are $\partial_{a_i}\lambda=z^{N-i}$ and $\partial_t\lambda=(-1)^nm\,t^{m-1}z^{-1}$. From the
identity $z^2\lambda'=\sum_{k=0}^{N-1}(N-k)a_kz^{n-k}-a_n$, where $a_0=1$, we obtain two rules
for reducing products modulo $z^2\lambda'$.
{
\[
z^{N}\equiv\frac1N\Big(a_nz^{-1}-\sum_{k=1}^{N-1}(N-k)a_kz^{N-k}\Big),\qquad
z^{-2}\equiv\frac1{a_n}\sum_{k=0}^{N-1}(N-k)a_kz^{N-1-k} .
\]
}
Consider first a product $z^{N-i}\,z^{N-j}$ of two polynomial basis elements. Under the first
rule, every contribution to the coefficient of $z^{-1}$ has a factor $\frac{a_n}N$. In the basis
containing $\partial_t\lambda$, the corresponding coefficient is therefore a polynomial in
$a_1,\dots,a_N$ multiplied by $a_n/\big((-1)^nm\,t^{m-1}\big)=t/m$. The mixed product
$z^{N-i}\cdot\partial_t\lambda=(-1)^nm\,t^{m-1}z^{N-i-1}$ already belongs to the span of
$z^{N-2},\dots,1,\partial_t\lambda$ with polynomial coefficients. For the square of
$\partial_t\lambda$, the second rule gives
{
\[
(\partial_t\lambda)^2=m^2t^{2m-2}z^{-2}\equiv(-1)^nm^2\,t^{m-2}\sum_{k=0}^{N-1}(N-k)a_kz^{N-1-k},
\]
}
which is polynomial in $t$ since $m\ge2$. Thus all the structure constants of the product are
polynomial in $(a_1,\dots,a_N,t)$.

\emph{Step 2.} We first compute the dual product in the coordinates $z_j$ given by the roots. We
then pass to the coordinates $p_j$ to express the result in terms of the reflecting hyperplanes.
In the quotient algebra the identity $E\circ(X*Y)=X\circ Y$ becomes
$\lambda\cdot(X*Y)\lambda\equiv X\lambda\,Y\lambda\pmod{z^2\lambda'}$. The element $\lambda$ is
invertible wherever $z\lambda(z)$ has simple roots and $t\ne0$. Consequently, whenever
$X\lambda\,Y\lambda/\lambda$ is a Laurent polynomial belonging to $V$, it coincides with
$(X*Y)\lambda$.

On the locus where $z\lambda(z)$ has distinct roots, its nonzero roots $z_1,\dots,z_n$ form a
system of local coordinates on $M'$. Differentiating $\lambda=z^{-1}\prod_k(z-z_k)$ with respect
to $z_i$ gives
{
\[
\partial_{z_i}\lambda=-z^{-1}\prod_{k\neq i}(z-z_k)=-\frac{\lambda}{z-z_i}\in V .
\]
}
For $i\ne j$ the quotient
$\partial_{z_i}\lambda\;\partial_{z_j}\lambda/\lambda=z^{-1}\prod_{k\neq i,j}(z-z_k)$ belongs to
$V$ because the polynomial in its numerator has degree $n-2=N-1$. Applying the partial fraction
identity already used for $B_n$, we find
{
\[
\partial_{z_i}*\partial_{z_j}=-\frac{\partial_{z_i}-\partial_{z_j}}{z_i-z_j}\qquad(i\neq j).
\]
}

Passing to the variables $z_i=p_i^m$ gives $\partial_{p_i}=m\,p_i^{m-1}\partial_{z_i}$. Thus for
$i\neq j$ we have
\[
\begin{gathered}
\partial_{p_i}*\partial_{p_j}=m^2p_i^{m-1}p_j^{m-1}\,\partial_{z_i}*\partial_{z_j} \\
=-\frac{m^2p_i^{m-1}p_j^{m-1}}{p_i^m-p_j^m}\Big(\frac{\partial_{p_i}}{m\,p_i^{m-1}}-\frac{\partial_{p_j}}{m\,p_j^{m-1}}\Big) \\
=-m\,\frac{p_j^{m-1}\partial_{p_i}-p_i^{m-1}\partial_{p_j}}{p_i^m-p_j^m}.
\end{gathered}
\]

To compare this expression with the formula in (ii), let $F(X,Y)$ denote the sum over the
reflecting hyperplanes in that formula, so that
{
\[
F(X,Y)=\sum_H\frac{\alpha_H(X)\,\alpha_H(Y)}{\alpha_H(p)}\,v_H ,
\]
}
By symmetry of $F$, we need only compare $F(\partial_{p_i},\partial_{p_j})$ with the preceding expression for $i<j$.
The form associated with $p_k=\zeta p_l$ depends only on $p_k$ and $p_l$. Thus a hyperplane can
contribute only when $\{k,l\}=\{i,j\}$. This leaves the $m$ hyperplanes $p_i=\zeta p_j$, for each of which we have
{
\[
\alpha_H(\partial_{p_i})=1,\qquad \alpha_H(\partial_{p_j})=-\zeta,\qquad \alpha_H(p)=p_i-\zeta p_j,\qquad
v_H=\partial_{p_i}-\zeta^{-1}\partial_{p_j},
\]
}
so its contribution is
{
\[
\frac{(1)(-\zeta)}{p_i-\zeta p_j}\big(\partial_{p_i}-\zeta^{-1}\partial_{p_j}\big)
=\frac{-\zeta\,\partial_{p_i}+\partial_{p_j}}{p_i-\zeta p_j}.
\]
}
Summing over the $m$th roots of unity gives
{
\begin{equation}\label{eq:Fij}
F(\partial_{p_i},\partial_{p_j})=-\partial_{p_i}\sum_{\zeta^m=1}\frac{\zeta}{p_i-\zeta p_j}+\partial_{p_j}\sum_{\zeta^m=1}\frac1{p_i-\zeta p_j}.
\end{equation}
}

To evaluate the two sums, we factor $\xi^m-\eta^m$ and choose compatible local
branches of the logarithm. We have
{
\[
\xi^m-\eta^m=\prod_{\zeta^m=1}(\xi-\zeta\eta),\qquad\text{so}\qquad
\log(\xi^m-\eta^m)=\sum_{\zeta^m=1}\log(\xi-\zeta\eta)
\]
}
Taking the derivatives with respect to $\xi$ and $\eta$, respectively, we obtain
{
\[
\frac{m\,\xi^{m-1}}{\xi^m-\eta^m}=\sum_\zeta\frac1{\xi-\zeta\eta},\qquad
\frac{-m\,\eta^{m-1}}{\xi^m-\eta^m}=\sum_\zeta\frac{-\zeta}{\xi-\zeta\eta}.
\]
}

Substituting $\xi=p_i$ and $\eta=p_j$ in these identities, we obtain from \eqref{eq:Fij}
{
\[
F(\partial_{p_i},\partial_{p_j})=-\partial_{p_i}\,\frac{m\,p_j^{m-1}}{p_i^m-p_j^m}+\partial_{p_j}\,\frac{m\,p_i^{m-1}}{p_i^m-p_j^m}
=-m\,\frac{p_j^{m-1}\partial_{p_i}-p_i^{m-1}\partial_{p_j}}{p_i^m-p_j^m},
\]
}
in agreement with the expression for $\partial_{p_i}*\partial_{p_j}$.

To compare the diagonal products, we first show that $E$ is the unit of $F$.
For $i<j$ and $\zeta^m=1$ we expand $\alpha_H(X)v_H$ and obtain
{
\[
(X^i-\zeta X^j)(\partial_{p_i}-\zeta^{-1}\partial_{p_j})
=X^i\partial_{p_i}+X^j\partial_{p_j}-\zeta^{-1}X^i\partial_{p_j}-\zeta X^j\partial_{p_i},
\]
}
where $X^k=dp_k(X)$. The coefficient of $z^{m-1}$ in $z^m-1$ vanishes for $m\ge2$, so $\sum_\zeta\zeta=0$.
Since inversion permutes the roots, we also have $\sum_\zeta\zeta^{-1}=0$. It follows that
{
\[
\sum_{\zeta^m=1}(X^i-\zeta X^j)(\partial_{p_i}-\zeta^{-1}\partial_{p_j})=m\,(X^i\partial_{p_i}+X^j\partial_{p_j}).
\]
}
Each index $k$ belongs to $n-1$ pairs. Thus, summing over $i<j$, we find
{
\[
\sum_H\alpha_H(X)\,v_H=m\,(n-1)\sum_kX^k\partial_{p_k}=h\,X .
\]
}

Since $E=\frac1h\sum_kp_k\partial_{p_k}$ and each form $\alpha_H$ is linear, we have $\alpha_H(E)=\frac1h\,\alpha_H(p)$.
Substitution cancels the denominator in each term, and the preceding identity gives
{
\[
F(E,Y)=\sum_H\frac{\alpha_H(E)\,\alpha_H(Y)}{\alpha_H(p)}\,v_H=\frac1h\sum_H\alpha_H(Y)\,v_H=Y .
\]
}

Moreover, $E*Y=E^{-1}\circ E\circ Y=Y$, so $E$ is the unit of both products. To recover the
diagonal coefficients of $*$, we substitute $E=\frac1h\sum_kp_k\partial_{p_k}$ in the unit identities
for both products. At points where $p_j\ne0$ we obtain
{
\[
\partial_{p_j}*\partial_{p_j}=\frac1{p_j}\Big(h\,\partial_{p_j}-\sum_{k\neq j}p_k\,\partial_{p_k}*\partial_{p_j}\Big),
\qquad
F(\partial_{p_j},\partial_{p_j})=\frac1{p_j}\Big(h\,\partial_{p_j}-\sum_{k\neq j}p_k\,F(\partial_{p_k},\partial_{p_j})\Big),
\]
}
The calculation for $i\ne j$ shows that the sums on the right agree. Thus
$F(\partial_{p_j},\partial_{p_j})=\partial_{p_j}*\partial_{p_j}$ wherever $p_j\ne0$. This equality extends as an identity of rational functions
and, together with the calculation for $i\ne j$, proves that $*=F$.

{
\emph{Step 3.} We now construct the natural connection for arbitrary $\beta\neq0$.
We work on $M'$, where $a_n\neq0$, and use the coordinates $(a_1,\dots,a_N,a_n)$ with indices in $\{1,\dots,N,n\}$.
Set $g_l:=z\,\partial_{a_l}\lambda$. Thus $g_l=z^{N-l+1}$ for $l\le N$ and $g_n=1$.

We first introduce the polynomials $\tilde A_{ij}$. Set $R_{ij}:=\partial_{a_i}\lambda\,\partial_{a_j}\lambda-c^l_{ij}\,\partial_{a_l}\lambda$, so that $z^2R_{ij}=g_ig_j-c^l_{ij}\,z\,g_l$ is a polynomial
of degree $\le2N$. By the definition of the product, $R_{ij}\equiv0$ modulo $z^2\lambda'$. We can
therefore write $R_{ij}=z^2\lambda'\,g$ for a Laurent polynomial $g$. In the expression
$z^2R_{ij}=z^2\lambda'\cdot z^2g$, the factor $z^2\lambda'$ has nonzero constant term $-a_n$. Thus, if the least
power in $z^2g$ were negative, it would also occur in the product.
It follows that $\tilde A_{ij}:=z^2g$ is a polynomial of degree $\le2N-(N+1)=N-1$. It is symmetric
in $i,j$, and we have
\[
R_{ij}=\lambda'\,\tilde A_{ij},\qquad g_ig_j=z^2\lambda'\,\tilde A_{ij}+\rho_{ij},\qquad \rho_{ij}=c^l_{ij}\,z\,g_l\in\mathrm{span}(z,z^2,\dots,z^n) .
\]
More explicitly, let $q_{ij}$ and $r_{ij}$ denote the quotient and the remainder
in the ordinary polynomial division of $g_ig_j$ by $z^2\lambda'$. Then $\tilde A_{ij}=q_{ij}-r_{ij}(0)/a_n$. In particular, $\tilde A_{ij}$ is polynomial
in $a_1,\dots,a_N$ and $a_n^{-1}$ and is independent of $\beta$.

Observe that $\tilde A_{ij}/z\in\mathrm{span}(z^{N-2},\dots,1,z^{-1})\subset V$ and that the operator $z\frac d{dz}$ preserves this space. We
can therefore define a torsion free connection $\nabla'$ on $M'$ by requiring that
\begin{equation}\label{eq:chrL}
\sum_{l\in\{1,\dots,N,n\}}\Gamma^l_{ij}(a)\,\partial_{a_l}\lambda
=-\Big(z\frac{d}{dz}+\beta\Big)\frac{\tilde A_{ij}}{z}.
\end{equation}
The Christoffel symbols are affine in $\beta$. In particular, comparing the
coefficients of $z^{-1}$ gives $\Gamma^n_{ij}=(1-\beta)\tilde A_{ij}(0)$.

\emph{(a)} We verify conditions (S1) through (S3). As in the case of
$B_n$, we have $P=z^\beta/\beta$, $\partial_P=z^{1-\beta}\partial_z$ and $\partial_P\Lambda=z^{1-\beta}\lambda'$. Condition (S1) reads $\partial_{a_N}\lambda=1$. To check (S2), we
use $E\lambda=\lambda-\frac1Nz\lambda'$ from Step~1 and the identity $\nu_0P\,\partial_P\Lambda=\frac\beta N\,\frac{z^\beta}\beta\,z^{1-\beta}\lambda'=\frac1N\,z\lambda'$ for $\nu_0=\frac\beta N$. Adding the two identities
gives $E\Lambda+\nu_0P\,\partial_P\Lambda=\Lambda$, which is (S2) with $\kappa=0$. For (S3), we use \eqref{eq:defconn}.
Since $\lambda$ is linear in the coordinates $a$, we have $\partial_{a_i}\partial_{a_j}\lambda=0$. Moreover,
\[
B_{ij}=\frac{R_{ij}}{\partial_P\Lambda}=\tilde A_{ij}\,z^{\beta-1},\qquad
\partial_PB_{ij}=z^{1-\beta}\partial_z\big(\tilde A_{ij}z^{\beta-1}\big)=\Big(z\frac d{dz}+\beta\Big)\frac{\tilde A_{ij}}z ,
\]
so that \eqref{eq:defconn} takes the form $\Gamma^l_{ij}\,\partial_{a_l}\lambda=-\partial_PB_{ij}$, which follows from
\eqref{eq:chrL}. Thus $\nabla'$ is the connection defined by $\Lambda$ in Theorem~\ref{thm:flat}.

\emph{(b)} To prove flatness, consider the local branches of the inverse
of $\Lambda$. They are $P_j=z_j(a,\lambda)^\beta/\beta$, where $z_j$ are the $n$ nonzero roots of
$z(\lambda(z)-\lambda)=z^n+a_1z^{n-1}+\dots+(a_N-\lambda)z+a_n$. Away from the discriminant these roots are local coordinates.
Since $z\mapsto z^\beta/\beta$ is locally invertible for $z\neq0$, the differentials $dP_j$ are
independent. Part (a) and Theorem~\ref{thm:flat} now imply that
$\nabla'$ is flat, $\nabla'e=0$ and $\nabla'\nabla'E=0$.

\emph{(c)} To find flat coordinates, we begin with the differentials
of the twisted residues in \eqref{eq:TL}.

\begin{lemma}\label{lem:phiL}
Consider the functionals $\varphi_m(f)=\res_{z=\infty}\lambda^{s_m-1}f\,z^{\beta-1}dz$ for $m=1,\dots,N$ and $\varphi_n(f)=\res_{z=0}\lambda^{\beta-1}f\,z^{\beta-1}dz$ on $V$. For $m\le N$ we have
$\varphi_m(z^{N-i})=0$ if $i>m$, together with $\varphi_m(z^{N-m})=1$ and $\varphi_m(z^{-1})=0$. The functional $\varphi_n$ vanishes on
polynomials and satisfies $\varphi_n(z^{-1})=a_n^{\beta-1}\ne0$. Moreover, $\partial u_m/\partial a_l=\varphi_m(\partial_{a_l}\lambda)$ for $l,m\in\{1,\dots,N,n\}$, and the
functions in \eqref{eq:TL} form a system of local coordinates on $M'$.
In particular, $\varphi_m(1)=\delta_{mN}$ for $m\le N$ and $\varphi_n(1)=0$.
\end{lemma}

\begin{proof}[Proof of Lemma~\ref{lem:phiL}]
As in the proof of Lemma~\ref{lem:phi}, we use the expansion
$\lambda^{s_m-1}z^{N-i}z^{\beta-1}=z^{m-1-i}\big(1+a_1z^{-1}+\dots+a_nz^{-n}\big)^{s_m-1}$ at infinity. Reading off the coefficient of $z^{-1}$ gives the first
two values. For $f=z^{-1}$ the largest exponent is $m-1-n<-1$, so this coefficient
vanishes. At zero we have $\lambda=a_nz^{-1}(1+O(z))$ and therefore $\lambda^{\beta-1}f\,z^{\beta-1}=a_n^{\beta-1}f\,(1+O(z))$, which gives the
stated values of $\varphi_n$ on polynomials and on $z^{-1}$.
The Jacobian of \eqref{eq:TL} with respect to $(a_1,\dots,a_N,a_n)$ is therefore
triangular, with last diagonal entry $\partial u_n/\partial a_n=\varphi_n(z^{-1})=a_n^{\beta-1}\ne0$. Thus the functions
\eqref{eq:TL} form a system of local coordinates on $M'$. A change
of branch of $a_n^\beta$ multiplies $u_n$ by a nonzero constant and fixes
the remaining coordinates. Consequently the connections defined
in these coordinates agree on overlaps and give a connection on $M'$. In computing the Jacobian we
used $\partial u_m/\partial a_l=\varphi_m(\partial_{a_l}\lambda)$, which follows by differentiating \eqref{eq:TL} under the
residue as in Lemma~\ref{lem:phi}.
\end{proof}

We show that the twisted residues \eqref{eq:TL} at the two poles
of $\lambda$ are flat coordinates of $\nabla'$. We may treat both poles together
by writing $u=\frac1s\res_q\lambda^s\omega$, where $\omega=z^{\beta-1}dz$ and either $(q,s)=(\infty,s_m)$ with $m\le N$, or $(q,s)=(0,\beta)$. Denote by
$\varphi=\res_q\lambda^{s-1}(\cdot)\,\omega$ the corresponding functional in Lemma~\ref{lem:phiL}. Thus $\partial_{a_l}u=\varphi(\partial_{a_l}\lambda)$.
The forms that occur in the following calculation have integer
exponents at $q$. Indeed, at $q=\infty$ we have $\lambda^{s_m-1}z^{\beta-1}=z^{m-N-1}(1+O(z^{-1}))$. At $q=0$ the expansion
$\lambda=a_nz^{-1}(1+O(z))$ gives $\lambda^{\beta-1}z^{\beta-1}=a_n^{\beta-1}(1+O(z))$. We can therefore integrate by parts to obtain
\[
\begin{aligned}
\Gamma^l_{ij}\,\partial_{a_l}u&=\varphi\big(\Gamma^l_{ij}\partial_{a_l}\lambda\big)\\
&=-\res_q\lambda^{s-1}\,\partial_z\big(\tilde A_{ij}z^{\beta-1}\big)\,dz\\
&=(s-1)\,\res_q\lambda^{s-2}\lambda'\,\tilde A_{ij}\,\omega=(s-1)\,\res_q\lambda^{s-2}R_{ij}\,\omega .
\end{aligned}
\]
For the term in $R_{ij}$ involving the structure constants, we obtain
$(s-1)\,c^l_{ij}\res_q\lambda^{s-2}\partial_{a_l}\lambda\,\omega=c^l_{ij}\,\partial_{a_l}\varphi(1)=0$. Indeed, $\varphi(1)$ is constant by Lemma~\ref{lem:phiL}. It follows that
\[
\Gamma^l_{ij}\,\partial_{a_l}u=(s-1)\,\res_q\lambda^{s-2}\,\partial_{a_i}\lambda\,\partial_{a_j}\lambda\,\omega=\partial_{a_i}\big(\res_q\lambda^{s-1}\partial_{a_j}\lambda\,\omega\big)=\partial_{a_i}\partial_{a_j}u ,
\]
which proves that $u$ is a flat coordinate. By Lemma~\ref{lem:phiL},
$(u_1,\dots,u_N,u_n)$ are coordinates on $M'$. Moreover, $\partial_{a_N}u_m=\varphi_m(1)$ equals $\delta_{mN}$ for $m\le N$ and $0$
for $m=n$. Thus $e=\partial_{a_N}=\partial_{u_N}$. The residues are homogeneous, with $\deg u_k=\frac kN$ for
$k\le N$ and $\deg u_n=\beta\deg a_n=\frac{n\beta}N$. It follows that
\[
E=\sum_{k\le N}\frac kN\,u_k\partial_{u_k}+\frac{n\beta}N\,u_n\partial_{u_n},\qquad [e,E]=e .
\]

\emph{(d)} We now verify the conditions for a flat F-manifold.
In the coordinates of part (c) we have
\[
N_0=\nu_0-\nabla'E=\mathrm{diag}\Big(\frac{\beta-1}N,\dots,\frac{\beta-N}N,\ \frac\beta N-\frac{n\beta}N\Big)
=\mathrm{diag}\Big(\frac{\beta-1}N,\dots,\frac{\beta-N}N,\ -\beta\Big),
\]
which is invertible for $\beta\notin\{0,1,\dots,N\}$ and, in particular, for $\beta=\frac1m$.
Theorem~\ref{thm:flat}(ii) then gives the symmetry of $\nabla'_lc^m_{hk}$. In the
coefficient coordinates the structure constants are independent of $\beta$
and the Christoffel symbols are affine in $\beta$. Thus the symmetry holds for
every $\beta\neq0$. Together with Step~1, this proves that $(\nabla',\circ,e)$ is a flat F-manifold
on $M'$ with linear Euler vector field $E$.

\emph{(e)} We now consider polynomiality at $\beta=\frac1m$. By Lemmas~\ref{lem:phiL}
and~\ref{lem:phi}, for $\beta=\frac1m$ the coordinates $u_k$, $k\le N$, are polynomial in $a_1,\dots,a_N$ and the inverse change of coordinates is also polynomial. Moreover, $u_n=m\,a_n^{1/m}=m\big((-1)^nt^m\big)^{1/m}=\epsilon\,m\,t$. Thus $(u_1,\dots,u_N,u_n)$
are polynomial coordinates on $\mathbb C^n$. By Step~1 the structure constants are
polynomial in these coordinates. The symmetry of $\partial_lc^m_{hk}$ allows us to integrate
twice to obtain a polynomial vector potential. The structure therefore
extends from $M'$ to all of $\mathbb C^n$.

\emph{Step 4.} We now identify the dual connection and the bi-flat
structure for $\beta=\frac1m$ and $\nu_0=\frac1h$. By Proposition~\ref{prop:dualflat}, wherever
$E\circ$ is invertible the connection $D_XY:=\nabla'_XY-\nabla'_{X*Y}E$ is flat. Moreover, $(D,*,E)$ is a flat F-manifold with $DE=0$, and $(\nabla',\circ,e,D,*,E)$ is a bi-flat F-manifold. Applying Step~3(a)
and Theorem~\ref{thm:char}(iv) at $\lambda=0$, we find that $w_j=P_j(\cdot,0)=z_j^\beta/\beta$, where $z_j$ are
the nonzero roots of $z\lambda(z)$, are flat coordinates for $\Gamma(0,\nu_0)=D+\nu_0\,*$. Since $\beta=\frac1m$,
these functions are $w_j=m\,p_j$ and are therefore linear in $p$. It follows
that the torsion free connection $D+\nu_0\,*$ coincides with the trivial connection
$\tilde\nabla=\partial$ in the coordinates $p$. In particular,
\[
D_XY=\partial_XY-\frac1h\,X*Y=\partial_XY-\frac1h\sum_H\frac{\alpha_H(X)\,\alpha_H(Y)}{\alpha_H(p)}\,v_H=\nabla^*_XY .
\]
This connection is invariant under $G(m,m,n)$ on $V^{\mathrm{reg}}$ and is of the type studied
by Dunkl and Kohno. The hyperplanes form a single orbit and occur in the
preceding formula with the common coefficient $\frac1h$. Thus their common weight is
$\frac2h$ when we express the connection in terms of the projections $\pi_H$. In particular, $(\nabla^*,*,E)$ is a flat F-manifold
and $\nabla^*E=0$. We can also verify the latter identity directly. Using $\partial_XE=\frac1hX$, $\alpha_H(E)=\frac1h\alpha_H(p)$
and $\sum_H\alpha_H\otimes v_H=h\,\Id$, we obtain
\[
\nabla^*_XE=\frac1hX-\frac1h\sum_H\frac{\alpha_H(X)\,\alpha_H(E)}{\alpha_H(p)}\,v_H=\frac1hX-\frac1{h^2}\sum_H\alpha_H(X)\,v_H=0 .
\]
Thus the flatness of $\nabla^*$ follows from the construction itself, independently
of any general flatness theorem for connections associated with reflection
arrangements.

Finally, as in Step~4 of the proof of Theorem~\ref{thm:Bn},
\eqref{eq:dual} and $\nabla'e=0$ give $\nabla^*_{X\circ Y}e=-\nabla'_{X*Y}E=\nabla^*_XY-\nabla'_XY$. We thus obtain $\nabla'_XY=\nabla^*_XY-\nabla^*_{X\circ Y}e=\nabla_XY$, which identifies
$\nabla'$ with the connection in (iv).
}

\emph{Step 5.} On $M^{\mathrm{reg}}\cap\{t\ne0\}$ we have constructed two flat F-manifolds $(\nabla,\circ,e)$
and $(\nabla^*,*,E)$. By Steps~1 through~4, their products satisfy $X*Y=E^{-1}\circ X\circ Y$ with $E\circ$
invertible, and their unit and Euler fields satisfy $[e,E]=e$, $\mathrm{Lie}_E\circ=\circ$, $\nabla\nabla E=0$ and
$\nabla^*E=0$. Step~4 shows that the two structures form a bi-flat pair related
by \eqref{eq:dual}. The compatibility condition follows from
Lemma~\ref{lem:AHE}.

By Step~3(e) the natural data extend polynomially to $M$. The formulas
in (ii) and (iii) are holomorphic on $M^{\mathrm{reg}}$, including the points where $t=0$.
To compare the structures at these points, define $L=E\circ$ on $M$ and let
$A$ denote multiplication by $e$ with respect to the arrangement product
$*$ on $M^{\mathrm{reg}}$. On $M^{\mathrm{reg}}\cap\{t\ne0\}$ we have $A=L^{-1}$ and therefore $LA=\Id$. By holomorphicity the
latter identity holds throughout $M^{\mathrm{reg}}$, so $L$ is invertible there. The
relation $X*Y=L^{-1}(X\circ Y)$ and the compatibility identities now extend from the dense
open set. Thus we obtain a bi-flat pair on all of $M^{\mathrm{reg}}$.

We next identify the superpotential. Choose local branches on $M^{\mathrm{reg}}\cap\{a_n\ne0\}$.
By Step~4 the functions $w_j=m\,p_j$ are flat coordinates for $\Gamma(0,\nu_0)=\nabla^*+\nu_0c^*$. They are
homogeneous of degree $\deg p_j=\frac1h=\nu_0$, which gives $\kappa=0$. To compute their derivatives
along $e$, recall that $z\lambda(z)=z^n+\dots+a_Nz+a_n$. Differentiating $(z\lambda)(z_j)=0$ with respect to $a_N$ gives
$(z\lambda)'(z_j)\,\partial_{a_N}z_j+z_j=0$. At a root we have $\lambda(z_j)=0$, so $(z\lambda)'(z_j)=z_j\lambda'(z_j)$. It follows that
{
\[
\partial_ew_j=z_j^{\beta-1}\partial_{a_N}z_j=-\frac{z_j^{\beta-1}}{\lambda'(z_j)}\neq0 .
\]
}
Since the flow $\phi_{-\lambda}$ of $e=\partial_{a_N}$ acts by $a_N\mapsto a_N-\lambda$, we have $\lambda(z,\phi_{-\lambda}(a))=\lambda(z,a)-\lambda$. Thus the zeros of
$\lambda(\cdot\,,\phi_{-\lambda}(a))$ are precisely the solutions of $\lambda(z,a)=\lambda$, and we obtain $w_j(\phi_{-\lambda}(a))=z_j(a,\lambda)^\beta/\beta=P_j(a,\lambda)$. Inverting
this relation with respect to $\lambda$ gives $\lambda(z)$ as a function of $P=z^\beta/\beta=m\,p$.
According to Definition~\ref{def:sp}, $\Lambda=\lambda(p^m)$ is therefore the
superpotential associated with $w_j$, of weight $\frac1h$ and constant $\kappa=0$.

Suppose, on the same open set, that $z^2\lambda'$ has simple roots. These are the
$n$ nonzero simple critical points on the curve with coordinate $z$, and
their critical values are finite. Choose one lift of each critical point
under $z=p^m$. Each point has $m$ lifts, and $P=mp$ is a local coordinate at each
of them. The differential of each critical value is the evaluation
map $X\mapsto(X\lambda)(z_*)$ at the corresponding point. Since $\mathbb C[z,z^{-1}]/(z^2\lambda')$ is the direct product of the
$n$ evaluation algebras, these differentials are independent characters
of the algebra. It follows that the critical values are canonical
coordinates. Moreover, they are generically distinct. Indeed, for
$\lambda=z^N+a_n/z$ with $a_n\ne0$, the critical values are $(N+1)z_*^N$ at the distinct solutions of
$z_*^{N+1}=a_n/N$.

We finally determine the weights for the connections and products
of the stated arrangement class. Since the mirrors form a single orbit, there is one weight
for the connection and one for the product. The identities $\nabla^*E=0$ and $E*X=X$
determine them. Once we fix the unit $e$, the reconstruction following
Lemma~\ref{lem:AHE} determines the natural structure. In the case
$m=2$ we obtain $\lambda(p^2)=p^{2n-2}+a_1p^{2n-4}+\dots+a_{n-1}+(-1)^nt^2p^{-2}$. With $\omega\propto dp$ this recovers the superpotential that
Dubrovin associated with the Frobenius manifold $D_n$ in \cite{Dub}.
\end{proof}

The argument in parts (a) through (d) of Step~3 and Theorem~\ref{thm:char}(iv)
also apply on $M'$ in rank two and for weights other than $\beta=\frac1m$. We thus
obtain the following family on the space of coefficients of Laurent polynomials.

\begin{corollary}\label{thm:Laurent}
For every $n\ge2$ and $\beta\ne0$, the product of Step~1, the unit $e=\partial_{u_N}$, the Euler
field $E=\sum_k\deg(a_k)a_k\partial_{a_k}$ and the connection with local flat coordinates \eqref{eq:TL}
define a flat F-manifold on $M'$. For $\nu_0=\beta/N$ and $\kappa=0$, the function $\lambda$
written in the local coordinate $P=z^\beta/\beta$ satisfies conditions (S1) through
(S3) wherever $z\ne0$ and $\lambda'(z)\ne0$. On the open set where $z\lambda(z)$ has distinct
roots, choose a local ordering $z_1,\dots,z_n$ of the roots and local branches of
the logarithm. Then $z_j^\beta/\beta$ are local flat coordinates for the connection
$\Gamma(0,\nu_0)$ in \eqref{eq:op}. In the coordinates $(a_1,\dots,a_N,a_n)$ the Christoffel symbols
of $\nabla$ have the form \eqref{eq:chrL}, where $\tilde A_{ij}$ is the polynomial
introduced in Step~3 of the proof of Theorem~\ref{thm:Gmmn}.
These symbols are affine in $\beta$ and polynomial in $a_1,\dots,a_N,a_n^{-1}$. In particular,
$\Gamma^n_{ij}=(1-\beta)\tilde A_{ij}(0)$.
\end{corollary}

The assertions about the Christoffel symbols follow from Step~3 of
the proof of Theorem~\ref{thm:Gmmn}. For instance, since $g_n=1$, we have
$\tilde A_{nn}=-1/a_n$ and therefore $\Gamma^n_{nn}=(\beta-1)/a_n$.

At $\beta\in\{1,\dots,N\}$ we use the analytic continuation of the normalized residues in
\eqref{eq:TL}. The Jacobian of the flat coordinates is still invertible.
{The argument of Step~3(d) also gives the symmetry of $\nabla_l c^m_{hk}$ in $l,h,k$ at these
values. Indeed, in the coordinates given by the coefficients, the product
is independent of $\beta$ and the Christoffel symbols are affine in $\beta$.
In the flat coordinates $u$ this is precisely the symmetry of $\partial_l c^m_{hk}$.}

When $m=2$, the extension by one coordinate described in Remark~\ref{rem:PS}
can be compared with the solution of the open WDVV equations of type
$D_n$ in Example~5.2 of \cite{PS}. The latter is obtained from $\lambda=z^{n-1}+a_1z^{n-2}+\dots+a_{n-1}-\frac{a_n^2}{4z}$ with
$\omega=dz/\sqrt{2z}$, or equivalently $\omega=dx$ after the change of variable $x^2=2z$. Up to rescaling, this
is the extension defined by the superpotential of Theorem~\ref{thm:Gmmn}.

The family in the corollary is defined on $M'$. We now study its extension
across $t=0$ after pulling it back to the quotient of $G(m,m,n)$.

\begin{remark}\label{rem:Gmmnbeta}
Let $m\ge2$ and $n\ge3$. On $t\ne0$ we pull back the structures of
Corollary~\ref{thm:Laurent} by $a_n=(-1)^nt^m$ and use the polynomial coordinates $(u_1,\dots,u_N,t)$.
The flat function $u_n$ is locally a nonzero constant multiple of $t^{m\beta}$.
Its flatness equation thus gives
\[
\Gamma^t_{tt}=\frac{\partial_t^2u_n}{\partial_tu_n}=\frac{m\beta-1}{t}.
\]
Thus the natural connection extends across $t=0$ if and only if $m\beta=1$.
For the dual connection we use flat functions of $\nabla^*+\nu_0*$ proportional to $p_j^{m\beta}$. The same calculation gives a pole along $p_j=0$
unless $m\beta=1$. A generic point of $p_j=0$ lies in $V^{\mathrm{reg}}$, where the arrangement
product $*$ is holomorphic. The pole therefore persists in the dual
connection. At $\beta=1/m$ both extensions are those of Theorem~\ref{thm:Gmmn}.

For $\beta\neq\frac1m$, Step~3 of the proof of Theorem~\ref{thm:Gmmn} and
Proposition~\ref{prop:dualflat} still apply on the part of $M'$ where
$z\lambda(z)$ has distinct roots. There we compute the dual connection using
the flat coordinates $z_j^\beta/\beta$ as in Step~4 of the proof of Theorem~\ref{thm:Bn}.
The resulting bi-flat F-manifolds have the dual product and Euler field
$E=\frac1h\sum_jp_j\partial_{p_j}$ of $G(m,m,n)$, together with the connection of Dunkl and Kohno for the
arrangement of $G(m,1,n)$. With the normalization of Section~\ref{sec:Gm1n},
the weights of this connection are $\frac\beta{n-1}$ on $p_i=\zeta p_j$ and $\frac1m-\beta$ on $p_j=0$.
In the family considered there the dual product has additional
terms on $p_j=0$, and the Euler field is $\frac1{mn}\sum_jp_j\partial_{p_j}$. Thus, within the Laurent
family, $\beta=1/m$ is the unique parameter for which the natural connection
is defined on the entire quotient and the dual data are defined on
the entire regular quotient. In rank two the group is $I_2(m)$. Its
reflecting hyperplanes form one orbit for odd $m$ and two orbits for even $m$.
\end{remark}
}

\subsubsection{\texorpdfstring{The vector potentials for $n=2,3,4$}{The vector potentials for n=2,3,4}}\label{sec:Laurentpot}

{
We give the vector potentials of Corollary~\ref{thm:Laurent} for generic
$\beta$ in ranks two, three and four. In \eqref{eq:TL} the pole at infinity
gives the flat coordinates $u_1,\dots,u_N$ and the pole at zero gives $u_n=a_n^\beta/\beta$. Thus, in
the formulas below, $a_n=(\beta u_n)^{1/\beta}$. The structure constants are $c^p_{kl}=\partial_k\partial_lF^p$, where $F^p$ are
the components of the vector potential. In each case the unit is $\partial_{u_N}$
and the Euler field is $E=\sum_{k\le N}\frac kNu_k\partial_{u_k}+\frac{n\beta}Nu_n\partial_{u_n}$. The only terms that need not be polynomial
in the flat coordinates are linear in $a_n$.

When $n=2$ we have $N=1$ and the unit is $\partial_{u_1}$. The components of the vector
potential are
{
\[
F^1=\tfrac12u_1^2+\frac{a_n}{1-\beta},\qquad F^2=u_1u_2 .
\]
}

When $n=3$ we have $N=2$ and the unit is $\partial_{u_2}$. The three components are
{
\begin{align*}
F^1&=u_1u_2+\tfrac1{12}(2\beta-1)\,u_1^3+\frac{2a_n}{1-\beta},\\
F^2&=\tfrac12u_2^2+u_1\,a_n+\tfrac1{48}\beta(1-\beta)\,u_1^4,\\
F^3&=u_2u_3+\tfrac\beta4\,u_1^2u_3 .
\end{align*}
}

Finally, when $n=4$ we have $N=3$ and the unit is $\partial_{u_3}$. In this case we obtain
{
\begin{align*}
F^1&=u_1u_3+\tfrac12u_2^2+\tfrac16(2\beta-1)\,u_1^2u_2+\tfrac1{216}(2\beta-1)(3\beta-1)\,u_1^4+\frac{3a_n}{1-\beta},\\
F^2&=u_2u_3+u_1\,a_n-\tfrac16(1-\beta)\,u_1u_2^2+\tfrac1{54}(1-\beta)\,u_1^3u_2
+\tfrac1{1080}(2\beta-1)(1-\beta)(1+\beta)\,u_1^5,\\
F^3&=\tfrac12u_3^2+u_2\,a_n+\tfrac16u_1^2\,a_n-\tfrac\beta{18}u_2^3+\tfrac1{18}\beta(1-\beta)\,u_1^2u_2^2
-\tfrac1{216}\beta(1-\beta)(1-2\beta)\,u_1^4u_2\\
&\quad+\tfrac1{4860}\beta(1-\beta)(1-2\beta+2\beta^2)\,u_1^6,\\
F^4&=u_3u_4+\tfrac\beta3\,u_1u_2u_4+\tfrac{\beta^2}{27}\,u_1^3u_4 .
\end{align*}
}

To verify the formulas, we compute the products of derivatives of $\lambda$
using the relation $\partial_k\lambda\,\partial_l\lambda\equiv c^p_{kl}\,\partial_p\lambda$ modulo $z^2\lambda'$. Here $\partial_k=\partial/\partial u_k$ with $z$ fixed, and the
representative on the right belongs to $V$. The coefficients so obtained
agree with $\partial_k\partial_lF^p$, as we see by using $\partial_{u_n}a_n=a_n^{1-\beta}$. In each of these ranks the last
component is $F^n=a_N\,u_n$, where $a_N$ is expressed in flat coordinates.

The residue coordinates are regular at $\beta=1$. At this value we can
remove the apparent pole of $F^1$ by adding an affine term. Indeed,
$N\,a_n/(1-\beta)$ has the same Hessian as $N(a_n-u_n)/(1-\beta)$, and the latter expression tends to
$Nu_n(\log u_n-1)$ on a local branch of the logarithm.

For $\beta=\frac1m$ we have $a_n=(u_n/m)^m$, and the potentials are therefore polynomial.
When $n=3,4$ they give the products of Theorem~\ref{thm:Gmmn}. For
$m=2$ we recover $D_2=A_1\times A_1$, $D_3=A_3$ and $D_4$ in ranks two, three and four,
respectively. In rank two the idempotent coordinates are $x_\pm=u_1\pm u_2$, in
which the product Frobenius metric takes the form $dx_+^2+dx_-^2$.

For $n=3$ and $\beta=\frac12$ we have $a_n=u_3^2/4$. The three components in this case are
$F^1=u_1u_2+u_3^2$, $F^2=\tfrac12u_2^2+\tfrac14u_1u_3^2+\tfrac1{192}u_1^4$ and $F^3=u_2u_3+\tfrac18u_1^2u_3$.

\begin{proof}[Proof of Theorem~\ref{thm:intro}]
A nontrivial irreducible well generated group has either one or two orbits
of reflecting hyperplanes. In the first case the assertion follows from
Corollary~\ref{cor:oneorbit}. In the second case, Proposition~\ref{prop:Gm1n}
gives the construction for $G(m,1,n)$ and Proposition~\ref{prop:F4} gives the
family for $F_4$. The results cited in Table~\ref{tab:coverage} cover the
even dihedral groups and the remaining exceptional groups of ranks two
and three. In each case the natural flat structure is polynomial on
the quotient and the dual data are defined on the regular quotient.
\end{proof}
}


\begin{thebibliography}{99}

\bibitem{Alc} A.~Alcolado, \emph{Extended Frobenius manifolds and the open WDVV equations}, Ph.D.\
thesis, McGill University, 2017.

\bibitem{Alm} G.~F.~Almeida, \emph{Open Hurwitz flat F-manifolds}, arXiv:2503.09258 (2025).

\bibitem{Ar} D.~Armstrong, \emph{Generalized noncrossing partitions and combinatorics of
  Coxeter groups}, Mem.\ Amer.\ Math.\ Soc.\ \textbf{202} (2009), no.~949.

\bibitem{AW} C.E.Arreche, N.F.Williams, \emph{Normal reflection subgroups of complex
reflection groups}, J.\ Inst.\ Math.\ Jussieu 22 (2023), no.~2, 879--917.

\bibitem{ABLR} A.~Arsie, A.~Buryak, P.~Lorenzoni, P.~Rossi, \emph{Riemannian F-manifolds, bi-flat
F-manifolds, and flat pencils of metrics}, Int.\ Math.\ Res.\ Not.\ IMRN \textbf{2022}, no.~21,
16730--16778.

\bibitem{ALDE} A.~Arsie, P.~Lorenzoni, \emph{From Darboux--Egorov system to bi-flat
$F$-manifolds}, J.\ Geom.\ Phys.\ \textbf{70} (2013), 98--116.

\bibitem{ALjmp}
A. Arsie and  P. Lorenzoni, \emph{Purely non-local Hamiltonian formalism, Kohno connections and $\vee$-systems},  J. Math. Phys. 55 (2014).

\bibitem{ALrefl} A.Arsie, P.Lorenzoni, \emph{Complex reflection groups, logarithmic connections
and bi-flat F-manifolds}, Lett.\ Math.\ Phys.\ 107 (2017), no.~10, 1919--1961.

\bibitem{Bn} A.Arsie, P.Lorenzoni, \emph{Flat F-manifolds, Miura invariants, and integrable
systems of conservation laws}, J.\ Integrable Syst.\ 3 (2018), no.~1, xyy004.
arXiv:1709.10300.

\bibitem{ALsurvey} A.Arsie, P.Lorenzoni, \emph{Bi-flat F-manifolds: A survey}, in
R.Donagi and T.Shaska (eds.), \emph{Integrable Systems and Algebraic Geometry},
Vol.~1, London Mathematical Society Lecture Note Series 458, Cambridge University
Press, 2020, Chapter~8, 167--191. 

\bibitem{AL} A.Arsie, P.Lorenzoni, \emph{Bi-flat F-structures as differential bicomplexes
and Gauss Manin connections}, Bull.\ Lond.\ Math.\ Soc.\ 57 (2025), no.~3, 786--808.

\bibitem{Be1} D.~Bessis, \emph{The dual braid monoid}, Ann.\ Sci.\ \'Ecole Norm.\ Sup.\ (4) \textbf{36}
  (2003), no.~5, 647--683.
  
\bibitem{Be2} D.~Bessis, \emph{Finite complex reflection arrangements are $K(\pi,1)$},
  Ann.\ of Math.\ (2) \textbf{181} (2015), 809--904.
  
\bibitem{BGU} {G.~Borot, A.~Giacchetto, and G.~Umer, \emph{Symmetries of F-Cohomological Field Theories and F-Topological Recursion}, Commun.\ Math.\ Phys.\ \textbf{406} (2025), no.~10, 248.}

\bibitem{Br} M.~Brou\'e, \emph{Introduction to Complex Reflection Groups and Their Braid
  Groups}, Lecture Notes in Math.\ \textbf{1988}, Springer, 2010.

\bibitem{Ch} C.~Chevalley, \emph{Invariants of finite groups generated by reflections},
  Amer.\ J.\ Math.\ \textbf{77} (1955), 778--782.  

\bibitem{CHL} W.~Couwenberg, G.~Heckman, E.~Looijenga, \emph{Geometric structures on the complement
of a projective arrangement}, Publ.\ Math.\ IH\'ES 101 (2005), 69--161.

\bibitem{Dub} B.~Dubrovin, \emph{Geometry of 2D topological field theories}, in: Integrable
Systems and Quantum Groups, Lecture Notes in Math.\ 1620, Springer, 1996.

\bibitem{DCoxeter} B. Dubrovin, {\em Differential Geometry of the Space of Orbits of a Coxeter Group},  Surveys in Differential Geometry , Vol. IV (1999), p. 181 -- 212.

\bibitem{Dalm} B.~Dubrovin, \emph{On almost duality for Frobenius manifolds}, in:
Geometry, Topology, and Mathematical Physics, Amer.\ Math.\ Soc.\ Transl.\ Ser.~2, 212 (2004).

\bibitem{FV} M. V. Feigin and A. P. Veselov, \emph{$\vee$-systems, holonomy Lie algebras and logarithmic vector fields},  IMRN, Volume 2018, Issue 7, 11 April 2018, p. 2070--2098.

\bibitem{Hu} J.\,E.~Humphreys, \emph{Reflection Groups and Coxeter Groups}, Cambridge
  Stud.\ Adv.\ Math.\ \textbf{29}, Cambridge Univ.\ Press, 1990.

\bibitem{IY} Y.~Ilyashenko and S.~Yakovenko, \emph{Lectures on Analytic Differential
Equations}, Graduate Studies in Mathematics, vol.~86, American Mathematical Society,
Providence, RI, 2008.

\bibitem{KMS1} M. Kato, T. Mano and J. Sekiguchi, \emph{Flat Structures without Potentials},  Rev. Roumaine Math. Pures Appl. {\bf 60} (2015), no. 4, 481--505

\bibitem{KMS2} M. Kato, T. Mano and J. Sekiguchi \emph{Flat structure on the space of isomonodromic
deformations}, SIGMA \textbf{16} (2020), 110. {\tt arXiv:1511.01608}.

\bibitem{KM} {Y.~Konishi and S.~Minabe, \emph{Almost duality for Saito structure and complex reflection groups II: the case of Coxeter and Shephard groups}, Pure Appl.\ Math.\ Q.\ \textbf{16} (2020), no.~3, 721--754.}

\bibitem{KMS} Y. Konishi, S.  Minabe, Y. Shiraishi, \emph{Almost duality for Saito structure and
complex reflection groups}, J.\ Integrable Syst.\ 3 (2018), no.~1, xyy003.

\bibitem{LT} G.\,I.~Lehrer, D.\,E.~Taylor, \emph{Unitary Reflection Groups},
  Cambridge Univ.\ Press, 2009.

\bibitem{OS} P.~Orlik, L.~Solomon, \emph{Unitary reflection groups and cohomology},
  Invent.\ Math.\ \textbf{59} (1980), 77--94.
  
\bibitem{OT} P.~Orlik, H.~Terao, \emph{Arrangements of Hyperplanes}, Grundlehren Math.\
  Wiss.\ \textbf{300}, Springer, 1992.  

\bibitem{PS} A.Proserpio, I.A.B.  Strachan, \emph{Dubrovin duality for open Hurwitz flat
F-manifolds}, arXiv:2512.08795.

\bibitem{SYS} K. Saito, T. Yano and J. Sekiguchi, \emph{On a certain generator system of the ring of
invariants of a finite reflection group}, Comm. in Algebra 8(4) (1980) 373 - 408.

\bibitem{Sa} K. Saito, \emph{On a linear structure of a quotient variety by a finite reflexion group},
RIMS Kyoto preprint 288, 1979.

\bibitem{ST} G.\,C.~Shephard, J.\,A.~Todd, \emph{Finite unitary reflection groups},
  Canad.\ J.\ Math.\ \textbf{6} (1954), 274--304.

\bibitem{Ve} A. P. Veselov, \emph{Deformations of the root systems and new solutions to generalised WDVV equations}, Physics Letters A, {\bf 261}, (1999), 297-302

\bibitem{Ve2} A. P. Veselov, \emph{On geometry of a special class of solutions to
generalized WDVV equations}, Integrability: the Seiberg-Witten and Whitham equations (Edinburgh, 1998), 125--135, Gordon and Breach. 

\end{thebibliography}
\end{document}